\documentclass{amsart}
\usepackage{amssymb,latexsym,amscd,hyperref}
\usepackage{tikz-cd}
\usepackage{pb-diagram}
\usepackage{MnSymbol}
\usepackage{ dsfont }
\theoremstyle{plain}
\newtheorem{theorem}{Theorem}
\newtheorem{corollary}[theorem]{Corollary}
\newtheorem{lemma}[theorem]{Lemma}
\newtheorem{example}[theorem]{Example}
\newtheorem{proposition}[theorem]{Proposition}
\newtheorem{algorithm}[theorem]{Algorithm}

\theoremstyle{definition}
\newtheorem{definition}[theorem]{Definition}

\makeatletter
\def\th@remark{%
	\thm@headfont{\bfseries}%
	\normalfont % body font
	\thm@preskip\topsep \divide\thm@preskip\tw@
	\thm@postskip\thm@preskip
}
\makeatother

\theoremstyle{remark}

\newtheorem{remark}[theorem]{Remark}

\numberwithin{theorem}{section} 
\makeatletter
\newlength{\currentindent}
\newcommand{\@minipagerestore}{\setlength{\parindent}{\currentindent}}
\makeatother

\newcommand {\C}{{\mathbb{C}}}
\newcommand {\R}{{\mathbb{R}}}
\newcommand {\Z}{{\mathbb{Z}}}
\newcommand {\N}{{\mathbb{N}}}
\newcommand {\F}{{\mathbb{F}}}
\newcommand {\PP}{\mathbb{P}}

\newcommand{\np}{ \vert \vert  {\mathfrak{p}} \vert \vert }

\newcommand\Fqtl{\mathbb{F}_q((t))}
\newcommand\Fqtg{\mathbb{F}_q(t)}

\newcommand{\Ker}      {\mathop{\rm {Ker}}}
\newcommand{\im}      {\mathop{\rm {Im}}}
\newcommand{\ASC}      {\mathop{\rm {ASC}}}

\newcommand {\OO}{{\mathcal{O}}}

\newcommand {\idp}{{\mathfrak{p}}}

\newcommand{\Gal}      {\mathop{\rm {Gal}}}

\newcommand{\Aut}      {\mathop{\rm {Aut}}}

\DeclareMathOperator{\re}{Re}

\newcommand{\Pic}      {\mathop{\rm {Pic}}}
\newcommand{\Div}      {\mathop{\rm {Div}}}
\newcommand{\Hom}      {\mathop{\rm {Hom}}}
\newcommand{\Inj}      {\mathop{\rm {Inj}}}
\newcommand{\Surj}      {\mathop{\rm {Surj}}}

\usepackage{indentfirst}

\begin{document}

\title[Asymptotics of elementary-abelian extensions]{On the asymptotics of elementary-abelian extensions of local and global function fields}

\author{Nicolas Potthast}
\address{University Paderborn, Department of Mathematics,
Warburger Str.~100, 33098 Paderborn, Germany} 

\email{nicolas.potthast@math.uni-paderborn.de}

\subjclass[2020]{Primary 11R45; Secondary 11R37 11S31}

\begin{abstract} 
	We determine the distribution of discriminants of wildly ramified elementary-abelian extensions of local and global function fields in characteristic $p$. For local and rational function fields, we also give precise formulae for the number of elementary-abelian extensions with a fixed discriminant divisor, which describe a local-global principle.
\end{abstract}

\maketitle
 
\renewcommand{\theenumi}{(\roman{enumi})}

\section{Introduction and notations}
In 1989, David Wright determined the distribution of discriminants of abelian extensions of global fields in \cite{Wr}. Wright, however, restricted his considerations to the case of Galois groups of order coprime to the characteristic. For the remaining case, he suggested that although the calculations are more complicated, they should give analogous results. 

Much progress has particularly been made in the number field case for which Gunter Malle proposed a precise conjecture -- also for non-abelian groups \cite{Malle1}, \cite{Malle2}. Meanwhile, the conjecture has been proven for instance for the symmetric groups $S_3$ by Harold Davenport and Hans Heilbronn \cite{DH71} and $S_4, S_5$ by work of Manjul Bhargava \cite{Bha1}, \cite{Bha2}. Additionally, Malle's conjecture can be generalised directly to global function fields under the assumption that the order of the Galois group is coprime to the characteristic of the field. Jordan S. Ellenberg and Akshay Venkatesh provide heuristic evidence for Malle's conjecture in this setting using geometric computations on Hurwitz schemes \cite{EV}. Furthermore, Jordan S. Ellenberg, TriThang Tran and Craig Westerland give an upper bound in the weak Malle conjecture for finite Galois extensions of the rational function field $\Fqtg$ in \cite{ETW}. 

The case of Galois groups whose order is divisible by the characteristic $p$, however, stayed %widely 
untouched until Thorsten Lagemann studied the cases of abelian $p$-groups over local function fields and cyclic $p$-groups over global function fields in \cite{La1}. The corresponding easier case when counting by conductor in the case of local function fields was settled by Jürgen Klüners and Raphael Müller \cite{KM}. 
Lagemann's results allowed him to confirm Wright's suggestion for cyclic $p$-groups on the one hand and to disprove it for all non-cyclic abelian $p$-groups using an estimate from the local results on the other hand. In particular, the naive generalisation of Malle's conjecture already fails for non-cyclic abelian $p$-groups over global function fields of characteristic $p$.

Additionally, Lagemann explored the asymptotic behaviour of elementary-abelian extensions of global function fields, counted by conductor, in \cite{La2}. His results also provide a lower and an upper bound for the asymptotics, counted by discriminant. He was later able to extend his results to all abelian $p$-groups, counted by conductor, and to provide a precise conjecture for the asymptotics of the abelian $p$-groups counted by discriminant \cite{La3}.
 %For simple cyclic groups $C_p$ he could even determine the exact asymptotics. 
 
The aim of this paper is to extend the results achieved by Wright and Lagemann and to compute the asymptotics of elementary-abelian extensions of global function fields of characteristic $p$, counted by discriminant. This case turns out to satisfy Lagemann's conjecture for all elementary-abelian $p$-groups, except for the Klein four group $C_2^2$ in the case of characteristic $2$, for which the conjecture breaks because of additional logarithmic factors. Additionally, we prove a local-global principle for the number of elementary-abelian extensions of the rational function field $\Fqtg$ with a fixed discriminant divisor and we compute the asymptotics of elementary-abelian extensions of local function fields. The asymptotic behaviour of local function fields was originally determined by Lagemann and was only published in his dissertation in German in \cite{La1}. In this article, we reproduce his results by a new approach which also includes local algebras in preparation for the case of rational function fields.

\subsection{Notations} \hfill \\
Now, we define some general notations. Let $G = C_p^r$ be an elementary-abelian group of rank $r$ and $F$ a local or global function field of characteristic $p$ over some finite field $\F_q$ of cardinality $q = p^n$. We denote the divisor group by $\Div_F$ and its neutral element by $\mathfrak{1}$. In analogy to the number field case, we write $\Div_F$ multiplicatively. If $F$ is a global function field, we denote the genus of $F$ by $g_F$, the Picard group by $\Pic_F$ and the Picard group of divisors of degree $0$ by $\Pic^0_F$, which we also refer to as the \emph{class group of $F$}, again in analogy to the number field case.

For a fixed field extension $K/F$, we denote its discriminant divisor by $\partial(K/F)$ and in the case of an abelian extension its conductor divisor by $\mathfrak{f}(K/F)$. We are mainly interested in the counting function
\begin{align*}
	Z(F,G;X):=\vert \left\{K/F \ : \  \Gal(K/F) \cong G,\ \vert \vert {\partial(K/F)} \vert \vert = X \right\} \vert
\end{align*}
of Galois field extensions $K/F$ in a given algebraic closure such that there is an isomorphism $\Gal(K/F) \cong G$ and the norm of the discriminant $\vert \vert {\partial(K/F)} \vert \vert = q^{\deg(\partial(K/F))}$ equals some $X = q^m$. The corresponding Dirichlet series
\begin{align*}
	\Phi(F,G;s) = \sum\limits_{\Gal(K/F) \cong G} \vert \vert {\partial(K/F)} \vert \vert^{-s} = \sum\limits_{m \geq 0} c_m q^{-ms}
\end{align*}
encodes the counting function in its coefficients.

\subsection{Main results of this article} \hfill \\
For the asymptotics of $C_p^r$-extensions of local function fields, the results are stated in Proposition \ref{LocalMeromorphicContinuation}, Corollary \ref{PhiRational} and Theorem \ref{LocalAsymptotics}. 

In the case of global function fields, Theorem \ref{ContinuationDirichletSeries} provides a meromorphic continuation of the Dirichlet series $\Phi(\Fqtg,C_p^r;s)$ for the important special case of a rational function field $\Fqtg$ as base field and Theorem \ref{ContinuationF} provides a meromorphic continuation of $\Phi(F,C_p^r;s)$ in the general case. By applying a Tauberian theorem for power series, these meromorphic continuations finally yield the main result of this article on the asymptotics of $C_p^r$-extensions of global function fields:

\begin{theorem} \label{MainResult}
	Let $F$ be a global function field of characteristic $p$ and $G = C_p^r$. Furthermore, let 
	\begin{align*}
	a_p(C_p^r) &:= \frac{1+r(p-1)}{p(p^r -1)}, \\
	b_p(C_p^r) &:= 
	\begin{cases}
		p-1 & \text{for } r=1, \\
		4 & \text{for } r=p=2,\\
		1 & \text{else,}
	\end{cases} \\
\intertext{and} \\
	L(C_p^r) &:= \begin{cases}
		(p-1)M & \text{for } r=1, \\
		12 & \text{for } r=p=2, \\
		p(p^r - 1) & \text{else,}
	\end{cases}
	\end{align*}
where $M$ is the least common multiple of the numbers $2, \ldots, p$. Finally, let $X = q^m$ and $m \to \infty$. 
	Then, there exists a real polynomial $P_{m \bmod L(C_p^r), F}$ of degree at most $b_p(C_p^r) - 1$ depending on the residue class of $m$ modulo $L(C_p^r)$ and the base field $F$ such that the number of $C_p^r$-extensions of $F$ when counted by discriminants satisfies
	\begin{align} \label{GA}
		Z(F, C_p^r; X) = P_{m \bmod L(C_p^r), F}(\log(X)) \cdot X^{a_p(C_p^r)} + O\left(X^{\frac{1 - \frac{1}{p} + r(p-1)}{p(p^{r} -1)} + \epsilon}\right)
	\end{align}
	for every $\epsilon > 0$. The polynomial $P_{m \bmod L(C_p^r), F}$ has degree $b_p(C_p^r) - 1$ for at least one residue class and its leading coefficient is always positive if $P_{m \bmod L(C_p^r), F}$ has this degree.
	\end{theorem}

The proof of Theorem \ref{MainResult} in the important special case $F = \Fqtg$ can be found in Section \ref{ProofSpecialCase}. The general case is proven at the end of Section \ref{ArbitraryBaseFields}.

\subsection{Overview of this article} \hfill \\
The general structure of the article and the results are described below. 

Firstly, we describe the discriminant of a given $C_p^r$-extension of a fixed base field $F$ in terms of the conductors of the intermediate extensions of degree $p$ by applying the well-known conductor-discriminant formula in Section \ref{CDF1}. In Section \ref{Parametrisation}, we provide a parametrisation of elementary-abelian extensions of local and global function fields in terms of the Artin-Schreier theory. Additionally, we describe the ramification behaviour of an elementary-abelian extensions by the ramification behaviour of its intermediate extensions of degree $p$. Section \ref{Combinatorics} contains some combinatorial definitions that we will need later for the counting functions. 

In Section \ref{CFFqtg}, we initially focus on extensions of local function fields and rational function fields since these fields obey a local-global principle and we are able to give precise formulae for the number of their $C_p^r$-extensions with a fixed discriminant divisor, which is generally difficult for arbitrary global function fields. This also has the advantage of avoiding further technical difficulties like non-trivial class groups and Selmer groups as discussed in \cite{La2}.
In Proposition \ref{RS}, we describe systems of representatives for the quotient spaces $F_{\idp}/\wp(F_{\idp})$ of a local function field $F_{\idp}$ by the Laurent series expansion and $\F_q(t)/\wp(\F_q(t))$ of a rational function field $\Fqtg$ by the partial fraction decomposition of rational functions where $\wp(x) := x^p - x$ is the so called \emph{Artin-Schreier operator}. These systems of representatives lead to an explicit description of $C_p$-extensions of $F_{\idp}$ and $\F_q(t)$ and their assigned conductors and discriminants and provide a feasible and explicit access to the class field theory of these fields. In fact, the systems of representatives also indicate a local-global principle for rational function fields. 
 
Subsequently, we precisely compute the number of $C_p^r$-extensions of $F_{\idp}$ and $\Fqtg$, respectively, with a fixed discriminant divisor and decompose the resulting arithmetic function as a linear combination of arithmetic functions in Proposition \ref{CountingFunction}. For $\Fqtg$, these arithmetic functions in the decomposition are multiplicative and yield a presentation of the corresponding Dirichlet series as a linear combination of Euler products in Proposition \ref{LCEuler}, whose Euler factors describe the local counts at the respective places. By analysing the Euler factors, we obtain a meromorphic continuation of these to the whole complex plane in Proposition \ref{LocalMeromorphicContinuation}, which also yields a meromorphic continuation of the Dirichlet series corresponding to the counting function of a local function field to the whole complex plane in Corollary \ref{PhiRational}. We use this to determine the asymptotic behaviour of elementary-abelian extensions of local function fields in Theorem \ref{LocalAsymptotics}. 

By analysing the Euler products in the case of a rational function field, we work out a meromorphic continuation to a sufficiently large domain in Proposition \ref{ContinuationEulerProduct}. Afterwards, we compute a meromorphic continuation for the Dirichlet series corresponding to the counting function of the field $\Fqtg$ in Theorem \ref{ContinuationDirichletSeries} and we prove the important special case $F = \Fqtg$ of Theorem \ref{MainResult}. 

Finally, in Section \ref{ArbitraryBaseFields} we consider elementary-abelian extensions of arbitrary global function fields. We essentially follow the approach presented in \cite{Wr}; for a different presentation see also \cite{Woo}. By this approach, we obtain a decomposition of the Dirichlet series corresponding to the counting function of the respective field as a linear combination of Euler products depending on the quotient $\prod\limits_{i=1}^{r} \OO_S^{\times}/ \left(\OO_S^{\times}\right)^p$ of $S$-units for some sufficiently large set of places $S$ generating the Picard group. This decomposition of the Dirichlet series is presented in Proposition \ref{WrightDecom} and Proposition \ref{PropDecom}. In Proposition \ref{TrivialCharacter}, we prove that the Euler products arising for the trivial class in the quotient of $S$-units coincide with those we have already considered in Section \ref{CFFqtg}. In the remainder of Section \ref{ArbitraryBaseFields}, we deal with the Euler products arising from non-trivial classes in this quotient and we provide a holomorphic continuation to a sufficiently large domain for them in Proposition \ref{AsympEst}. Finally, we obtain a meromorphic continuation of the Dirichlet series corresponding to the counting function of the respective field in Theorem \ref{ContinuationF} and we prove Theorem \ref{MainResult} for arbitrary global function fields $F$.

\section*{Acknowledgements}
%The acknowledgements are omitted in accordance with the double anonymous review procedure.
My special thanks go to Jürgen Klüners for his excellent supervision and helpful discussions about this work. I would also like to thank Fabian Gundlach for proofreading an earlier version and for very productive discussions on Wright's paper and the last section of this article. Finally, I would like to thank the anonymous referee for the thorough review and the many helpful comments that improved the article. \\
The work was funded by the Deutsche Forschungsgemeinschaft (DFG, German Research Foundation) — SFB-TRR 358/1 2023 — 491392403.

\newpage

\section{Application of the conductor-discriminant formula} \label{CDF1}
For a given $C_p^r$-extension of a fixed local or global function field $F$, we can describe the discriminant of this extension in terms of the conductors of its intermediate fields.

\begin{proposition} \label{CDF}
Let $K/F$ be a Galois extension with group $G = C_p^r$. Then, the discriminant is
\begin{align}
\partial(K/F) = \prod\limits_{\substack{F \leq F' \leq K \\ \Gal(F'/F) \cong C_p}} \mathfrak{f}(F'/F)^{p-1} = \prod\limits_{\substack{F \leq F' \leq K \\ \Gal(F'/F) \cong C_p}} \partial(F'/F).
\end{align}
\end{proposition}

\begin{proof}
	Let $\chi \in \hat{G} \cong C_p^r$ be a non-trivial character. By the fundamental homomorphism theorem, we have $G/\Ker(\chi) \cong \im(\chi) \leq \C^{\times}$. 
	Since $\im(\chi)$ is a finite subgroup, we deduce that it is cyclic and therefore isomorphic to $C_p$. In particular, we have $K^{\Ker(\chi)} = F'$ with $\Gal(F'/F) \cong C_p$. 
	
	Let $\chi'$ be another non-trivial character with $\ker(\chi) = \ker(\chi')$. Then, the fundamental homomorphism theorem yields $\chi' = \bar{\chi'} \circ \bar{\chi}^{-1} \circ \chi$ wherein $\bar{\chi}: G/\Ker(\chi) \rightarrow \im(\chi)$ and $\bar{\chi'}: G/\Ker(\chi') \rightarrow \im(\chi')$ denote the respective canonical isomorphisms. The map $\bar{\chi'} \circ \bar{\chi}^{-1}$ is an automorphism of $\im(\chi) = \im(\chi') \cong C_p$. 
	
	We conclude that $p-1$ non-trivial characters correspond to each $C_p$ intermediate field. By an application of the abelian version of the conductor-discriminant formula \cite[Propositions 11.9, 11.10, p. 534f]{Neu}, we complete the proof. 	
\end{proof}

%We remark that the previous Proposition applies to global as well as local function fields.

\section{Description of Artin-Schreier extensions of function fields} \label{Parametrisation}
For a given local or global function field $F$ of characteristic $p$, we aim to describe its finite Galois extensions with elementary-abelian Galois group, the so called \emph{Artin-Schreier extensions}. 

As our base field $F$ has characteristic $p$ and contains the finite field $\F_p$, we can naturally consider $F$ as an $\F_p$-vector space. Let $\wp \colon \bar{F} \rightarrow \bar{F}, x \mapsto x^p -x$ denote the so called \emph{Artin-Schreier operator}. Since the operator is $\F_p$-linear, $\wp(F)$ is an $\F_p$-subspace of $F$ and we can consider the quotient space $F/\wp(F)$.

\begin{proposition} \label{QuotientSpace}
	The map $U \mapsto F(\wp^{-1}(U))$ is a bijection between $\F_p$-subspaces of $F/\wp(F)$ and abelian extensions of $F$ of exponent $p$. 
	
	In particular, subspaces of $U$ correspond to subfields of $F(\wp^{-1}(U))$ and \linebreak $r$-dimensional $\F_p$-subspaces of $F/\wp(F)$ are in bijection with the Artin-Schreier extensions with Galois group $C_p^r$.
\end{proposition}

\begin{proof}
For a proof see \cite[Theorem 8.3, p. 296]{Lang}.
\end{proof}

\begin{corollary} \label{CpGen}
	Let $K/F$ be a $C_p^r$-extension and $\overline{u_i}$ for $1 \leq i \leq r$ generators of the corresponding subspace $U \leq F/\wp(F)$. Then, the $\frac{p^r -1}{p-1}$ intermediate $C_p$-extensions correspond to the $1$-dimensional subspaces of $U$ generated by the elements
	\begin{align} \label{CpGenerator}
	\overline{u_k} + \sum\limits_{j=k+1}^{r} \lambda_j \overline{u_j} \ \text{with} \ \lambda_j \in \F_p
	\end{align}
for $1 \leq k \leq r$.
\end{corollary}

\begin{proof}
It is a well-known fact that we have $\frac{p^r -1}{p-1} = \sum\limits_{j=0}^{r-1} p^{j}$ intermediate $C_p$-extensions. By Proposition \ref{QuotientSpace}, these correspond to $1$-dimensional subspaces in $U$. Since we have $p^{r-k}$ elements for a given index $k$ in \eqref{CpGenerator}, we have a total of $\sum\limits_{j=0}^{r-1} p^{j}$ elements given. Due to the normalization of the first coefficient, it is obvious that each element generates a different $1$-dimensional subspace.
\end{proof}

In the case of $C_p$-extensions of $F$, we provide a more detailed description including the different exponent and the conductor exponent for a given place $\idp \in \PP_F$:

\begin{proposition} \label{ASC}
Let $F$ be a function field of characteristic $p$. Furthermore, let $u \in F$ and $F' = F(y_u)$ with $\wp(y_u) =u$. Then, $F'/F$ is a $C_p$-extension if and only if $u \not \in \wp(F)$ and $F' = F$ if $u \in \wp(F)$. 

For $\idp \in \PP_F$, we define the integer $\ASC(\idp,u)$, the so called \emph{Artin-Schreier conductor}, by
\begin{align} 
\ASC(\idp,u) = \begin{cases}
	m,  & \exists z \in \wp(F) \ \text{with} \ \nu_{\idp}(u-z) = -m < 0 \ \text{and} \ p \nmid m, \\
	-1, & \exists z \in \wp(F) \ \text{with} \ \nu_{\idp}(u-z) \geq 0.
\end{cases} \label{ASCFormula}
\end{align}
We note that $\ASC(\idp,u)$ is well-defined by \cite[Lemma 3.7.7, p. 125]{St}. Then, the following statements apply:
\begin{enumerate}
	\item[a)] The place $\idp$ is unramified in $F'/F$ if and only if $\ASC(\idp,u) =-1$.
	\item[b)] The place $\idp$ is totally ramified in $F'/F$ if and only if $\ASC(\idp,u) > 0$. Denote by $\mathfrak{P}$ the unique place of $F'$ lying over $\idp$. Then, we have
	\begin{align*}
	d(\mathfrak{P} \vert {\idp}) = (p-1)(\ASC(\idp,u) + 1)
	\end{align*}
for the different exponent.
\item[c)] In the situation of b), we have
\begin{align*}
	c(\mathfrak{P} \vert {\idp}) = \ASC(\idp,u) + 1
\end{align*}
for the conductor exponent.
\end{enumerate}

\end{proposition} 

\begin{proof}
The proofs of a) and b) can be found in \cite[Proposition 3.7.8, p. 127]{St}. Since the different exponent of a totally ramified extension coincides with the discriminant exponent, we obtain c) by applying Proposition \ref{CDF}.
\end{proof}

\begin{remark} \label{ReducedRep}
	Since every representative of a class $\overline{u} \in F/\wp(F)$ has the same Artin-Schreier conductor, we can also write $\ASC(\idp, \overline{u})$ for the class. In particular, we can choose a representative for a given class with $z=0$ in \eqref{ASCFormula}. We call every representative with this property \emph{reduced}.
\end{remark}

In order to describe the ramification behaviour of a $C_p^r$-extension at a particular place $\idp$, we now introduce the concept of a chain which we will use to specify the conductor exponents of the extension.

\begin{definition} \label{Chain}
	Let $h \in \N_0$. We call a tuple $\mathfrak{C} = (c_1, \ldots, c_h) \in \N^h$, which fulfils the property $c_1 \geq \ldots \geq c_h > 0$, a \emph{chain} and we call $h$ the \emph{length} of the chain. Note that we only have the empty chain $()$ in the case $h = 0$. 
	
	Additionally, for any $x \in \N_0$ the subchain of $\mathfrak{C}$ for whom all components of the chain are $ > x$ is denoted by $\mathfrak{C}_{>x}$. Similarly, $\mathfrak{C}_{x}$ denotes the excerpt of $\mathfrak{C}$ where all components are equal to $x$.
\end{definition}

\begin{proposition} \label{ConductorChain}
%	Let $K/F$ be a $C_p^r$-extension and $\idp \in \PP_F$ a fixed place. Then, we can choose generating intermediate $C_p$-extensions $F_i$ for $1 \leq i \leq r$ with conductor exponents $c_i(\idp)$ so that we have a decreasing chain $\mathfrak{C}(\idp) := \left(c_1(\idp), \ldots, c_r(\idp)\right) \in \N_0^{r}$ such that
%	\begin{align} \label{ConductorStep}
%	c_i(\idp) \geq c_{i+1}(\idp) \quad \text{for} \quad 1 \leq i \leq r-1
%	\end{align}
%and the $p^{r-i}$ intermediate $C_p$-extensions of $K/F$ in $F_i \cdots F_r$, that are not already contained in $F_{i+1} \cdots F_r$ with $F_{i+1} \cdots F_r = \varnothing$ for $i = r$, have the conductor exponent $c_i(\idp)$. We call the elements of $\mathfrak{C}(\idp)$ \emph{conductor steps} and $r$ the \emph{length} of the chain. \\
%\\
%\\
Let $K/F$ be a $C_p^r$-extension and $\idp \in \PP_F$ a fixed place. Then, we can choose generating intermediate $C_p$-extensions $F_i$ for $1 \leq i \leq r$ with conductor exponents $c_i(\idp) \in \N_0$ such that
\begin{align} \label{ConductorStep}
	c_i(\idp) \geq c_{i+1}(\idp) \quad \text{for} \quad 1 \leq i \leq r-1
\end{align}
and the $p^{r-i}$ intermediate $C_p$-extensions of $K/F$ in $F_i \cdots F_r$, that are not already contained in $F_{i+1} \cdots F_r$ with $F_{i+1} \cdots F_r = F$ for $i = r$, have the conductor exponent $c_i(\idp)$. 

In particular, there is a maximum index $h \in \N_0$ with $c_h(\idp) > 0$ which gives us the chain $\mathfrak{C}(K, \idp) := \left(c_1(\idp), \ldots, c_h(\idp)\right) \in \N^{h}$. This chain completely describes the ramification behaviour of the extension $K/F$ at the place $\idp$. 

\end{proposition}

\begin{proof}
	We prove the claim by induction on $r$. The case $r=1$ is trivial. 
	
	Assume the property is true for all intermediate extensions of $K/F$. Let $F_r/F$ be an intermediate $C_p$-extension with conductor exponent $c_{F_r}(\idp)$ such that $c_{F_r}(\idp)$ is the minimal conductor exponent in $K/F$. Furthermore, let $L/F$ be an intermediate extension with Galois group isomorphic to $C_p^{r-1}$ and corresponding tuple $\left(c_1(\idp), \ldots, c_{r-1}(\idp) \right) \in \N_0^{r-1}$ of conductor exponents such that $LF_r = K$. Let $h \in \N_0$ denote the maximum index with $c_h(\idp) > 0$. That is $c_{h+1}(\idp) = \ldots = c_{r-1}(\idp) = 0$ and $\mathfrak{C}(L, \idp) = \left(c_1(\idp), \ldots, c_h(\idp)\right) \in \N^{h}$. Finally, let $u_i \in F$ be reduced representatives of generators of the corresponding $1$-dimensional subspaces to $F_i$ for $1 \leq i \leq r$. 
	
	%\underline{Case 1:} Assume, we have an index $i$ with $c_{i-1}(\idp) > c_r(\idp) > c_i(\idp)$. We permute the indices as $\tau(j) = j$ for $1 \leq j \leq i-1$, $\tau(r) = i$ and $\tau(j) = j+1$ for $i \leq j \leq r-1$ and obtain a new decreasing chain with $c_{i-1}(\idp) > c_{i}(\idp) > c_{i+1}(\idp)$. By Proposition \ref{ASC}, we have 
	%\begin{align*}
	%	\nu_{\idp}(u_{i-1}) < \nu_{\idp}(u_{i}) < \nu_{\idp}(u_{i+1}).
	%\end{align*}
	%For each $i+1 \leq j \leq r$, we have $p^{r-1-(j-1)} = p^{r-j}$ intermediate $C_p$-extensions with conductor exponent $c_j(\idp)$ by our induction hypothesis. We apply Corollary \ref{CpGen} for the index $i$ and obtain by the ultrametric inequality
	%\begin{align*}
	%\nu_{\idp}\left(u_{i} + \sum\limits_{j=i+1}^{r} \lambda_j u_j\right) = \nu_{\idp}(u_{i}).
	%\end{align*}
	%Therefore, we have $p^{r-i}$ intermediate $C_p$-extensions with conductor exponent $c_{i}(\idp)$. Analogously, we obtain $p^{r-j}$ intermediate $C_p$-extensions for $1 \leq j \leq i-1$, since we have $p^{r-1-j}$ by our induction hypothesis and the additional generator $u_r$ yields another factor $p$. \\
	\underline{Case 1:} Assume, we have $c_{r-1}(\idp) > c_{F_r}(\idp)$ which implies $h = r-1$. By Proposition \ref{ASC}, it also applies $\nu_{\idp}(u_{r-1}) < \nu_{\idp}(u_{r})$. For $1 \leq i \leq r-1$, we apply Corollary \ref{CpGen} and obtain
	\begin{align*}
		\nu_{\idp}\left(u_{i} + \sum\limits_{j=i+1}^{r} \lambda_j u_j\right) = 	\nu_{\idp}\left(u_{i} + \sum\limits_{j=i+1}^{r-1} \lambda_j u_j\right) = \nu_{\idp}(u_i)
	\end{align*}
	by the ultrametric inequality. Therefore, we have $p^{r-1-i} \cdot p = p^{r-i}$ intermediate $C_p$-extensions with conductor exponent $c_{i}(\idp)$ in $F_i \cdots F_r$, that are not already contained in $F_{i+1} \cdots F_r$, by the induction hypothesis and the additional factor $p$ coming from the generator $u_r$. 
	Due to the assumption $\nu_{\idp}(u_{r-1}) < \nu_{\idp}(u_{r})$ in this case, there is exactly one intermediate $C_p$-extension with conductor exponent $c_{F_r}(\idp)$. Consequently, we set $c_r(\idp) := c_{F_r}(\idp)$. 
	
	 %If we have $c_r(\idp) = 0$, we have the chain $\mathfrak{C}(K, \idp) = \mathfrak{C}(L, \idp) = \left(c_1(\idp), \ldots, c_{r-1}(\idp)\right) \in \N^{r-1}$. Otherwise, we have $\mathfrak{C}(K, \idp) = \left(c_1(\idp), \ldots, c_{r}(\idp)\right) \in \N^{r}$.  \\ 
	\underline{Case 2:} Assume, we have an index $i < r$ with $c_{i-1}(\idp) > c_{i}(\idp) = \ldots = c_{F_r}(\idp)$. Analogously to the previous case, we have $\nu_{\idp}(u_{i}) = \ldots = \nu_{\idp}(u_{r}) > \nu_{\idp}(u_{i-1})$
	so that the indices $1 \leq j \leq i-1$ can be treated as in Case 1. Therefore, we further assume $i \leq j \leq r$ without loss of generality.
	By the ultrametric inequality, canceling might occur in this case and we have
	\begin{align*} %\label{Canceling}
		\nu_{\idp}\left(u_{j} + \sum\limits_{k=j+1}^{r} \lambda_k u_k\right) \geq \nu_{\idp}\left(u_{j} + \sum\limits_{k=j+1}^{r-1} \lambda_k u_k\right) = \nu_{\idp}(u_j).
	\end{align*}
	Therefore, we consider 
	\begin{align*}
	M := \max\limits_{i \leq j \leq r} \left(\max\limits_{(\lambda_{j+1}, \ldots, \lambda_r) \in \F_p^{r-j}} \nu_{\idp}\left(u_{j} + \sum\limits_{k = j+1}^{r} \lambda_k u_k\right) \right) \geq \nu_{\idp}(u_r).
	\end{align*}
	If we have $M = \nu_{\idp}(u_r)$, then canceling does not occur and we can again proceed as we did in Case 1 and we set $c_r(\idp) := c_{F_r}(\idp)$. However, if we have $M > \nu_{\idp}(u_r)$, we replace $u_r$ with a generator for which the maximum $M$ is reached. Now, we can once again apply Case 1 and finally choose $c_r(\idp) := -M + 1$ for $M < 0$ respective $c_r(\idp) := 0$ for $M \leq 0$.
	
	In both cases, the following applies to the chains: If we have $c_r(\idp) = 0$, we obtain the chain $\mathfrak{C}(K, \idp) = \mathfrak{C}(L, \idp) = \left(c_1(\idp), \ldots, c_{h}(\idp)\right) \in \N^{h}$. Otherwise, we have $\mathfrak{C}(K, \idp) = \left(c_1(\idp), \ldots, c_{r}(\idp)\right) \in \N^{r}$.
%	\\
%	Let us first assume that we have equality in \eqref{Canceling}. Then, analogous calculations as in Case 1 apply for the indices $1 \leq j \leq i-1$ and $k+1 \leq j \leq r$. For the remaining indices $i \leq j \leq k$, we also have another factor $p$ by the additional generator $u_r$. \\
%	If we have canceling in \eqref{Canceling}, we consider $\max\limits_{(\lambda_i, \ldots, \lambda_k) \in \F_p^{k-i+1}} \nu_{\idp}\left(u_{r} + \sum\limits_{j=i}^{k} \lambda_j u_j\right)$ and replace $u_r$ by this maximum and start over. We repeat this exchange process until we can either apply Case 1 or the previous sub-case without canceling. Finally, we note that the exchange process terminates in finitely many steps.
\end{proof}

\begin{remark}
The choice of the generating $C_p$-extensions in Proposition \ref{ConductorChain} depends on the given place $\idp$. In general, it is not possible to choose the same set of $C_p$-fields for all places $\idp \in \PP_F$ simultaneously. In particular, the chains of two different places within a fixed $C_p^r$-extension may have a different shape.
\end{remark}

\begin{remark} \label{RelationCondExpRamGroups}
In fact, the conductor exponents encode the jumps of the higher ramification groups of the respective place in the upper numbering which are given by $c_i(\idp) - 1$, the Artin-Schreier conductors. More precisely, $F_i \cdots F_r$ is the subfield of $K$ which is fixed by the higher ramification group $G_{\idp}^{(c_i(\idp))}$ if we have $c_{i-1}(\idp) > c_i(\idp) > 0$. 

If $0$ appears as a conductor exponent, this means that an unramified component is included in the extension. We do not distinguish here between inertia and splitting, as this is not of interest for our further investigations. 

Note that we always have $G_{\idp}^{(0)} = G_{\idp}^{(1)}$ as we do not allow tamely ramified parts in our extensions under consideration.
\end{remark}

%\begin{remark} \label{C_pNonMultiplicative}
%In \cite[Lemma 5.2]{La2}, Thorsten Lagemann presents formulae for the number of $C_p^r$-extensions $K/F$ of a global function field $F$ with a prescribed conductor divisor $\mathfrak{m}$. The formula given in case c) in his Lemma apparently describes a multiplicative formula. Case a) and case b) with $\widetilde{c}_{\mathfrak{m}_0} \neq 0$, however, deviate from multiplicativity. The origin of these deviations lies in a non-trivial $p$-torsion of the Picard group and associated non-trivial Selmer groups $S_{\mathfrak{m}}$ for some divisors $\mathfrak{m}$. The details can be found in \cite{La2}. \\
%For the field $\F_q(t)$ of rational functions, the setting $g_{\F_q(t)} = 0$, $\widetilde{c}_{\mathfrak{1}} = 0$ and $\Pic[p] \cong \{1\}$ yields that all modules covered by case b) in the Lemma also behave multiplicatively and $\mathfrak{m} = \mathfrak{1}$ is the only module that has to be treated differently due to an inclusion-exclusion principle for the trivial extension. In fact, analogous results as for the rational function field can be obtained for a local function field by considering only one place $\idp$ and $\mathcal{M} = \{\mathfrak{1}\}$. \\
%In the special case of $C_p$-extensions, Lagemann's Lemma with slight adjustments also describes the number of $C_p$-extensions counted by discriminant as the discriminant of an $C_p$-extension with conductor $\mathfrak{m}$ is $\partial = \mathfrak{m}^{p-1}$ by Proposition \ref{CDF}.
%\end{remark}

\section{Conductor chains and compositions} \label{Combinatorics}

%In dependence on the respective chain $\mathfrak{C}(K, \idp)$ of length $h$ in Proposition \ref{ConductorChain}, we can associate a uniquely determined composition $\omega_{\mathfrak{C}(K, \idp)}$ of $h$. This association is going to be crucial for our counting function later on.
In this section, we would like to assign a clearly defined composition $\omega_{\mathfrak{C}(K, \idp)}$ of $h$ to a chain $\mathfrak{C}(K, \idp)$ of length $h$ according to Proposition \ref{ConductorChain}. On the one hand, this composition makes it easier to describe the associated discriminant of the chain, and on the other hand, chains with the same composition are subject to the same combinatorial behaviour when counting.

\begin{definition} \label{Composition}
	A \emph{composition} of $h \in \N_0$ is an ordered partition of $h$. That is a tuple $\omega = (a_1, \ldots, a_{\lambda})$ of positive integers such that $\sum\limits_{i=1}^{\lambda} a_i = h$. We call $\lambda_{\omega} := \lambda$ the \emph{length} of $\omega$ and $a_1(\omega) := a_1, \ldots, a_{\lambda}(\omega) := a_{\lambda}$ its \emph{blocks}. We further define $A_i(\omega) := \sum\limits_{k=1}^{i} a_k(\omega)$ for $0 \leq i \leq \lambda_{\omega}$. If the composition is clearly determined, we henceforth suppress the dependency on $\omega$ to shorten the notation. 
	
	Let $\Omega_h$ denote the set of compositions of $h \in \N_0$. We note that the set $\Omega_0$ only contains the empty composition $()$. For $h \geq 1$, we recall the classical fact $\vert {\Omega_h} \vert = 2^{h-1}$. Furthermore, let $\Omega^f := \bigcup\limits_{h = 0}^{f} \Omega_h$ be the set of all compositions of numbers $\leq f$ for any $f \in \N_0$. 
	
	For $x, y \in \N_0$ satisfying $x\geq y$, we recall the \emph{Gaussian binomial} $\binom{x}{y}_p := \frac{\prod\limits_{i = x-y+1}^{x} \left(p^i - 1 \right)}{\prod\limits_{i=1}^{y} \left(p^i -1\right)}$, which is a natural number and gives the number of $y$-dimensional $\F_p$-subspaces of the $x$-dimensional $\F_p$-vector space $\F_p^x$. 
	%	\begin{align*}
		%		\gamma_{\omega} := \frac{\vert {\Aut\left(C_p^f\right)} \vert}{(p-1)^f \prod\limits_{i=1}^{\lambda_{\omega}} \prod\limits_{j=1}^{a_i} \sum\limits_{k=1}^{j} p^{f - A_{i-1} - k}} \in \N_0
		%	\end{align*}
	For a fixed composition $\omega \in \Omega_h$, we define
	\begin{align*}
		\gamma_{\omega} := \prod\limits_{i=1}^{\lambda_{\omega}}  \binom{A_{i}(\omega)}{A_{i-1}(\omega)}_p = \prod\limits_{i=1}^{\lambda_{\omega}}  \binom{A_{i}(\omega)}{a_i(\omega)}_p \in \N,
	\end{align*}
	which is the number of flags $V_1 \leq \ldots \leq V_{\lambda_{\omega}}$ in the vector space $\F_p^h$ satisfying $\dim(V_i) = A_i(\omega)$ for all $1 \leq i \leq \lambda_{\omega}$. Accordingly, $\gamma_{\omega}$ is a polynomial expression in $p$, which we call \emph{combinatorial constant} of $\omega$. 
\end{definition}

\begin{lemma} \label{gammainvariance}
Let  $\omega = (a_1, \ldots, a_{\lambda}) \in \Omega_h$ be a composition and let the symmetric group $S_{\lambda}$ act on $\omega$ by permutation of the blocks. That is $\sigma(\omega) := \left(a_{\sigma(1)}, \ldots, a_{\sigma(\lambda)}\right)$ for any $\sigma \in S_{\lambda}$. Then, the combinatorial constant $\gamma_{\omega}$ is invariant under the action of $S_{\lambda}$, meaning that $\gamma_{\omega} = \gamma_{\sigma(\omega)}$ for all $\sigma \in S_{\lambda}$.
\end{lemma}

\begin{proof}
The claim follows by the following short computation:
\begin{align*}
\gamma_{\omega} = \prod\limits_{i=1}^{\lambda_{\omega}}  \binom{A_{i}(\omega)}{a_i(\omega)}_p = \frac{\prod\limits_{i=1}^{\lambda_{\omega}} \prod\limits_{j= A_{i-1}(\omega)+1}^{A_i(\omega)}\left(p^j-1\right)}{\prod\limits_{i=1}^{\lambda_{\omega}} \prod\limits_{j=1}^{a_i(\omega)} \left(p^j-1\right)} = \frac{\prod\limits_{i=1}^{h} \left(p^j-1\right)}{\prod\limits_{i=1}^{\lambda_{\omega}} \prod\limits_{j=1}^{a_i(\omega)} \left(p^j-1\right)}.
\end{align*}
\end{proof}

\begin{definition} \label{ChainComp}
To any chain $\mathfrak{C} = (c_1, \ldots, c_h) \in \N^h$, we associate a uniquely determined composition $\omega_{\mathfrak{C}} = (a_1, \ldots, a_{\lambda}) \in \Omega_h$ such that 
\begin{align*}
(c_1, \ldots, c_h) = (\underbrace{b_1, \ldots, b_1}_{a_1}, \ldots, \underbrace{b_{\lambda}, \ldots, b_{\lambda}}_{a_{\lambda}}) \quad \text{with} \quad  b_1 > \ldots > b_{\lambda} > 0. 
\end{align*}
More precisely, we have $c_{A_i} > c_{A_i + 1} = \ldots = c_{A_{i+1}} > c_{A_{i+1}+1}$ for $0 \leq i \leq \lambda_{\omega_{\mathfrak{C}}}-1$.
\end{definition}

\begin{example}
	Let us consider a chain $\mathfrak{C} = (c_1, c_2, c_3) \in \N^3$ with the property $c_1 > c_2 = c_3$. Then, we associate the composition $(1,2) \in \Omega_3$ to this chain. % if we have $c_3(\idp) > 0$. In the case $c_3(\idp) = 0$, the assigned composition is $(1) \in \Omega_1$.
\end{example}

Naturally, we would like to apply the previous Definitions to an Artin-Schreier extension $K/F$ and its associated conductor exponents, having Proposition \ref{ConductorChain} in mind. 

In Remark \ref{RelationCondExpRamGroups}, we have already discussed the relationship between the associated chain $\mathfrak{C}(K, \idp) = (c_1(\idp), \ldots, c_h(\idp))$ with some $h \leq r$ and the higher ramification groups. Now, we would also like to interpret the associated composition $\omega_{\mathfrak{C}(K, \idp)} = \left(a_1, \ldots, a_{\lambda_{\omega_{\mathfrak{C}(K, \idp)}}} \right) \in \Omega_h$ in terms of the higher ramification groups. 

We consider the filtration of higher ramification groups in the upper numbering and we coarsen it in the sense that we only consider all those ramification groups where jumps actually occur, say $\{1\} = T_0 \subsetneq T_1 \subsetneq T_2 \subsetneq \ldots \subsetneq T_{\lambda} \subseteq G \cong C_p^r$. Then, we have $T_i = G_{\idp}^{(c_{A_i + 1}(\idp))}$ by our considerations in Remark \ref{RelationCondExpRamGroups}. Furthermore, we have $\lambda = \lambda_{\omega_{\mathfrak{C}(K, \idp)}}$, $\dim_{\F_p}(G/T_{\lambda}) = r - h$ and $a_i = \dim_{\F_p}(T_i/T_{i-1})$ as well as $A_i = \dim_{\F_p}(T_i/T_0) = \dim_{\F_p}(T_i)$ for $1 \leq i \leq \lambda_{\omega_{\mathfrak{C}(\idp)}}$.

\begin{lemma} \label{DiscExponent} %Kann auch für beliebiges $F$ formuliert werden
	Let $K/F$ be a $C_p^r$-extension and $\idp \in \PP_{F}$ a fixed place. Furthermore, let $\mathfrak{C}(K, \idp) = (c_1(\idp), \ldots, c_h(\idp)) \in \N^h$ with $h \leq r$ be the chain of conductor exponents according to Proposition \ref{ConductorChain} and $\omega_{\mathfrak{C}(K, \idp)} = \left(a_1, \ldots, a_{\lambda}\right) \in \Omega_h$ its associated composition. 
	Then, the exponent $d_{\idp}(K/F)$ of the $\idp$-part of the discriminant $\partial(K/F)$ is
	\begin{align*}
		d(\mathfrak{C}(K, \idp)) := (p-1) \sum\limits_{i=1}^{\lambda_{\omega_{\mathfrak{C}(K,\idp)}}} \sum\limits_{j=1+A_{i-1}}^{A_i} p^{r-j} c_{A_i}(\idp).
	\end{align*}
\end{lemma}

\begin{proof}
	The claim follows from the Propositions \ref{CDF} and \ref{ConductorChain}.
\end{proof}

Next, we provide an abbreviated notation of typical components that will appear in our counting functions.

\begin{definition} \label{DefinitionCountingFunction}
	For $c \in \N_0$ and $j \in \N_0$, we define
	%	\begin{align*}
		%		\widehat{Z}_{\idp, p^j}(c) := \np^{r_{c}} - p^j \np^{r_{c-1}}
		%	\end{align*}
	\begin{align*}
		\widehat{Z}_{\idp, p^j}(c) := \left\{
		\begin{array}{ll}
			1, & c = 0, \\
			\np^{r(c)} - p^j \np^{r(c-1)},& c>0. \\
		\end{array}
		\right.
	\end{align*}
	where $r(c) := c-1 - \left\lfloor \frac{c-1}{p} \right\rfloor$ with the floor function $\lfloor \cdot \rfloor$. For a chain $\mathfrak{C} = (c_1, \ldots, c_h) \in \N^h$ and its associated composition $\omega_{\mathfrak{C}} = (a_1, \ldots, a_{\lambda}) \in \Omega_h$, we further define
	\begin{align*}
		\widehat{Z}_{\idp}(\mathfrak{C}) := \prod\limits_{i=1}^{\lambda_{\omega_{\mathfrak{C}}}} \prod\limits_{j=0}^{a_i-1} \widehat{Z}_{\idp, p^j}(c_{A_i}).
	\end{align*}
Note that we have $\widehat{Z}_{\idp}(()) = 1$ for the empty chain $()$.
	%Additionally, for any $x \in \N_0$ the expression $\widehat{Z}_{\idp}(\mathfrak{C}_{>x})$ denotes the part of $\widehat{Z}_{\idp}(\mathfrak{C})$ for whom the conductor exponents are $ > x$ and $\widehat{Z}_{\idp}(\mathfrak{C}_{x})$ denotes the part with conductor exponents which are equal to $x$.
\end{definition}

\section{Artin-Schreier extensions of local and rational function fields} \label{CFFqtg}

In this section, we restrict ourselves to Artin-Schreier extensions of a local function field or a rational function field $\Fqtg$ over a finite field $\F_q$ to temporarily avoid technical difficulties like non-trivial class groups and Selmer groups. For these special cases, we obtain a multiplicativity in the counting functions presented in Proposition \ref{CountingFunction}, which breaks if arbitrary global function fields are admitted as our base field. We note that a deviation from multiplicativity in the general case can also be observed in Lagemann's formulae for counting by conductors in \cite[Lemma 5.2]{La2}. The general case of arbitrary global function fields as base field is discussed later in Chapter \ref{ArbitraryBaseFields} using the adelic methods developed by Wright in \cite{Wr}.

\subsection{A system of representatives for Artin-Schreier extensions} \hfill\\
In the case of a local function field $F_{\idp}$ that is the local completion of some global function field $F$ at a fixed place $\idp \in \PP_F$, we can explicitly describe the quotient space $F_{\idp}/\wp(F_{\idp})$ by applying the characterisation of $F_{\idp}$ as a field of formal Laurent series. 
Similarly, the base field $\F_q(t)$ allows an explicit description of the quotient space $\F_q(t)/\wp(\F_q(t))$ using the partial fraction decomposition of rational functions. 
In fact, these explicit descriptions reveal that nearly a local-global principle holds for local and rational function fields.

\begin{proposition} \label{PFD}
%Let $z(t) = \frac{f(t)}{g(t)} \in \F_q(t)$ with coprime polynomials $f(t), g(t) \in \F_q[t]$ and $g(t) \neq 0$. \\
Let $z \in \F_q(t)$ be a rational function and $\pi_{\idp}$ a uniformiser for a given place $\idp$, that is either a monic irreducible polynomial in $\F_q[t]$ or $\pi_{\infty} = t^{-1}$ for the infinite place. Then, $z$ has a unique representation as
\begin{align*}
z = z_0 + \sum\limits_{\idp \in \PP_{\F_q(t)}} \sum\limits_{i=1}^{n(\idp)} \frac{z_{i,\idp}}{\pi_{\idp}^{i}}
\end{align*}
with $n(\idp) \in \N \cup \{-1\}$, $n(\idp) \in \N$ for only finitely many places $\idp$, $z_0 \in \F_q$, $z_{i, \idp} \in \F_q[t]$ with $\deg(z_{i, \idp}) < \deg(\idp)$ and $z_{n(\idp), \idp} \neq 0$. Note that we can interpret the polynomials $z_{i, \idp}$ as canonical representatives of the finite field $\F_{q^{\deg(\idp)}}$ as we have $\F_q[t]/(\pi_{\idp}) \cong \F_{q^{\deg(\idp)}}$ for $\idp \neq \infty$ and $\mathcal{O}_{\infty}/(\pi_{\infty}) \cong \F_q = \F_{q^{\deg(\infty)}}$ for $\idp = \infty$.
\end{proposition}

\begin{proof}
A proof is given in \cite[Theorems 5.2 and 5.3, p. 189]{Lang}.
\end{proof}

\begin{proposition} \label{RS}
Let $\mathfrak{R}_{\F_q}$ be a minimal system of representatives of the quotient space $\F_q/\wp(\F_q)$. Such a minimal system of representatives has cardinality $p$. Then, we have:
\begin{enumerate}
	\item[a)] For each $z \in \Fqtg$, the class $\overline{z}$ in the quotient space $\Fqtg/\wp(\Fqtg)$ is uniquely represented by an element in the set
	\begin{align*}
		\mathfrak{R}_{\Fqtg} = \left\{\overline{z_0} + \sum\limits_{\idp \in \PP_{\F_q(t)}} \sum\limits_{\substack{i=1 \\ p \ \nmid \ i}}^{\ASC(\idp, \overline{z})} \frac{z_{i,\idp}}{\pi_{\idp}^{i}} \ : \ \overline{z_0} \in \mathfrak{R}_{\F_q}, z_{i,\idp} \in \F_{q^{\deg(\idp)}}, z_{\ASC(\idp, \overline{z}), \idp} \neq 0 \right\}.
	\end{align*} 
In particular, $\mathfrak{R}_{\Fqtg}$ is a minimal system of representatives for the quotient space $\Fqtg/\wp(\Fqtg)$.
	\item[b)] For each $z \in F_{\idp}$, where $F_{\idp}$ denotes the local completion of some global function field $F$ with field of constants $\F_q$ at a fixed place $\idp \in \PP_{F}$ with uniformiser $\pi_{\idp}$, the class $\overline{z}$ in the quotient space $F_{\idp}/\wp(F_{\idp})$ is uniquely represented by an element in the set
	\begin{align*}
		\mathfrak{R}_{F_{\idp}} = \left\{\overline{z_0} + \sum\limits_{\substack{i=1 \\ p \ \nmid \ i}}^{\ASC(\idp, \overline{z})} \frac{z_{i,\idp}}{\pi_{\idp}^{i}} \ : \ \overline{z_0} \in \mathfrak{R}_{\F_q}, z_{i,\idp} \in \F_{q^{\deg(\idp)}}, z_{\ASC(\idp, \overline{z}), \idp} \neq 0 \right\}.
	\end{align*}
In particular, $\mathfrak{R}_{F_{\idp}}$ is a minimal system of representatives for the quotient space $F_{\idp}/\wp(F_{\idp})$.
\end{enumerate}
Note that all the representatives are reduced simultaneously at all places of the respective field in the sense of Proposition \ref{ASC} and Remark \ref{ReducedRep} as we have $p \nmid \ASC(\idp, \overline{z})$.
\end{proposition}

\begin{proof}
The cardinality of the constant field $\F_q$ is $q = p^n$. It is easy to verify that $\vert {\wp(\F_q)} \vert = p^{n-1}$ applies which yields $\vert {\mathfrak{R}_{\F_q}} \vert = \frac{p^n}{p^{n-1}} = p$. 

\underline{Proof of a)}: For $z \in \F_q(t)$, we consider the partial fraction decomposition given by Proposition \ref{PFD} and its class $\overline{z}$ in the quotient space $\F_q(t)/\wp(\F_q(t))$. We first show that $z$ is represented by an element in $\mathfrak{R}_{\F_q(t)}$. Let $z = z_0 + \sum\limits_{\idp \in \PP_{\F_q(t)}} \sum\limits_{i=1}^{n(\idp)} \frac{z_{i,\idp}}{\pi_{\idp}^{i}}$ be the partial fraction decomposition. The element $z_0 \in \F_q$ has a unique representative $\overline{z_0} \in \mathfrak{R}_{\F_q}$. For all but finitely many places, we have $\nu_{\idp}(z) \geq 0$ and therefore $n(\idp) = \ASC(\idp, \overline{z}) = -1$. 

We consider a place $\idp$ with $\nu_{\idp}(z) < 0$, now. If we have $p \nmid i$ for all indices, then we have $n(\idp) = \ASC(\idp, \overline{z})$ and the $\idp$-part of $z$ already matches with the representation in $\mathfrak{R}_{\F_q(t)}$. Otherwise, we recursively cancel all indices divisible by $p$, starting with the largest index $i = p \ell$ for some $\ell \in \N$. Since the Frobenius homomorphism is bijective on finite fields, there is a unique canonical representative $g_{p\ell, \idp} \in \F_q[t]/(\pi_{\idp})$ satisfying $g^p_{p\ell, \idp} \equiv -z_{p\ell, \idp} \mod \pi_{\idp}$. By our choice, we have $\deg(g^p_{p\ell, \idp}) < \deg(\pi^p_{\idp})$ and obtain a representation $g^p_{p\ell, \idp} = \sum\limits_{j=0}^{p-1} h_{j, \idp} \cdot \pi_{\idp}^j$ with unique polynomials $h_{j, \idp}$ satisfying $\deg(h_{j, \idp}) < \deg(\pi_{\idp})$ and $h_{0, \idp} = -z_{pl, \idp}$. For the $\idp$-part of $z$, we compute
\begin{align*}
& \sum\limits_{i=1}^{n(\idp)} \frac{z_{i,\idp}}{\pi_{\idp}^{i}} + \wp\left(g_{p\ell, \idp} \cdot \pi_{\idp}^{-\ell}\right) \\
=& \sum\limits_{i=1}^{n(\idp)} \frac{z_{i,\idp}}{\pi_{\idp}^{i}} + \left(\sum\limits_{j=1}^{p-1} h_{j, \idp} \cdot \pi_{\idp}^j  - z_{p\ell, \idp}\right) \cdot \pi_{\idp}^{-p\ell} - g_{p\ell, \idp} \cdot \pi_{\idp}^{-\ell} \\
=& \sum\limits_{\substack{i=1 \\ i \neq p\ell}}^{n(\idp)} \frac{z_{i,\idp}}{\pi_{\idp}^{i}} + \sum\limits_{j=1}^{p-1} h_{j, \idp} \cdot \pi_{\idp}^{j-p\ell} - g_{p\ell, \idp} \cdot \pi_{\idp}^{-\ell}.
\end{align*}
We note that we canceled the index $i = p\ell$, adjusting only the smaller indices. Therefore, we can modify the $\idp$-part of $z$ in finitely many steps to $\sum\limits_{\substack{i=1 \\ p \ \nmid \ i}}^{\ASC(\idp, \overline{z})} \frac{z_{i,\idp}}{\pi_{\idp}^{i}}$. Note that the modification of the $\idp$-part of $z$ did not change the $\mathfrak{q}$-part of $z$ for any place $\mathfrak{q} \neq \idp$. By an analogous procedure for the remaining finitely many places $\idp$ with $\nu_{\idp}(z) < 0$, we conclude that $z$ is indeed represented by an element of $\mathfrak{R}_{\F_q(t)}$. 

Finally, we prove that the system of representatives is indeed minimal. Let $u, v \in \mathfrak{R}_{\F_q(t)}$ with $u \neq v$. We have to show that $u$ and $v$ represent different classes in $\F_q(t)/\wp(\F_q(t))$. Assuming $u - v \in \wp(\F_q(t))$, we find an element $x \in \F_q(t)$ with $u-v = x^p -x$. If we further assume $x \in \F_q$, we are done by the minimality of $\mathfrak{R}_{\F_q}$. Otherwise for $x \in \F_q(t) \setminus \F_q$, we have at least one place $\idp$ with $\nu_{\idp}(x) < 0$. By $\nu_{\idp}(x^p-x) = p \cdot \nu_{\idp}(x)$, we end with a contradiction, since the $\idp$-part of $u-v$ does not have indices divisible by $p$ and $p \nmid \nu_{\idp}(u-v)$ is valid in particular. 

\underline{Proof of b)}: Since $F_{\idp} \cong \F_{q^{\deg(\idp)}}((t))$ is a field of formal Laurent series, $z$ has a representation as a Laurent series expansion of the form $z = \sum\limits_{i = n(\idp)}^{\infty} z_{i, \idp} \cdot \pi_{\idp}^i$. For the principal part and the constant coefficient at $i=0$, we can proceed as in part a). It remains to show that the part $\sum\limits_{i = 1}^{\infty} z_{i, \idp} \cdot \pi_{\idp}^i$ can be eliminated. Let $x_i := - {\sum\limits_{j=0}^{\infty}} (z_{i, \idp} \cdot \pi_{\idp}^i)^{p^j}$. Then, we have
\begin{align*}
\wp(x_i) = - {\sum\limits_{j=0}^{\infty}} (z_{i, \idp} \cdot \pi_{\idp}^i)^{p^{j+1}} + \sum\limits_{j=0}^{\infty} (z_{i, \idp} \cdot \pi_{\idp}^i)^{p^j} = z_{i, \idp} \cdot \pi_{\idp}^i \in \wp(F_{\idp})
\end{align*}
which implies $\wp\left(\sum\limits_{i=1}^{\infty} x_i\right) = \sum\limits_{i = 1}^{\infty} z_{i, \idp} \cdot \pi_{\idp}^i \in \wp(F_{\idp})$ and completes the proof.
\end{proof}

\begin{corollary} \label{VSIso}
	We have the following isomorphisms of $\F_p$-vector spaces
	\begin{enumerate}
		\item[a)]
	\begin{align*}
		\Fqtg/\wp(\Fqtg) \cong \F_q/\wp(\F_q) \oplus \bigoplus\limits_{\idp} \  \Fqtg_{\idp}^{\text{ram}}/\wp\left(\Fqtg_{\idp}^{\text{ram}}\right),
	\end{align*}
\item[b)] \begin{align*}
	F_{\idp}/\wp\left(F_{\idp}\right) \cong \F_q/\wp(\F_q) \oplus F_{\idp}^{\text{ram}}/\wp\left(F_{\idp}^{\text{ram}}\right)
\end{align*}
	\end{enumerate}
	where $\Fqtg_{\idp}^{\text{ram}}/\wp\left(\Fqtg_{\idp}^{\text{ram}}\right)$, respectively $F_{\idp}^{\text{ram}}/\wp\left(F_{\idp}^{\text{ram}}\right)$, denote the respective $\F_p$-subspace encoding the ramification in a fixed prime $\idp$. This means that the extension corresponding to a subspace $U$ of the quotient space on the left-hand side by Proposition \ref{QuotientSpace} is unramified at $\idp$ if and only if the image of $U$ under the canonical projection to the subspace $\Fqtg_{\idp}^{\text{ram}}/\wp\left(\Fqtg_{\idp}^{\text{ram}}\right)$, respectively $F_{\idp}^{\text{ram}}/\wp\left(F_{\idp}^{\text{ram}}\right)$, is trivial.
\end{corollary}

\begin{proof}
	The isomorphisms follow directly by Proposition \ref{RS}. By combining this with Propositions \ref{ASC} and \ref{ConductorChain}, we get the statement about the ramification.
\end{proof}

\begin{remark}
%For local function fields, we obtain an analogous minimal system of representatives, if we restrict ourselves to one place $\idp$. The proof, however, is slightly different in this case, as all non-negative indices for a given Laurent series have to be eliminated. \\
For the global function field $\Fqtg$, the systems of representatives given in Proposition \ref{RS} nearly describe a local-global principle. The only interfering aspect is the coefficient $z_0 \in \mathfrak{R}_{\F_q}$ which occurs globally, just as locally, only once, but not for every place.
\end{remark}

\subsection{The counting functions for $C_p^r$-extensions of $\F_q(t)$ and $F_{\idp}$} \hfill\\
Now, we determine the number of $C_p^r$-extensions of a rational function field $\F_q(t)$ and a local function field $F_{\idp}$ with a given discriminant. As the $C_p^r$-extensions are in bijection to $r$-dimensional $\F_p$-subspaces of $\Fqtg/\wp(\Fqtg)$, respectively $F_{\idp}/\wp(F_{\idp})$, by Proposition \ref{QuotientSpace}, our approach is to count injective $\F_p$-linear mappings $\F_p^r \rightarrow \Fqtg/\wp(\Fqtg)$, respectively $\F_p^r \rightarrow F_{\idp}/\wp(F_{\idp})$, that respect the prescribed discriminant. To simplify the counting, we first drop the requirement of injectivity and carry out a correction afterwards. This correction is performed by an inclusion-exclusion principle for finite abelian groups, which is presented in Lemma \ref{Inclusion-Exclusion} and goes back to Delsarte \cite{De}. 

By Proposition \ref{RS} and Corollary \ref{VSIso}, we have to deal with the unramified part $\F_q/\wp(\F_q)$ and the ramified part $\bigoplus\limits_{\idp} \  \Fqtg_{\idp}^{\text{ram}}/\wp(\Fqtg_{\idp}^{\text{ram}})$, respectively $F_{\idp}^{\text{ram}}/\wp(F_{\idp}^{\text{ram}})$. Note that the ramified part for $\Fqtg$ behaves multiplicatively. 

Additionally, we will first introduce a slightly modified counting function in Definition \ref{CF1}, where we count by the discriminant divisor instead of the norm of the discriminant, as well as some notations for the counting and the inclusion-exclusion principle in Definition \ref{varphifunction}.

\begin{lemma} \label{Inclusion-Exclusion}
Let $A$ and $G$ be finite abelian $p$-groups. Furthermore, let $\Hom(G,A)$ be the set of group homomorphisms and $\Inj(G,A)$ the set of injective group homomorphisms. Then, we have
\begin{align*}
\vert {\Inj(G,A)} \vert = \sum\limits_{pG \trianglelefteq H \trianglelefteq G} \mu\left(G/H\right) \vert {\Hom(H,A)} \vert,
\end{align*}
where the sum ranges over all subgroups lying between $pG$ and $G$ and $\mu$ denotes the multiplicative Möbius function for finite abelian groups, which is defined in \cite{De} as follows:
\begin{align*}
\mu(H) = \begin{cases}
	(-1)^k p^{\frac{k(k-1)}{2}} \ &\text{if} \ H \cong C_p^k \ \text{for some} \ k \geq 0,\\
	0 \ &\text{else.}
\end{cases}
\end{align*}
\end{lemma}

\begin{proof}
All ingredients required for the proof can be found in \cite{De}. 

For a more detailed treatment, see also \cite[Lemma A.1]{La3}, where Lagemann provides a version for counting surjective homomorphisms which we have to dualise for our purposes, since we are not counting surjective but injective linear mappings. In fact, by duality of finite abelian groups, we have a bijection between $\Inj(G,A)$ and $\Surj(\hat{A},\hat{G})$, where $\Surj(\hat{A},\hat{G})$ denotes the set of surjective group homomorphisms from the dual group $\hat{A}$ to the dual group $\hat{G}$. Hence, we can apply the equation on the Möbius inversion in the proof of \cite[Lemma A.1]{La3} to $\Surj(\hat{A},\hat{G})$, which yields
\begin{align*}
\vert {\Inj(G,A)} \vert = \vert {\Surj\left(\hat{A},\hat{G}\right)} \vert = \sum\limits_{p\hat{G} \trianglelefteq \tilde{H} \trianglelefteq \hat{G}} \mu\left(\hat{G}/\tilde{H}\right) \vert {\Hom(\hat{A},\tilde{H})} \vert.
\end{align*}
Since $G$ is finite, there exists a (non-canonical) isomorphism $\iota_G: G \rightarrow \hat{G}$. Hence, for each $\tilde{H} \trianglelefteq \hat{G}$ containing $p\hat{G}$, we obtain a subgroup $H := \iota_G^{-1}(\tilde{H}) \trianglelefteq G$ containing $pG$ such that $\tilde{H} \cong H$. Likewise, $H$ and its dual group $\hat{H}$ are (not canonically) isomorphic which yields $\tilde{H} \cong H \cong \hat{H}$ and, therefore, $\vert {\Hom(\hat{A},\tilde{H})} \vert = \vert {\Hom(\hat{A},\hat{H})} \vert = \vert {\Hom(H,A)} \vert$ by duality. The previous consideration also yields $\mu\left(\hat{G}/\tilde{H}\right) = \mu\left(G/H\right)$. This results in
\begin{align*}
	\vert {\Inj(G,A)} \vert = \sum\limits_{p\hat{G} \trianglelefteq \tilde{H} \trianglelefteq \hat{G}} \mu\left(\hat{G}/\tilde{H}\right) \vert {\Hom(\hat{A},\tilde{H})} \vert = \sum\limits_{pG \trianglelefteq H \trianglelefteq G} \mu\left(G/H\right) \vert {\Hom(H,A)} \vert.
\end{align*}
\end{proof}

\begin{definition} \label{CF1}
Let $F$ be a global function field. $Z(F, C_p^r, \partial)$ denotes the number of $C_p^r$-extensions $K/F$ such that $\partial(K/F) = \partial$ for some prescribed discriminant divisor $\partial$. 

For the local completion $F_{\idp}$ of $F$ at some place $\idp \in \PP_F$, we define $Z(F_{\idp}, C_p^r, \partial)$ as the number of $C_p^r$-extensions $K_{\idp}/F_{\idp}$ such that $\partial(K_{\idp}/F_{\idp}) = \partial$ for some prescribed discriminant divisor $\partial$, which in this case is simply a power of $\idp$.
\end{definition}

\begin{definition} \label{varphifunction}
For any $0 \leq f \leq r$, we define the constant
\begin{align*}
	e_f := \frac{\vert C_p^f \vert}{\vert {\Aut(C_p^r)} \vert} \cdot \binom{r}{f}_p \cdot \mu\left(C_p^{r-f}\right),
\end{align*}
which depends on the group $C_p^r$, the number of subgroups of type $C_p^f$ and the Möbius function $\mu$ for finite abelian groups. 

Additionally, we define the functions
\begin{align*}
	\varphi_f\left(\idp^{n(\idp)}\right) := %\sum\limits_{f=1}^{i} 
	\sum\limits_{h = 0}^{f} \binom{f}{h}_p \quad \sum\limits_{\omega' \in \Omega_h} \gamma_{\omega'} \sum\limits_{\substack{\mathfrak{C} \in \N^h: \\
			\omega_{\mathfrak{C}} = \omega', \\ 
			n(\idp) = d(\mathfrak{C}) . }}
	%	{\substack{\text{chain} \ \mathfrak{C}(\idp)_r \in \N_0^r %\times \{0\}^{r-i} \\
			%			\\ \text{with} \ \omega(\mathfrak{C}(\idp)_i) \ = \ \omega
			%			\\ \text{and} \ J(\mathfrak{C}(\idp)_r) \leq i}} 
	\widehat{Z}_{\idp}(\mathfrak{C}).
\end{align*}
and continue them multiplicatively so that we have $\varphi_f(\partial) = \prod\limits_{\idp^{n(\idp)} \vert \vert {\partial}} \varphi_f\left(\idp^{n(\idp)}\right)$.
Note that the combinatorial constant $\gamma_{\omega'}$ does not depend on the respective place $\idp$ under consideration, $\varphi_f(\mathfrak{1}) = 1$ and $\varphi_0$ is the trivial multiplicative function.
\end{definition}

\begin{proposition} \label{CountingFunction} 
	We have the following decompositions for the counting functions defined in Definition \ref{CF1}:
	\begin{enumerate}
		\item[a)] 
		Let $\partial = \prod\limits_{\idp \in \PP_{\Fqtg}} \idp^{n(\idp)}$ be a fixed discriminant divisor. Then, we have
		\begin{align*}
			Z(\Fqtg, C_p^r, \partial) = \sum\limits_{f=0}^{r} e_f \varphi_f(\partial) = \sum\limits_{f=0}^{r} e_f \prod\limits_{\idp^{n(\idp)} \vert \vert {\partial}} \varphi_f\left(\idp^{n(\idp)}\right).
		\end{align*}
		\item[b)] 
		Let $\partial = \idp^{n(\idp)}$ be a fixed discriminant divisor. Then, we have
		\begin{align*}
			Z(F_{\idp}, C_p^r, \partial) = \sum\limits_{f=0}^{r} e_f \varphi_f(\partial) = \sum\limits_{f=0}^{r} e_f \varphi_f\left(\idp^{n(\idp)}\right).
		\end{align*}
	\end{enumerate}
\end{proposition}

\begin{proof}
%By Corollary \ref{VSIso}, part a) of the Proposition directly implies part b) if we only consider one place $\idp \in \PP_{\Fqtg}$ and the local completion $\Fqtg_{\idp}$ of $\Fqtg$ at $\idp$. Hence, from now on we assume $F = \Fqtg$. \\
\underline{a)}: By Proposition \ref{QuotientSpace}, $\F_p$-subspaces of $\Fqtg/\wp(\Fqtg)$ are in bijection to the Artin-Schreier extensions of $\Fqtg$. Therefore, we may consider $\F_p$-linear mappings $\iota: \F_p^r \rightarrow \Fqtg/\wp(\Fqtg)$ and count them in order to count Artin-Schreier extensions. Since we are only interested in $\F_p$-linear mappings, we will tacitly assume the linearity of the mappings from now on. Hence, we just write $\iota: \F_p^r \rightarrow \Fqtg/\wp(\Fqtg)$ for an $\F_p$-linear map. 

We define the discriminant $\partial(\iota)$, respectively chain $\mathfrak{C}(\iota)$, associated to a map $\iota$ as the discriminant, respectively chain, associated to the extension $\Fqtg(\wp^{-1}(\im(\iota)))/\Fqtg$ corresponding to the image $\im(\iota)$ of $\iota$ in $\Fqtg/\wp(\Fqtg)$ by Proposition \ref{QuotientSpace}.
%Note that we can assign objects such as discriminants or chains, which we have previously defined for extensions $K/\Fqtg$, to a mapping $\iota$ by referring to the respective property of the image of $\iota$ in the quotient space $\Fqtg/\wp(\Fqtg)$. For example, we have $\partial(\iota) = \partial(\im(\iota)) = \partial(F(\wp^{-1}(\im(\iota))))$ where the second equality is due to Proposition \ref{QuotientSpace}. \\
More precisely, $r$-dimensional subspaces of $\Fqtg/\wp(\Fqtg)$ correspond to $C_p^r$-extensions by Proposition \ref{QuotientSpace}. Hence, we have to count injective mappings $\iota$ satisfying $\partial(\iota) = \partial$ and normalise by $\vert {\Aut(\F_p^r)} \vert$ in order to compute $Z(\Fqtg, C_p^r, \partial)$ which yields
\begin{align*}
Z(\Fqtg, C_p^r, \partial) = \frac{1}{\vert {\Aut(\F_p^r)} \vert} \cdot \left\vert\left\{\iota: \F_p^r \rightarrow \Fqtg/\wp(\Fqtg) \ : \ \iota \ \text{injective}, \ \partial(\iota) = \partial \right\}\right\vert.
\end{align*}
Now, we apply Lemma \ref{Inclusion-Exclusion}, with which it is possible to drop the requirement of injectivity. We note that although $\Fqtg/\wp(\Fqtg)$ itself is not finite-dimensional, we can find a suitable finite-dimensional $\F_p$-subspace of $\Fqtg/\wp(\Fqtg)$ that contains all images of mappings $\iota$ satisfying $\partial(\iota) = \partial$. We obtain
\begin{align} \label{DecompDelsarte}
Z(\Fqtg, C_p^r, \partial) = \frac{1}{\vert {\Aut(\F_p^r)} \vert} \sum\limits_{\F_p^f \leq \F_p^r} \mu\left(\F_p^r/\F_p^f \right) \cdot \left\vert\left\{\iota: \F_p^f \rightarrow \Fqtg/\wp(\Fqtg) \ : \ \partial(\iota) = \partial \right\}\right\vert,
\end{align}
where the sum ranges over all $\F_p$-subspaces of $\F_p^r$. By Corollary \ref{VSIso}, we can decompose the mapping $\iota$ as a tuple of mappings where each component refers to one of the direct summands. That is one mapping $\iota_{\F_q}: \F_p^f \rightarrow \F_q/\wp(\F_q)$ for the subspace encoding unramified parts and one mapping $\iota_{\idp}: \F_p^r \rightarrow \Fqtg_{\idp}^{\text{ram}}/\wp\left(\Fqtg_{\idp}^{\text{ram}}\right)$ for each place $\idp \in \PP_{\F_q(t)}$ which encodes the ramification in that place. This results in \eqref{DecompDelsarte}=
\begin{align} \label{DecompDelsarte2}
&\frac{1}{\vert {\Aut(\F_p^r)} \vert} \cdot \sum\limits_{\F_p^f \leq \F_p^r} \mu\left(\F_p^{r-f} \right) \cdot \left\vert\left\{\iota_{\F_q}: \F_p^f \rightarrow \F_q/\wp(\F_q) \right\}\right\vert \notag \\
& \quad \quad \quad \quad \quad \quad \quad \quad \cdot \prod\limits_{\idp \in \PP_{\Fqtg}} \left\vert\left\{\iota_{\idp}: \F_p^f \rightarrow \Fqtg_{\idp}^{\text{ram}}/\wp\left(\Fqtg_{\idp}^{\text{ram}}\right) \ : \ \partial(\iota_{\idp}) = \idp^{n(\idp)} \right\}\right\vert.
\end{align}
Now, we simplify this expression. The number of $f$-dimensional $\F_p$-subspaces in $\F_p^r$ is $\binom{r}{f}_p$ and the number of mappings $\iota_{\F_q}$ is $p^f$ which yields \eqref{DecompDelsarte2}=
\begin{align*}
&\sum\limits_{f=0}^r \frac{\binom{r}{f}_p \cdot \mu\left(\F_p^{r-f}\right) \cdot p^f}{\vert {\Aut(\F_p^r)} \vert} \cdot \prod\limits_{\idp \in \PP_{\Fqtg}} \left\vert\left\{\iota_{\idp}: \F_p^f \rightarrow \Fqtg_{\idp}^{\text{ram}}/\wp\left(\Fqtg_{\idp}^{\text{ram}}\right) \ : \ \partial(\iota_{\idp}) = \idp^{n(\idp)} \right\}\right\vert \\
= & \sum\limits_{f=0}^r e_f \cdot \prod\limits_{\idp \in \PP_{\Fqtg}} \left\vert\left\{\iota_{\idp}: \F_p^f \rightarrow \Fqtg_{\idp}^{\text{ram}}/\wp\left(\Fqtg_{\idp}^{\text{ram}}\right) \ : \ \partial(\iota_{\idp}) = \idp^{n(\idp)} \right\}\right\vert.
\end{align*}
Hence, it remains to prove the following local equation for each place $\idp \in \PP_{\Fqtg}$:
\begin{align} \label{LocalEq}
\left\vert\left\{\iota_{\idp}: \F_p^f \rightarrow \Fqtg_{\idp}^{\text{ram}}/\wp\left(\Fqtg_{\idp}^{\text{ram}}\right) \ : \ \partial(\iota_{\idp}) = \idp^{n(\idp)} \right\}\right\vert = \varphi_f\left(\idp^{n(\idp)}\right).
\end{align}

By Proposition \ref{ConductorChain} and Definition \ref{ChainComp}, each mapping $\iota_{\idp}$ has a chain $\mathfrak{C}(\iota_{\idp}) \in \N^h$ satisfying $d(\mathfrak{C}(\iota_{\idp})) = n(\idp)$ and an associated composition $\omega_{\mathfrak{C}(\iota_{\idp})} \in \Omega_h$ for some $0 \leq h \leq f$. We recall the formula for the discriminant of a chain given in Lemma \ref{DiscExponent}. 

If we only consider chains $\mathfrak{C} \in \N^h$ for any $0 \leq h \leq f$ satisfying $d(\mathfrak{C}) = n(\idp)$, then the sets $\left\{\iota_{\idp}: \F_p^f \rightarrow \Fqtg_{\idp}^{\text{ram}}/\wp\left(\Fqtg_{\idp}^{\text{ram}}\right) \ : \ \mathfrak{C}(\iota_{\idp}) = \mathfrak{C} \right\}$ form a partition of the set $\left\{\iota_{\idp}: \F_p^f \rightarrow \Fqtg_{\idp}^{\text{ram}}/\wp\left(\Fqtg_{\idp}^{\text{ram}}\right) \ : \ \partial(\iota_{\idp}) = \idp^{n(\idp)} \right\}$,
which yields
\begin{align} \label{partition}
& \left\vert\left\{\iota_{\idp}: \F_p^f \rightarrow \Fqtg_{\idp}^{\text{ram}}/\wp\left(\Fqtg_{\idp}^{\text{ram}}\right) \ : \ \partial(\iota_{\idp}) = \idp^{n(\idp)} \right\}\right\vert \notag \\
= \quad &\sum\limits_{h = 0}^{f} \sum\limits_{\substack{\mathfrak{C} \in \N^h: \\
		n(\idp) = d(\mathfrak{C})  }}  \left\vert\left\{\iota_{\idp}: \F_p^f \rightarrow \Fqtg_{\idp}^{\text{ram}}/\wp\left(\Fqtg_{\idp}^{\text{ram}}\right) \ : \ \mathfrak{C}(\iota_{\idp}) = \mathfrak{C} \right\}\right\vert. 
%= \quad &\sum\limits_{h = 0}^{f} \sum\limits_{\omega' \in \Omega_h} \sum\limits_{\substack{\mathfrak{C} \in \N^h: \\
%		\omega_{\mathfrak{C}} = \omega', \\ 
%		n(\idp) = d(\mathfrak{C}) . }} 
%	\left\vert\left\{\iota_{\idp}: \F_p^f \rightarrow \Fqtg_{\idp}^{\text{ram}}/\wp\left(\Fqtg_{\idp}^{\text{ram}}\right) \ : \  \mathfrak{C}(\iota_{\idp}) = \mathfrak{C}, \ \omega_{\mathfrak{C}(\iota_{\idp})} = \omega' \right\}\right\vert.
\end{align}
Therefore, it remains to prove the identity 
\begin{align} \label{inductionhypo}
\left\vert\left\{\iota_{\idp}: \F_p^f \rightarrow \Fqtg_{\idp}^{\text{ram}}/\wp\left(\Fqtg_{\idp}^{\text{ram}}\right) \ : \ \mathfrak{C}(\iota_{\idp}) = \mathfrak{C} \right\}\right\vert = \binom{f}{h}_p \cdot \gamma_{\omega_{\mathfrak{C}}} \cdot \widehat{Z}_{\idp}(\mathfrak{C})
\end{align}
since equation \eqref{inductionhypo} yields \eqref{partition}=
\begin{align*}
\sum\limits_{h = 0}^{f} \sum\limits_{\substack{\mathfrak{C} \in \N^h: \\
		n(\idp) = d(\mathfrak{C})  }}  \binom{f}{h}_p \cdot \gamma_{\omega_{\mathfrak{C}}} \cdot \widehat{Z}_{\idp}(\mathfrak{C}) = %\sum\limits_{f=1}^{i} 
	\sum\limits_{h = 0}^{f} \binom{f}{h}_p \quad \sum\limits_{\omega' \in \Omega_h} \gamma_{\omega'} \sum\limits_{\substack{\mathfrak{C} \in \N^h: \\
	\omega_{\mathfrak{C}} = \omega', \\ 
	n(\idp) = d(\mathfrak{C}) . }} \widehat{Z}_{\idp}(\mathfrak{C}) = \varphi_f\left(\idp^{n(\idp)}\right).
\end{align*}

Now, let $\mathfrak{C} = (c_1, \ldots, c_h) \in \N^h$ be a fixed chain of length $h$ with the associated composition $\omega_{\mathfrak{C}} = \left(a_1, \ldots, a_{\lambda_{\omega_{\mathfrak{C}}}}\right) \in \Omega_h$. Any map $\iota_{\idp}: \F_p^f \rightarrow \Fqtg_{\idp}^{\text{ram}}/\wp\left(\Fqtg_{\idp}^{\text{ram}}\right)$ has a unique factorisation through $\F_p^f/\ker(\iota_{\idp})$ by the fundamental homomorphism theorem. For a chain of length $h$, we have $\dim_{\F_p}(\im(\iota_{\idp})) = h$ which yields $\F_p^f/\ker(\iota_{\idp}) \cong \F_p^h$. Hence, we obtain
\begin{align*}
& \left\vert\left\{\iota_{\idp}: \F_p^f \rightarrow \Fqtg_{\idp}^{\text{ram}}/\wp\left(\Fqtg_{\idp}^{\text{ram}}\right) \ : \ \mathfrak{C}(\iota_{\idp}) = \mathfrak{C} \right\}\right\vert \\
= \quad & \sum\limits_{\substack{V \leq \F_p^f: \\ \dim_{\F_p}(V) = f-h }} \left\vert\left\{\overline{\iota_{\idp}}: \F_p^f/V \rightarrow \Fqtg_{\idp}^{\text{ram}}/\wp\left(\Fqtg_{\idp}^{\text{ram}}\right) \ : \ \overline{\iota_{\idp}} \ \text{injective}, \mathfrak{C}(\overline{\iota_{\idp}}) = \mathfrak{C} \right\}\right\vert \\
= \quad & \binom{f}{h}_p \cdot \left\vert\left\{\iota_{\idp}: \F_p^h \rightarrow \Fqtg_{\idp}^{\text{ram}}/\wp\left(\Fqtg_{\idp}^{\text{ram}}\right) \ : \ \mathfrak{C}(\iota_{\idp}) = \mathfrak{C} \right\}\right\vert.
\end{align*}
%a
%\begin{align} \label{fact}
%&\left\vert\left\{\iota_{\idp}: \F_p^f \rightarrow \Fqtg_{\idp}^{\text{ram}}/\wp\left(\Fqtg_{\idp}^{\text{ram}}\right) \ : \ \mathfrak{C}(\iota_{\idp}) = \mathfrak{C} \right\}\right\vert \notag \\
%= \quad & \left\vert\left\{\pi: \F_p^f \rightarrow F_p^f/\ker(\iota_{\idp}) \right\}\right\vert \cdot \left\vert\left\{\overline{\iota_{\idp}}: F_p^f/\ker(\iota_{\idp}) \rightarrow \Fqtg_{\idp}^{\text{ram}}/\wp\left(\Fqtg_{\idp}^{\text{ram}}\right) \ : \ \mathfrak{C}(\overline{\iota_{\idp}}) = \mathfrak{C} \right\}\right\vert.
%\end{align}
%For a chain of length $h$, we have $\dim_{\F_p}(\im(\iota_{\idp})) = h$ which yields $\F_p^f/\ker(\iota_{\idp}) \cong \F_p^h$ such that \eqref{fact} = 
%\begin{align*}
%\binom{f}{h}_p \cdot \left\vert\left\{\iota_{\idp}: \F_p^h \rightarrow \Fqtg_{\idp}^{\text{ram}}/\wp\left(\Fqtg_{\idp}^{\text{ram}}\right) \ : \ \mathfrak{C}(\iota_{\idp}) = \mathfrak{C} \right\}\right\vert.
%\end{align*}

Thus, it remains to prove 
\begin{align} \label{induchypo2}
\left\vert\left\{\iota_{\idp}: \F_p^h \rightarrow \Fqtg_{\idp}^{\text{ram}}/\wp\left(\Fqtg_{\idp}^{\text{ram}}\right) \ : \ \mathfrak{C}(\iota_{\idp}) = \mathfrak{C} \right\}\right\vert = \gamma_{\omega_{\mathfrak{C}}} \cdot \widehat{Z}_{\idp}(\mathfrak{C})
\end{align}
for any $0 \leq h \leq f$ and all chains $\mathfrak{C}$ of length $h$ which we do by induction on $\lambda$, the number of blocks of the associated composition $\omega_{\mathfrak{C}}$. The case $h=0$ is clear. Hence, we fix some $h \geq 1$ and our induction hypothesis is that \eqref{induchypo2} holds for all chains $\mathfrak{C}$ with $\lambda_{\omega_{\mathfrak{C}}} < \lambda$. 

For the initial step, we consider $\lambda = 1$. That is $c_1 = \ldots = c_h$ and $\omega_{\mathfrak{C}} = (h)$. By Proposition \ref{RS}, we have to choose $h$ elements $\overline{z_1}, \ldots, \overline{z_h} \in \mathfrak{R}_{\Fqtg}$, the images of the standard basis vectors of $\F_p^h$, such that $\ASC(\idp, \overline{z_j}) = c_1 - 1$ for all $1 \leq j \leq h$ and the leading coefficients of the $\idp$-parts of the $\overline{z_j}$ are $\F_p$-linearly independent. Let $\sum\limits_{\substack{i=1 \\ p \ \nmid \ i}}^{\ASC(\idp, \overline{z_j})} \frac{z_{j,i,\idp}}{\pi_{\idp}^{i}}$ be the $\idp$-part of $\overline{z_j}$ which has precisely $r(c_1)$ summands. As we have $z_{j,i,\idp} \in \F_{q^{\deg(\idp)}}$, we have $q^{\deg(\idp)} = \np$ possibilities for each coefficient which results in $\np^{r(c_1)}$ possibilities for the $\idp$-part of each $\overline{z_j}$. However, we still have to exclude the linear dependency of the leading coefficients. For $j=1$, we only have to satisfy the condition $z_{1, c_1-1,\idp} \neq 0$ which yields $\np^{r(c_1)} - \np^{r(c_1 - 1)} = \widehat{Z}_{\idp,1}(c_1)$ possibilities for $\overline{z_1}$ where we recall Definition \ref{DefinitionCountingFunction}. For $1 < j \leq h$, we have to satisfy the condition $z_{j, c_1-1,\idp} \not \in \langle z_{1, c_1-1,\idp}, \ldots, z_{j-1, c_1-1,\idp} \rangle_{\F_p}$ which yields $\np^{r(c_1)} - p^{j-1 }\np^{r(c_1 - 1)} = \widehat{Z}_{\idp,p^{j-1}}(c_1)$. Therefore, we obtain
\begin{align*}
\widehat{Z}_{\idp}(\mathfrak{C}) = \prod\limits_{j=0}^{h-1} \widehat{Z}_{\idp, p^j}(c_1) = \prod\limits_{i=1}^{\lambda} \prod\limits_{j=0}^{a_i-1} \widehat{Z}_{\idp, p^j}(c_{A_i})
\end{align*}
for $\mathfrak{C}$ and we are done with the case $\lambda = 1$.

Now, we assume $\lambda \geq 2$. Let $\mathfrak{C}' = (c_{a_1 + 1}, \ldots, c_f) \in \N^{f-a_1}$ be the subchain of $\mathfrak{C}$ which we obtain by removing the components belonging to the first block $a_1$ of $\omega_{\mathfrak{C}}$. The chain $\mathfrak{C}'$ has the composition $\omega_{\mathfrak{C}'} = \left(a_2, \ldots, a_{\lambda_{\omega_{\mathfrak{C}}}}\right)$. Note that $\mathfrak{C}'$ is not the empty chain by our assumption. Furthermore, let $\mathfrak{C}'' = (c_1, \ldots, c_{a_1})$ denote the subchain of $\mathfrak{C}$ corresponding to the first block with $\omega_{\mathfrak{C}''} = (a_1)$. 

For each mapping $\iota_{\idp}: \F_p^f \rightarrow \Fqtg_{\idp}^{\text{ram}}/\wp\left(\Fqtg_{\idp}^{\text{ram}}\right)$ such that $\mathfrak{C}(\iota_{\idp}) = \mathfrak{C}$, there is a uniquely determined $(f-a_1)$-dimensional $\F_p$-subspace $V_{\mathfrak{C}',\iota_{\idp}}$ of the image of $\iota_{\idp}$ which has the chain $\mathfrak{C}'$. The uniqueness is due to the fact that each other $(f-a_1)$-dimensional $\F_p$-subspace of the image of $\iota_{\idp}$ contains at least one $1$-dimensional $\F_p$-subspace with conductor exponent $c_1$. Let $V_{\iota_{\idp}} := \iota^{-1}_{\idp}(V_{\mathfrak{C}',\iota_{\idp}}) \cong \F_p^{f-a_1} \leq \F_p^f $ denote the preimage and let $W \cong \F_p^{a_1}$ be an arbitrary complement of $V_{\iota_{\idp}}$ in $\F_p^f$. Then, we have $\F_p^f \cong V_{\iota_{\idp}} \oplus W$ and the mapping $\iota_{\idp}$ is uniquely described by the mappings ${\iota_{\idp}}{\vert_{V_{\iota_{\idp}}}}$ and ${\iota_{\idp}}{\vert_{W}}$. We obtain
\begin{align} \label{Comp2}
&\left\vert\left\{\iota_{\idp}: \F_p^f \rightarrow \Fqtg_{\idp}^{\text{ram}}/\wp\left(\Fqtg_{\idp}^{\text{ram}}\right) \ : \ \mathfrak{C}(\iota_{\idp}) = \mathfrak{C} \right\}\right\vert \notag \\
= \quad &\sum\limits_{\substack{V \leq \F_p^f: \\ \dim_{\F_p}(V) = f-a_1}} \left\vert\left\{\iota_{\idp}: \F_p^f \rightarrow \Fqtg_{\idp}^{\text{ram}}/\wp\left(\Fqtg_{\idp}^{\text{ram}}\right) \ : \ \mathfrak{C}(\iota_{\idp}) = \mathfrak{C}, \ V_{\iota_{\idp}} = V \right\}\right\vert \notag \\
%= \quad &\sum\limits_{\substack{V \leq \F_p^f: \\ \dim_{\F_p}(V) = f-a_1}} \left\vert\left\{ {\iota_{\idp}}{\vert_{V}}: V \rightarrow \Fqtg_{\idp}^{\text{ram}}/\wp\left(\Fqtg_{\idp}^{\text{ram}}\right) \ : \ \mathfrak{C}({\iota_{\idp}}{\vert_{V}}) = \mathfrak{C}'\right\}\right\vert \\
%& \quad \quad \quad \quad \quad \quad \quad \quad \cdot \left\vert\left\{ {\iota_{\idp}}{\vert_{W}}: W \rightarrow \Fqtg_{\idp}^{\text{ram}}/\wp\left(\Fqtg_{\idp}^{\text{ram}}\right) \ : \ \mathfrak{C}({\iota_{\idp}}{\vert_{W}}) = \mathfrak{C}''\right\}\right\vert \\
= \quad & \binom{f}{f-a_1}_p \cdot \left\vert\left\{\iota_{\idp}\vert_{V}: \F_p^{f-a_1} \rightarrow \Fqtg_{\idp}^{\text{ram}}/\wp\left(\Fqtg_{\idp}^{\text{ram}}\right) \ : \ \mathfrak{C}(\iota_{\idp}\vert_{V}) = \mathfrak{C}' \right\}\right\vert \notag \\
\quad \quad \quad \quad & \cdot \left\vert\left\{\iota_{\idp}\vert_{W}: \F_p^{a_1} \rightarrow \Fqtg_{\idp}^{\text{ram}}/\wp\left(\Fqtg_{\idp}^{\text{ram}}\right) \ : \ \mathfrak{C}(\iota_{\idp}\vert_{W}) = \mathfrak{C}'' \right\}\right\vert.
\end{align}
By our induction hypothesis, we have \eqref{Comp2}=
\begin{align*}
\binom{f}{f-a_1}_p \cdot \gamma_{\omega_{\mathfrak{C}'}} \cdot \widehat{Z}_{\idp}(\mathfrak{C}') \cdot \gamma_{\omega_{\mathfrak{C}''}} \cdot \widehat{Z}_{\idp}(\mathfrak{C}'') = \binom{f}{f-a_1}_p \cdot \gamma_{\omega_{\mathfrak{C}'}} \cdot \widehat{Z}_{\idp}(\mathfrak{C})
\end{align*}
as we have $\gamma_{\omega_{\mathfrak{C}''}} = 1$ by Definition \ref{Composition} and $\widehat{Z}_{\idp}(\mathfrak{C}') \cdot \widehat{Z}_{\idp}(\mathfrak{C}'') = \widehat{Z}_{\idp}(\mathfrak{C})$ by Definition \ref{DefinitionCountingFunction}. Furthermore, we have $\binom{f}{f-a_1}_p \cdot \gamma_{\omega_{\mathfrak{C}'}} = \gamma_{\sigma(\omega_{\mathfrak{C}})}$ for $\sigma(\omega_{\mathfrak{C}}) = (a_2, \ldots, a_{\lambda}, a_1)$ which yields $\binom{f}{f-a_1}_p \cdot \gamma_{\omega_{\mathfrak{C}'}} = \gamma_{\omega_{\mathfrak{C}}}$ by Lemma \ref{gammainvariance}. 

\underline{b)}: The derivation of equation \eqref{DecompDelsarte2} also applies to a local field $F_{\idp}$ as well as the local equation \eqref{LocalEq} and its proof.
\end{proof}

\newpage

\subsection{The corresponding Dirichlet series for the base field $\F_q(t)$} \hfill \\

\begin{definition} \label{DefinitionDirichletSeries}
For $0 \leq f \leq r$ and a local or global function field $F$, we define the Dirichlet series
\begin{align*}
	\Phi_f(F, C_p^r;s) := \sum\limits_{\partial} \varphi_f(\partial) \cdot \vert \vert {\partial} \vert \vert^{-s} = \prod\limits\limits_{\idp \in \mathbb{P}_F} \left( 1 + \sum\limits\limits_{n(\idp) \geq 1} \varphi_f\left(\idp^{n(\idp)} \right) \cdot \np^{-n(\idp)s} \right)
\end{align*}
corresponding to the function $\varphi_f(\partial)$ defined in Definition \ref{varphifunction}. Note that the representation as an Euler product follows by the multiplicativity of $\varphi(\partial)$ and the norm function. We further denote the respective Euler factors by $\Phi_{f, \idp}(F, C_p^r;s)$.
\end{definition}

We remark that we have defined the Dirichlet series for an arbitrary base field $F$ with foresight regarding Section \ref{ArbitraryBaseFields}. At present, however, we are still only interested in the special case $F = \Fqtg$. 

As a next step, we investigate the Dirichlet series $\Phi(\Fqtg,C_p^r;s) = \sum\limits_{\partial} Z(\Fqtg, C_p^r, \partial) \cdot \vert \vert {\partial} \vert \vert^{-s}$ corresponding to the counting function $Z(\Fqtg, C_p^r, \partial)$. 

\begin{proposition} \label{LCEuler} %Dirichlet series as linear combination of Euler products
For the Dirichlet series $\Phi(\Fqtg, C_p^r;s)$, %= \sum\limits_{\partial} Z(\Fqtg, C_p^r, \partial) \cdot \vert \vert {\partial} \vert \vert^{-s}$ 
we obtain the decomposition
\begin{align*}
	\Phi(\Fqtg, C_p^r;s) = \sum\limits_{f=0}^{r} e_f \cdot \Phi_f(\Fqtg, C_p^r;s).
\end{align*}
\end{proposition}

\begin{proof}
This is a direct consequence of Proposition \ref{CountingFunction} and the decomposition of $Z(\Fqtg, C_p^r, \partial)$ presented there.
\end{proof}

\subsection{Meromomorphic continuation of the Euler factors} \label{LocalContinuation} \hfill \\
Henceforth, we consider the Euler factors and compute meromorphic continuations of them which will finally lead us to a meromorphic continuation of the Dirichlet series $\Phi(\Fqtg, C_p^r;s)$. With regard to Section \ref{ArbitraryBaseFields}, however, we make sure in the following calculations that we do not use any arithmetic properties that are special for the base field $\Fqtg$. 

Since we are initially only considering the Euler factors, we suppress $\idp$ wherever it is possible in order to shorten the notation.
By the definition of $\varphi_f\left(\idp^{n(\idp)}\right)$ given in Definition \ref{varphifunction}, we get
\begin{align} \label{EqEuler}
	\Phi_{f, \idp}(F, C_p^r;s) &= 1 + \sum\limits_{n \geq 1}	\varphi_f\left(\idp^{n}\right) \cdot \np^{-ns} \nonumber \\
	&= 1 + \sum\limits_{n \geq 1} \sum\limits_{h = 1}^{f} \binom{f}{h}_p \quad \sum\limits_{\omega' \in \Omega_h} \gamma_{\omega'} \sum\limits_{\substack{\mathfrak{C} \in \N^h: \\
			\omega_{\mathfrak{C}} = \omega', \\ 
			n = d(\mathfrak{C}) . }}
	%	{\substack{\text{chain} \ \mathfrak{C}_r \in \N_0^r %\times \{0\}^{r-i} \\
	%			\\ \text{with} \ \omega(\mathfrak{C}_i) \ = \ \omega
	%			\\ \text{and} \ J(\mathfrak{C}_r) \leq i}} 
	\widehat{Z}_{\idp}(\mathfrak{C}) \cdot \np^{-ns} \nonumber \\
	&= 1 + \sum\limits_{h = 1}^{f} \binom{f}{h}_p \quad \sum\limits_{\omega' \in \Omega_h} \gamma_{\omega'} \sum\limits_{n \geq 1} \sum\limits_{\substack{\mathfrak{C} \in \N^h: \\
			\omega_{\mathfrak{C}} = \omega', \\ 
			n = d(\mathfrak{C}) . }}
	%	{\substack{\text{chain} \ \mathfrak{C}_r \in \N_0^r %\times \{0\}^{r-i} \\
	%			\\ \text{with} \ \omega(\mathfrak{C}_i) \ = \ \omega
	%			\\ \text{and} \ J(\mathfrak{C}_r) \leq i}} 
	\widehat{Z}_{\idp}(\mathfrak{C}) \cdot \np^{-ns}.
\end{align}
%Certainly, we can swap the summations over the compositions and $n$ in equation \ref{EqEuler}. \\

We now fix a composition $\omega' \in \Omega_h$ and we only consider chains satisfying the condition $\omega_{\mathfrak{C}} = \omega'$. To avoid clutter, we also omit the dependence on $\mathfrak{C}$ in our notation and we abbreviate $\omega_{\mathfrak{C}}$ by $\omega$ %$\lambda(\omega(\mathfrak{C}))$ by $\lambda$ 
for instance. 
First of all, we aim to simplify the term $\widehat{Z}_{\idp}(\omega)$. %We consider a generic composition $\omega_i \in \Omega_i$ and carry out the computations in general. \\ \\
Since we have $r(c) = c-1 - \left\lfloor \frac{c-1}{p} \right\rfloor$ by Definition \ref{DefinitionCountingFunction}, we represent $c_i$ as $c_i = pk_i + \ell_i +1$ with $k_i \geq 0$ and $1 \leq \ell_i \leq p-1$. This representation allows us to explicitly compute the value of the floor function:
\begin{align*}
r(c_i) = (p-1)k_i + \ell_i.
\end{align*}
We remark that the value $\ell_i = 0$ is not admitted because the Artin-Schreier conductor defined in Proposition \ref{ASC} is not divisible by $p$. Therefore, the conductor exponents are not congruent to $1 \bmod p$ by Proposition \ref{ASC}(c). Indeed, consistently the formula yields $\widehat{Z}_{\idp, 1}(c) = \np^{r(c)} - \np^{r(c-1)} = 0$ in the case $c \equiv 1 \mod p$. 

By the assumptions on the chain $\mathfrak{C}$, we have $c_{A_i} > c_{A_i+1} = \ldots = c_{A_{i+1}}$ for all indices $1 \leq i \leq \lambda -1$. With regard to $c_i = pk_i + \ell_i + 1$, we have to distinguish two cases for each index to fulfil the proposed inequalities. On the one hand we can have $k_{A_i} > k_{A_{i+1}}$ so that $1 \leq \ell_{A_i}, \ell_{A_{i+1}} \leq p-1$ are independent of one another and on the other hand we can have $k_{A_i} = k_{A_{i+1}}$ with the restriction $1 \leq \ell_{A_{i+1}} < \ell_{A_i} \leq p-1$. 

Next, we want to refine the composition $\omega_{\mathfrak{C}}$ assigned to the chain $\mathfrak{C}$ in such a way that we can also describe these dependencies of the $k$- and $\ell$-values. This leads us to the concept of \emph{2-level compositions}.

\begin{definition} \label{2-level-composition}
	Let $\rho = (b_1, \ldots, b_{\lambda_{\rho}}) \in \Omega_h$ be a composition. A \emph{2-level composition of $h$} is a tuple of compositions 
	\begin{align*}
	\theta := (\theta_{1}, \ldots, \theta_{\lambda_{\rho}}) = ((a_{1}, \ldots, a_{\Lambda_{1}}), \ldots, (a_{1 + \Lambda_{{\lambda_{\rho}-1}}}, \ldots, a_{\Lambda_{\lambda_{\rho}}})) 
	\end{align*}
where $\theta_i := (a_{1 + \Lambda_{i-1}}, \ldots, a_{\Lambda_i}) \in \Omega_{b_i}$ is a composition belonging to the $b_i$-block of $\rho$ for all $1 \leq i \leq \lambda_{\rho}$. As before, we denote the lengths of these compositions by $\lambda_{\theta_i}$ and we define $\Lambda_{i} := \sum\limits_{j=1}^{i} \lambda_{\theta_j}$. We call the $b_i$ \emph{outer blocks}, the $a_i$ \emph{inner blocks}, the $\theta_i$ \emph{inner compositions} and $\rho$ the \emph{outer composition} of $\theta$.
Let further $\theta_{\omega} := (a_{1}, \ldots, a_{\Lambda_{\lambda_{\rho}}})$ be the composition with all inner blocks of $\theta$ as its blocks. 

Furthermore, we define the monic polynomial $g_{\theta}(X) := \prod\limits_{i=1}^{\Lambda_{\lambda_{\rho}}} \prod\limits_{j=0}^{a_i-1} \left(X - p^j \right) \in \Z[X]$ of degree $h$ which we call the \emph{structure polynomial of $\theta$}. We note $g_{\theta}(\np) \neq 0$ if and only if $a_i \leq n \cdot \deg(\idp)$ for all $1 \leq i \leq \Lambda_{\lambda_{\rho}}$ where we recall $\np = q^{\deg(\idp)} = p^{n \cdot \deg(\idp)}$. 

We call a 2-level composition \emph{admissible} if the length of all inner compositions is at most $p-1$. That is $\lambda_{\theta_i} \leq p-1$ for all $1 \leq i \leq \lambda_{\rho}$. 

Finally, $\mathfrak{Z}_h$ denotes the set of 2-level compositions of $h$ and $\mathfrak{Z}^{ad}_h$ denotes the set of admissible 2-level compositions of $h$. We note that the set $\mathfrak{Z}_0$ only contains the empty 2-level composition $(())$. For $h \geq 1$, we have $\vert {\mathfrak{Z}_h} \vert = 3^{h-1}$. This can be proven in a similar way to the classical fact about compositions mentioned in Definition \ref{Composition}.
\end{definition}

To every chain $\mathfrak{C} = (c_1, \ldots, c_h) \in \N^h$ such that $c_i \not \equiv 1 \bmod p$ for $1 \leq i \leq h$, we can assign a uniquely determined admissible 2-level composition of $h$, which we call $\theta_{\mathfrak{C}}$, in the following way: We consider the tuple $(k_1, \ldots, k_h)$ of $k$-values of $\mathfrak{C}$ and associate a composition to this tuple as we did in Definition \ref{ChainComp} for the chains. However, we note that the $k$-values can also have the value 0 and that we also take this value into account when we associate the composition. The composition obtained in this way is the outer composition $\rho = (b_1, \ldots, b_{\lambda_{\rho}}) \in \Omega_h$. 

Now, we can consider an outer block $b_i$ of the outer composition. By construction, all $k$-values of components of $\mathfrak{C}$ belonging to this block coincide. We further consider the $\ell$-values of these components of $\mathfrak{C}$ which form a chain $\mathfrak{C}_{b_i}$ and we associate a composition to $\mathfrak{C}_{b_i}$ as in Definition \ref{ChainComp}. This composition is the inner composition $\theta_i$. 

In fact, the 2-level composition $\theta_{\mathfrak{C}}$ created in this way is automatically admissible as the $\ell$-values are integral and bounded by $1$ and $p-1$.
Additionally, we can recover the composition $\omega_{\mathfrak{C}}$ originally associated to $\mathfrak{C}$ by considering the composition with all inner blocks of $\theta_{\mathfrak{C}}$ as its blocks. This means $\omega_{\mathfrak{C}} = \theta_{\mathfrak{C}_{\omega}}$.

\begin{example}
Let us consider a chain $\mathfrak{C} \in \N^4$ satisfying $c_1 = c_2 > c_3 > c_4 > 0$, which yields $\omega_{\mathfrak{C}} = (2,1,1) \in \Omega_4$. If we further have $k_2 > k_3$ with $1 \leq \ell_2, \ell_3 \leq p-1$ and $k_3 = k_4$ with $1 \leq \ell_4 < \ell_3 \leq p-1$ (assuming $p \geq 3$), this setting corresponds to the admissible 2-level composition $((2),(1,1)) \in \mathfrak{Z}^{ad}_4$. Note that in the case $p=2$, neither the previous configuration of the $k$- and $\ell$-values nor the 2-level composition $((2),(1,1))$ are admissible. If instead we have $k_3 > k_4$, we obtain the 2-level composition $((2),(1),(1)) \in \mathfrak{Z}^{ad}_4$.
\end{example}

We apply the notation provided in Definition \ref{2-level-composition} and obtain \eqref{EqEuler}=

\begin{align} \label{EinsetzenConfig}
1 + \sum\limits_{h = 1}^{f} \binom{f}{h}_p \quad \sum\limits_{\theta \in \mathfrak{Z}^{ad}_h} \gamma_{\theta_{\omega}} \sum\limits_{n \geq 1} \sum\limits_{\substack{\mathfrak{C} \in \N^h: \\
		\omega_{\mathfrak{C}} = \theta_{\omega}, \\ 
		n = d(\mathfrak{C}) . }}
%	{\substack{\text{chain} \ \mathfrak{C}_r \in \N_0^r %\times \{0\}^{r-i} \\
%			\\ \text{with} \ \omega(\mathfrak{C}_i) \ = \ \omega
%			\\ \text{and} \ J(\mathfrak{C}_r) \leq i}} 
\widehat{Z}_{\idp}(\mathfrak{C}) \cdot \np^{-ns}.
\end{align}

From now on, we choose and fix any $\theta \in \mathfrak{Z}^{ad}_h$ and carry out the computations for this $\theta$. Using the notation provided in Definition \ref{2-level-composition}, %In particular, $\omega = (a_1, \ldots, a_{\Lambda_{\lambda_{\rho}}}) \in \Omega_h$ is the corresponding composition of the associated $h$-configuration. 
we compute

\begin{align} \label{SideCalculation1}
\widehat{Z}_{\idp}(\mathfrak{C}) &= \prod\limits_{i=1}^{\Lambda_{\lambda_{\rho}}} \prod\limits_{j=0}^{a_i-1} \widehat{Z}_{\idp, p^j}(c_{A_i}) \nonumber \\ 
 &= \prod\limits_{i=1}^{\Lambda_{\lambda_{\rho}}} \prod\limits_{j=0}^{a_i-1} \left(\np^{c_{A_i}-1 - \left\lfloor \frac{c_{A_i}-1}{p} \right\rfloor} - p^j \np^{c_{A_i}-2 - \left\lfloor \frac{c_{A_i}-2}{p} \right\rfloor} \right) \nonumber \\
&= \prod\limits_{i=1}^{\Lambda_{\lambda_{\rho}}} \prod\limits_{j=0}^{a_i-1} \left(\np^{(p-1)k_{A_i} + \ell_{A_i}} - p^j \np^{(p-1)k_{A_i} + \ell_{A_i} -1} \right) \nonumber \\
&= \left(\prod\limits_{i=1}^{\Lambda_{\lambda_{\rho}}} \prod\limits_{j=0}^{a_i-1} \left(\np - p^j \right)\right) \cdot \left(\prod\limits_{i=1}^{\Lambda_{\lambda_{\rho}}} \np^{a_i(p-1)k_{A_i}} \cdot \np^{a_i(\ell_{A_i}-1)}\right) \nonumber \\
&= g_{\theta}(\np) \cdot \prod_{i=1}^{\lambda_{\rho}} \np^{\sum\limits_{j=1 + \Lambda_{i-1}}^{\Lambda_i} a_j (p-1)k_{\Lambda_i}} \cdot \prod\limits_{i=1}^{\Lambda_{\lambda_{\rho}}} \np^{a_i(\ell_{A_i}-1)} \nonumber \\
&= g_{\theta}(\np) \cdot \prod_{i=1}^{\lambda_{\rho}} \np^{b_i (p-1) k_{\Lambda_i}} \cdot \prod\limits_{i=1}^{\Lambda_{\lambda_{\rho}}} \np^{a_i(\ell_{A_i}-1)}.
\end{align}
%Transform and simplify the Euler products by $c_i(\idp) = pk_i + \ell_i$...
%First, we consider the case $k_i > k_{i+1}$ for all $1 \leq i \leq \lambda -1$ in which we have a summation over $k_i$ and $\ell_i$ for every block of the corresponding composition and we do not have any further restrictions to the $\ell_i$. We call this case generic for our composition and we obtain
Now, we consider the sub-series in equation \eqref{EinsetzenConfig} arising from $\theta \in \mathfrak{Z}^{ad}_h$. We plug in equation \eqref{SideCalculation1} and Proposition \ref{DiscExponent} and obtain

\begin{align} \label{eq1}
	&\sum\limits_{n \geq 1} \sum\limits_{\substack{\mathfrak{C} \in \N^h: \\
			\omega_{\mathfrak{C}} = \theta_{\omega}, \\ 
			n = d(\mathfrak{C}) . }}
	%	{\substack{\text{chain} \ \mathfrak{C}_r \in \N_0^r %\times \{0\}^{r-i} \\
	%			\\ \text{with} \ \omega(\mathfrak{C}_i) \ = \ \omega
	%			\\ \text{and} \ J(\mathfrak{C}_r) \leq i}} 
	\widehat{Z}_{\idp}(\mathfrak{C}) \cdot \np^{-ns} \nonumber \\
= \quad &\sum\limits_{k_{\Lambda_1} \geq 0} \sum\limits_{{k_{\Lambda_2}} = 0}^{{k_{\Lambda_1}}-1} \cdots \sum\limits_{k_{\Lambda_{\lambda_{\rho}}} = 0}^{k_{\Lambda_{\lambda_{\rho} -1}} -1} %Hier weitermachen!!!
\sum\limits_{ \ell_{A_1}, \ldots, \ \ell_{A_{\Lambda_{\lambda_{\rho}}}} } 
g_{\theta}(\np) \cdot \prod_{i=1}^{\lambda_{\rho}} \np^{b_i (p-1) k_{\Lambda_i}} \cdot
\prod\limits_{i=1}^{\Lambda_{\lambda_{\rho}}} \np^{a_i(\ell_{A_i}-1)} \nonumber \\
&\cdot \np^{-{\left[(p-1) \sum\limits_{i=1}^{\Lambda_{\lambda_{\rho}}} \sum\limits_{j=1+A_{i-1}}^{A_i} p^{r-j} (pk_{A_i} + \ell_{A_i} + 1) \right]}s}.
\end{align}
By $\sum\limits_{\ell_{A_1}, \ldots, \ \ell_{A_{\Lambda_{\lambda_{\rho}}}}}$, we mean a summation over all integral points in $\left[ 1,p-1\right]^{\Lambda_{\lambda_{\rho}}}$ respecting the restrictions which are encoded in the 2-level composition $\theta$. %That is $p-1 \geq \ell_{A_{\Lambda_{i-1} + 1}} \geq \ldots \geq \ell_{A_{\Lambda_{i}}} \geq 1$ for $1 \leq i \leq \lambda_{\rho}$.

\begin{definition} \label{Constants}
For $1 \leq i \leq \lambda_{\rho}$, we define the function
\begin{align*}
\alpha_i(s) := b_i(p-1) - (p-1) \sum\limits_{j=1+B_{i-1}}^{B_i} p^{r+1-j} s
\end{align*}
where $B_i := \sum\limits_{k=1}^{i} b_k$.
We further define 
\begin{align*}
\psi_{\ell}(\theta) := \sum\limits_{\ell_{A_1}, \ldots, \ \ell_{A_{\Lambda_{\lambda_{\rho}}}}} g_{\theta}(\np) \cdot \prod\limits_{i=1}^{\Lambda_{\lambda_{\rho}}} \np^{\beta_i(s) }
\end{align*} 
with $\beta_i(s) := a_i(\ell_{A_i}-1) -(p-1) \sum\limits_{j=1+A_{i-1}}^{A_i} p^{r-j} (\ell_{A_i} + 1)s$ where $A_i := \sum\limits_{k=1}^{i} a_k$.
\end{definition}

By ordering the occurring terms by their respective variables and using the previously defined functions, we deduce \eqref{eq1} = 

\begin{align} \label{eq2}
	\sum\limits_{k_{\Lambda_1} \geq 0} \np^{\alpha_1(s) k_{\Lambda_1}}  \sum\limits_{k_{\Lambda_2} = 0}^{k_{\Lambda_1} - 1}  \np^{\alpha_2(s) k_{\Lambda_2}} 
	\cdots \sum\limits_{k_{\Lambda_{\lambda_{\rho}}} = 0}^{k_{\Lambda_{\lambda_{\rho} -1}} -1} 
	\np^{\alpha_{\lambda_{\rho}}(s) k_{\Lambda_{\lambda_{\rho}}}} \cdot \psi_{\ell}(\theta).
	%\\ &\sum\limits_{\ell_{A_1}, \ldots, \ \ell_{A_{\Lambda_{\lambda_{\rho}}}}} g_{\omega}(\np) \cdot \prod\limits_{i=1}^{\Lambda_{\lambda_{\rho}}} \np^{a_i(\ell_{A_i}-1) -(p-1) \sum\limits_{j=1+A_{i-1}}^{A_i} p^{r-j} (\ell_{A_i} + 1)s } 
%	&\sum\limits_{\ell_1 = 1}^{p-1} \np^{a_1(\ell_1-1) - (p-1)\sum\limits_{j=1}^{A_1} p^{r-j} (\ell_1 +1) s} \cdots 
%	\sum\limits_{\ell_{\lambda} = 1}^{p-1} \np^{a_{\lambda}(\ell_{\lambda}-1) - (p-1)\sum\limits_{j=1+A_{\lambda - 1}}^{A_{\lambda}} p^{r-j} (\ell_{\lambda} +1) s}
\end{align}
%By $\sum\limits_{\ell_{A_1}, \ldots, \ \ell_{A_{\Lambda_{\lambda_{\rho}}}}}$ we mean a summation for each variable $\ell_j$ which respects the restrictions given by the prescribed $k\ell$-configuration.

As a next step, we want to explicitly compute the nested geometric series above. The following Lemma provides a strong tool to represent \eqref{eq2} in terms of meromorphic functions.

\begin{lemma} \label{DecompositionGS} %decomposition of the geometric series
Let $X > 1$ be real and let $\alpha_i(s)$ for $1 \leq i  \leq J$ be complex polynomials of degree 1 with real coefficients whose leading coefficients are negative. Then, we have
\begin{align*}
&\quad \sum\limits_{k_1 \geq 0} \left(X^{\alpha_1(s)}\right)^{k_1} \sum\limits_{k_2 = 0}^{k_1 -1} \left(X^{\alpha_2(s)}\right)^{k_2} \ldots \sum\limits_{k_J = 0}^{k_{J-1} -1} \left(X^{\alpha_J(s)}\right)^{k_J} \\
=& \quad X^{\sum\limits_{i=1}^{J-1} (J-i) \alpha_i(s)} \cdot \prod\limits_{i=1}^{J} \left(1 - X^{\sum\limits_{j=1}^{i} \alpha_j(s)}\right)^{-1}
= \frac{1}{1 - X^{\sum\limits_{j=1}^{J} \alpha_j(s)}} \cdot \prod\limits_{i=1}^{J-1} \frac{X^{\sum\limits_{j=1}^{i} \alpha_j(s)}}{1 - X^{\sum\limits_{j=1}^{i} \alpha_j(s)}}
\end{align*}
as a meromorphic continuation.
\end{lemma}

\begin{proof}
According to our assumption, the left-hand side is at least holomorphic for all $s \in \C$ for which $\re({\alpha_i(s)}) < 0$ applies for all $1 \leq i \leq J$. This means that $\re(s)$ is greater than the maximum of all real parts of the zeros of the $\alpha_i(s)$ for $1 \leq i \leq J$.  

Now, we prove the claim by induction on $J$. The case $J=1$ is clear by the well-known identity for geometric series. For $J+1$ we compute
\begin{alignat*}{2}
& \sum\limits_{k_1 \geq 0} \left(X^{\alpha_1(s)}\right)^{k_1} \sum\limits_{k_2 = 0}^{k_1 -1} \left(X^{\alpha_2(s)}\right)^{k_2} \ldots \sum\limits_{k_J = 0}^{k_{J-1} -1} \left(X^{\alpha_J(s)}\right)^{k_J} \sum\limits_{k_{J+1} = 0}^{k_{J} -1} \left(X^{\alpha_{J+1}(s)}\right)^{k_{J+1}} \\
={}& \sum\limits_{k_1 \geq 0} \left(X^{\alpha_1(s)}\right)^{k_1} \sum\limits_{k_2 = 0}^{k_1 -1} \left(X^{\alpha_2(s)}\right)^{k_2} \ldots \sum\limits_{k_J = 0}^{k_{J-1} -1} \left(X^{\alpha_J(s)}\right)^{k_J} \frac{\left(X^{\alpha_{J+1}(s)}\right)^{k_{J}} -1}{X^{\alpha_{J+1}(s)} -1} \\
\stackrel{\text{I.H.}}{=}{}& \frac{X^{\sum\limits_{i=1}^{J-1} (J-i) \alpha_i(s)}}{X^{\alpha_{J+1}}(s) - 1} \cdot \prod\limits_{i=1}^{J-1} \left(1 - X^{\sum\limits_{j=1}^{i} \alpha_i(s)}\right)^{-1} \cdot \left(\frac{1}{1 - X^{\sum\limits_{i=1}^{J+1} \alpha_i(s)}} - \frac{1}{1 - X^{\sum\limits_{i=1}^{J} \alpha_i(s)}}\right) \\
={}& \frac{X^{\sum\limits_{i=1}^{J-1} (J-i) \alpha_i(s)}}{X^{\alpha_{J+1}}(s) - 1} \cdot \prod\limits_{i=1}^{J-1} \left(1 - X^{\sum\limits_{j=1}^{i} \alpha_i(s)}\right)^{-1} \cdot \frac{\left(1 - X^{\sum\limits_{i=1}^{J} \alpha_i(s)} \right) - \left(1 - X^{\sum\limits_{i=1}^{J+1} \alpha_i(s)} \right)}{\left(1 - X^{\sum\limits_{i=1}^{J} \alpha_i(s)} \right) \cdot \left(1 - X^{\sum\limits_{i=1}^{J+1} \alpha_i(s)} \right)} \\
={}& X^{\sum\limits_{i=1}^{J} (J+1-i) \alpha_i(s)} \cdot \prod\limits_{i=1}^{J+1} \left(1 - X^{\sum\limits_{j=1}^{i} \alpha_i(s)}\right)^{-1}.
\end{alignat*}
The second equality is clear.
\end{proof}

\begin{definition} \label{LocalMeromorphicFunction}
For $1 \leq j \leq f$, we define the complex function
\begin{align*}
	\delta_{j, \idp}(s) := 1 - \np^{\eta_{j}(s)}
\end{align*}
with $\eta_{j}(s) := j(p-1) -  p^{r+1-j} (p^j-1) s$.
\end{definition}

We apply Lemma \ref{DecompositionGS} to equation \eqref{eq2} and use the notation provided in Definition \ref{LocalMeromorphicFunction} to obtain \eqref{eq2}=

%\begin{align} %\label{eq3}
%	&\prod\limits_{j=1}^{\lambda_{\rho}} \left(1 - \np^{\sum\limits_{h=1}^j \alpha_h(s)}\right)^{-1} \np^{\sum\limits_{h=1}^{\lambda_{\rho} -1}(\lambda_{\rho} -h) \alpha_h (s)} \cdot \psi_{\ell}(\theta).
%%	 \sum\limits_{\ell_1 = 1}^{p-1} \np^{a_1(\ell_1-1) - (p-1)\sum\limits_{j=1}^{A_1} p^{r-j} (\ell_1 +1) s} \cdots \nonumber \\
%%	 &\sum\limits_{\ell_{\lambda} = 1}^{p-1} \np^{a_{\lambda}(\ell_{\lambda}-1) - (p-1)\sum\limits_{j=1+A_{\lambda - 1}}^{A_{\lambda}} p^{r-j} (\ell_{\lambda} +1) s}
%\end{align}

\begin{align} \label{eq3}
	\frac{1}{1 - \np^{\eta_h(s)}} \cdot \prod\limits_{j=1}^{\lambda_{\rho} - 1} \frac{\np^{\eta_{B_j}(s)}}{1 - \np^{\eta_{B_j}(s)}} \cdot \psi_{\ell}(\theta).
	%	 \sum\limits_{\ell_1 = 1}^{p-1} \np^{a_1(\ell_1-1) - (p-1)\sum\limits_{j=1}^{A_1} p^{r-j} (\ell_1 +1) s} \cdots \nonumber \\
	%	 &\sum\limits_{\ell_{\lambda} = 1}^{p-1} \np^{a_{\lambda}(\ell_{\lambda}-1) - (p-1)\sum\limits_{j=1+A_{\lambda - 1}}^{A_{\lambda}} p^{r-j} (\ell_{\lambda} +1) s}
\end{align}
%If we have $k_i = k_{i+1}$ for some index $1 \leq i \leq \lambda -1$, then we unite the blocks corresponding to the indices $i$ and $i+1$ and we obtain a block with length $a_i + a_{i+1}$. This case behaves similarly to the generic case of the composition $(a_1, \ldots, a_{i-1}, a_i + a_{i+1}, \ldots, a_{\lambda})$ regarding its description as a rational function. We, however, have to keep the further restrictions to the $\ell_i$ in mind. With this reduction, we can restrict ourselves to consider the generic cases for the compositions with all possible ensembles of the $\ell_i$.
Inserting equation \eqref{eq3} in equation \eqref{EinsetzenConfig} yields
\begin{align} \label{PhiMeromorphic}
	\Phi_{f, \idp}(F, C_p^r;s) = 1 +  \sum\limits_{h = 1}^{f} \binom{f}{h}_p \quad \sum\limits_{\theta \in \mathfrak{Z}^{ad}_h} \gamma_{\theta_{\omega}} \cdot \delta_{h, \idp}^{-1}(s) \prod\limits_{j = 1}^{\lambda_{\rho} - 1} \frac{\np^{\eta_{B_j}(s)}}{1 - \np^{\eta_{B_j}(s)}} \cdot \psi_{\ell}(\theta). 
\end{align}

\begin{proposition} \label{LocalMeromorphicContinuation}
	The Euler factor $\Phi_{f, \idp}(F, C_p^r;s)$ has the meromorphic continuation
	\begin{align} \label{LMC}
		\Phi_{f, \idp}(F, C_p^r;s) = \left(\prod\limits_{j=1}^{f} \delta_{j, \idp}^{-1}(s)\right) \cdot \Psi_{f, \idp}(s)
	\end{align}
	with $\Psi_{f, \idp}(s) :=$
	\begin{align*}
		\prod\limits_{j=1}^{f} \delta_{j, \idp}(s) +  \sum\limits_{h = 1}^{f} \binom{f}{h}_p \quad \sum\limits_{\theta \in \mathfrak{Z}^{ad}_h} \gamma_{\theta_{\omega}} \prod\limits_{\substack{j=1 \\
				j \not \in \{B_1, \ldots, B_{\lambda_{\rho}}\} }}^{f} \delta_{j, \idp}(s) \cdot \prod\limits_{j=1}^{\lambda_{\rho} - 1} \np^{\eta_{B_j}(s)} \cdot \psi_{\ell}(\theta).
	\end{align*}
	The function $\Psi_{f, \idp}(s)$ is a polynomial in $\np^{-s}$, whose coefficients also depend on $\np$, and is therefore holomorphic on the whole complex plane. 
	
	Furthermore, $\Psi_{f, \idp}(s)$ does not vanish for $s \in \R$ satisfying $s \geq \frac{f(p-1)}{p^{r+1-f}(p^f - 1)}$. In particular, $\Phi_{f, \idp}(F, C_p^r;s)$ has a single pole at $\frac{f(p-1)}{p^{r+1-f}(p^f - 1)}$.
\end{proposition}

\begin{proof}
	We separate the meromorphic terms in equation \eqref{PhiMeromorphic} and obtain equation \eqref{LMC}. It is clear by the definitions of $\delta_{j, \idp}(s), \eta_{B_j}(s)$ and $\psi_{\ell}(\theta)$ that $\Psi_{f, \idp}(s)$ is a polynomial in $\np^{-s}$, whose coefficients depend on $\np$. 
	
	 %The holomorphy of $\Psi_{f, \idp}(s)$ is clear, since $\Psi_{f, \idp}(s)$ is composed as a finite sum respective product of holomorphic complex functions.
	Now, we prove the non-vanishing of $\Psi_{f, \idp}(s)$ for $\frac{f(p-1)}{p^{r+1-f}(p^f - 1)} \leq s \in \R$. By the definition of $\eta_j(s)$ in Definition \ref{LocalMeromorphicFunction}, we deduce $0 > \eta_j(s) \in \R$ and therefore $0 < \delta_{j,\idp}(s) \in \R$ for $1 \leq j < f$ and $s \geq \frac{f(p-1)}{p^{r+1-f}(p^f - 1)}$. For $\eta_f(s)$, we have $0 \geq \eta_f(s) \in \R$ and therefore $0 \leq \delta_{f,\idp}(s) \in \R$ for $s \geq \frac{f(p-1)}{p^{r+1-f}(p^f - 1)}$, where equality applies for $s = \frac{f(p-1)}{p^{r+1-f}(p^f - 1)}$ only.
	
For a fixed $\theta \in \mathfrak{Z}^{ad}_h$, we find that the value of $\prod\limits_{j=1}^{\lambda_{\rho} - 1} \np^{\eta_{B_j}(s)} \cdot \psi_{\ell}(\theta)$ for some real $s \geq \frac{f(p-1)}{p^{r+1-f}(p^f - 1)}$ is real and non-negative. Additionally, it does not vanish if $g_{\theta}(\np)$ %contained in $\psi_{\ell}(\theta)$ 
is not zero which is at least true for the 2-level composition whose outer composition is the all-one composition $(1, \ldots, 1)$ of length $h$. 

For $s = \frac{f(p-1)}{p^{r+1-f}(p^f - 1)}$, all summands with $h < f$ contain the factor $\delta_{f,p}(s)$ and vanish, as $\delta_{f,p}(s)$ vanishes at this point. Nonetheless, the summand for the 2-level composition whose outer composition is the all-one composition $(1, \ldots, 1)$ of length $f$ is positive.
Accordingly, we have $\Psi_{f, \idp}(s) > 0$ for $\frac{f(p-1)}{p^{r+1-f}(p^f - 1)} \leq s \in \R$. 

Finally, $\delta_{f, \idp}^{-1}(s)$ has a single pole at $s = \frac{f(p-1)}{p^{r+1-f}(p^f - 1)}$ so that this also applies to $\Phi_{f, \idp}(F, C_p^r;s)$.
\end{proof}

\begin{corollary} \label{PhiRational}
The Dirichlet series 
	\begin{align} \label{Rational}
\Phi(F_{\idp}, C_p^r;s) &= \sum\limits_{f=0}^{r} e_f \cdot \Phi_{f, \idp}(F, C_p^r;s) \notag \\
&= \prod\limits_{j=1}^{r} \delta_{j, \idp}^{-1}(s) \sum\limits_{f=0}^{r} e_f \prod\limits_{j=f+1}^{r} \delta_{j, \idp}(s) \cdot \Psi_{f, \idp}(s).
	\end{align}
corresponding to the counting function of the local function field $F_{\idp}$ is a rational function in $q^{-s}$. In particular, its coefficients satisfy a linear recurrence.
\end{corollary}

\begin{proof}
Firstly, let us recall Proposition \ref{CountingFunction}. Proposition \ref{LocalMeromorphicContinuation} thus gives equation \eqref{Rational}. Additionally, Proposition \ref{LocalMeromorphicContinuation} yields that $\Phi_{f, \idp}(F, C_p^r;s)$ and thus $\Phi(F_{\idp}, C_p^r;s)$ is rational in $q^{-s}$ whereby we recall $\np^{-s} = \left(q^{-s}\right)^{\deg(\idp)}$. 

If we substitute $q^{-s} = t$, we can interpret the Dirichlet series $\Phi(F_{\idp}, C_p^r;s)$ as a power series in the variable $t$. %A comparison of coefficients yields a linear recurrence on the Dirichlet series' coefficients in terms of the underlying rational function. 
A comparison of the coefficients leads to a linear recursion for the coefficients of the Dirichlet series in relation to the underlying rational function.
\end{proof}

\begin{remark}
	We note that the decomposition of $\Phi_{f, \idp}(F, C_p^r;s)$ into meromorphic terms, which we obtained from Lemma \ref{DecompositionGS}, can also be described by a properly modified Igusa function of degree $f$, which is defined in \cite[Definition 2.2]{CSV}, for example. We can naturally extend the definition of the Igusa functions to the degree $0$ by defining $I_0 \equiv 1$.
%	Inserting equation \eqref{eq3} in equation \eqref{EinsetzenConfig}, yields $\Phi_{f, \idp}(\Fqtg, C_p^r;s) =$
%	\begin{align} \label{PhiMeromorphic}
%		1 +  \sum\limits_{h = 0}^{f} \binom{f}{h}_p \quad \sum\limits_{\theta \in \mathfrak{Z}^{ad}_h} \gamma_{\theta_{\omega}} \cdot \delta_{h, \idp}^{-1}(s) \prod\limits_{j = 1}^{\lambda_{\rho} - 1} \frac{\np^{\eta_{B_j}(s)}}{1 - \np^{\eta_{B_j}(s)}} \cdot \psi_{\ell}(\theta). 
%	\end{align}
	Indeed, the sum in \eqref{PhiMeromorphic} coincides with 
	\begin{align}
		I_f\left(p; \np^{\eta_1(s)}, \ldots, \np^{\eta_f(s)}\right) + \text{lower order terms}
	\end{align}
	defined in \cite[Definition 2.2]{CSV} since we find the factor $\np^{\eta_h(s)}$ for a fixed $\theta \in \mathfrak{Z}^{ad}_h$, if we consider the summand of $\psi_{\ell}(\theta)$ corresponding to $\ell_i = p-1$ for all $1 \leq i \leq h$ and the leading term $\np^h$ of the strucure polynomial of $\theta$. The remaining summands of $\psi_{\ell}(\theta)$ correspond to the ``lower order terms".
\end{remark}

\subsection{Asymptotics of elementary-abelian extensions of local function fields}  %\hfill \\

\begin{corollary} \label{LocalPoles}
For $1 \leq j \leq r$, the meromorphic complex function $\delta_{j, \idp}^{-1}(s)$ has simple poles at
\begin{align*}
s = \frac{j(p-1)}{p^{r+1-j} (p^j - 1)} + \frac{2\pi i}{\log\left(\np\right)} \cdot \frac{k}{p^{r+1-j} (p^j - 1)} + \frac{2 \pi i}{\log\left(\np\right)} \cdot \Z
\end{align*}
for $1 \leq k \leq p^{r+1-j} (p^j - 1)$ and no further poles.
\end{corollary}

\begin{proof}
By the definition of $\delta_{j,\idp}(s)$ in Definition \ref{LocalMeromorphicFunction}, we have to solve the equation
\begin{align*}
\np^{j(p-1)} = \np^{p^{r+1-j} (p^j - 1)s}.
\end{align*}
Because of $\np = e^{\log\left(\np\right)}$, we obtain the equation
\begin{align*}
e^{\log\left(\np\right) \cdot j(p-1)} = e^{ \log\left(\np\right) \cdot p^{r+1-j} (p^j - 1)s}.
\end{align*}
Taking into account the periodicity of the complex exponential function, we obtain the claimed simple poles.
\end{proof}

By Proposition \ref{LocalMeromorphicContinuation}, all those values of simple poles of $\delta_{j,\idp}^{-1}(s)$ for $1 \leq j \leq f$ computed in Corollary \ref{LocalPoles} are candidates for simple poles of the Euler factor $\Phi_{f, \idp}(F, C_p^r;s)$. Whether or not $\Psi_{f, \idp}(s)$ vanishes at the respective value determines whether $\Phi_{f, \idp}(F, C_p^r;s)$ actually has a pole or is holomorphically continuable. At least, we have already proven that the rightmost candidate on the real line is definitely a simple pole of $\Phi_{f, \idp}(F, C_p^r;s)$. 

By the meromorphic continuation of $\Phi(F_{\idp}, C_p^r;s)$ presented in Corollary \ref{PhiRational} and the computation of the poles of the functions $\delta_{j, \idp}^{-1}(s)$ in Corollary \ref{LocalPoles}, we can compute all poles of the Dirichlet series $\Phi(F_{\idp}, C_p^r;s)$. However, we still have to decide which ones are furthest to the right. This is possible with the following Lemma.

\begin{lemma} \label{LocalInequalityAbscissa}
For $2 \leq f \leq r$, we have 
\begin{align*} %\label{Ineq1}
\frac{(f-1)(p-1)}{p^{r+1 - (f-1)} (p^{f-1} -1)} < \frac{f(p-1)}{p^{r+1-f}(p^f -1)}.
\end{align*}
\end{lemma}

\begin{proof}
We have
\begin{align*}
& \frac{(f-1)(p-1)}{p^{r+1 - (f-1)} (p^{f-1} -1)} < \frac{f(p-1)}{p^{r+1-f}(p^f -1)} \\
\Longleftrightarrow \quad &(f-1) \cdot (p^f - 1) < f \cdot p^{} \cdot (p^{f-1} - 1) \\
\Longleftrightarrow \quad & 1 + f \cdot (p-1) < p^f.
\end{align*}
The last inequality is Bernoulli's inequality.
\end{proof}

%\begin{theorem} \label{LocalAsymptotics} %deducing the local asymptotics by the Euler factors ???
%	Let $L_{\text{loc}}(C_p^r) = p(p^r -1)$, $X = q^m$ and $m \to \infty$ with $m \equiv e \mod L_{\text{loc}}(C_p^r)$ for a fixed $0 \leq e < L_{\text{loc}}(C_p^r)$. \\
%	Then, there exists at least one arithmetic progression $m \equiv e \mod L_{\text{loc}}(C_p^r)$ so that there exists a constant $c(\F_q((t)), C_p^r,e) \in \R_{> 0}$ and the number of $C_p^r$-extensions of a local function field $\F_q((t))$ when counted by discriminants satisfies the following asymptotic equivalence:
%	\begin{align*}
%		Z(\F_q((t)), C_p^r; X) \sim 
%		c(\F_q((t)), C_p^r, e) \cdot X^{\frac{r(p-1)}{p(p^r-1)}}.
%	\end{align*}
%	Moreover, we have $Z(\F_q((t)), C_p^r; X) \ - \  c(\F_q((t)), C_p^r, e) \cdot X^{\frac{r(p-1)}{p(p^r-1)}} \in  O\left(X^{\frac{(r-1)(p-1)}{p^2(p^{r-1} -1)}}\right)$ for $r \geq 2$ and $Z(\F_q((t)), C_p; X) \ - \  c(\F_q((t)), C_p, e) \cdot X^{\frac{p-1}{p(p-1)}} \in  O\left(X^{\epsilon}\right)$ for every $\epsilon > 0$.
%\end{theorem}

\begin{theorem} \label{LocalAsymptotics} 
	Let $\Fqtl$ be a local function field, $L_{\text{loc}}(C_p^r) := p(p^r -1)$, $X = q^m$ and $m \to \infty$. Then, there exists a constant $c(\F_q((t)), C_p^r, m \bmod L_{\text{loc}}(C_p^r)) \in \R_{\geq 0}$ depending on the residue class of $m$ modulo $L_{\text{loc}}(C_p^r)$ such that the number of $C_p^r$-extensions of a local function field $\F_q((t))$ when counted by discriminants satisfies
	\begin{align*}
	&Z(\F_q((t)), C_p^r; X) \\
	= \ &c(\F_q((t)), C_p^r, m \bmod L_{\text{loc}}(C_p^r)) \cdot X^{\frac{r(p-1)}{p(p^r-1)}} + \begin{cases}
		O(1) & \text{for } r=1, \\
		O\left(X^{\frac{(r-1)(p-1)}{p^2(p^{r-1} -1)}}\right) & \text{for } r >1.
	\end{cases}
	\end{align*}
	The constant is positive for at least one residue class.
%	Moreover, the constant is positive for at least one residue class and we have the asymptotic equivalence
%	\begin{align*}
%				Z(\F_q((t)), C_p^r; X) \sim 
%				c(\F_q((t)), C_p^r, e \mod L_{\text{loc}}(C_p^r)) \cdot X^{\frac{r(p-1)}{p(p^r-1)}}
%	\end{align*}
%	for $m \to \infty$ with $m \equiv e \mod L_{\text{loc}}(C_p^r)$ and some fixed $0 \leq e < L_{\text{loc}}(C_p^r)$ such that $c(\F_q((t)), C_p^r, e \mod L_{\text{loc}}(C_p^r)) > 0$.
\end{theorem}

\begin{proof}
%%Without loss of generality, we may assume $\idp = (t)$ so that we have $\deg(\idp) = 1$, $\np = q$ and $\Fqtg_{(t)} \cong \Fqtl$. 
%Our aim is to apply \cite[Lemma A.2, Lemma A.4 and Theorem A.5]{La2}, which is a Tauberian theorem for power series generalising \cite[Theorem 17.4, p. 311]{Ro}. We note that the error estimate in Lagemann's theorem can be improved by considering the position of the second rightmost pole(s). For example, we refer to \cite[Theorem IV.9]{FS}. \\
By Corollary \ref{PhiRational}, the Dirichlet series $\Phi(\Fqtl, C_p^r;s)$ corresponding to the counting function of the local function field $\Fqtl$ is a rational function in $q^{-s}$ and has the form given in equation \eqref{Rational}. Substituting $z = q^{-s}$, we can transform $\Phi(\Fqtl, C_p^r;s)$ into a power series in $z$ which we denote by $\phi(\Fqtl, C_p^r;z) = \sum\limits_{m \geq 0} c_m z^m$. By Proposition \ref{LocalMeromorphicContinuation}, Corollary \ref{LocalPoles} and Lemma \ref{LocalInequalityAbscissa}, $\phi(\Fqtl, C_p^r;z)$ has radius of convergence $R = q^{- \frac{r(p-1)}{p(p^r - 1)}}$ and is meromorphically continuable on the whole complex plane except for a discrete set of simple poles. According to Corollary \ref{LocalPoles}, the candidates for poles on the circle of convergence are given by $\alpha_j := R\xi^{-j}$ for some primitive $L_{\text{loc}}(C_p^r)$-th root of unity $\xi$ and $1 \leq j \leq L_{\text{loc}}(C_p^r)$ and we have at least one simple pole on the radius of convergence at $z = R$ by Proposition \ref{LocalMeromorphicContinuation}. 

Additionally, by Corollary \ref{LocalPoles}, the candidates for subsequent poles lie on the circles with radii $q^{-\frac{(r-1)(p-1)}{p^2(p^{r-1} -1)}}$ for $r \geq 2$ and all these candidates have pole order at most $1$. In the case $r=1$, we do not have any further poles. 

Now, we apply a Tauberian theorem \cite[Theorem IV.9]{FS} to obtain
\begin{align*}
c_m = \sum\limits_{j=1}^{L_{\text{loc}}(C_p^r)} \tau_j \alpha_j^{-m} + \begin{cases}
	O(1) & \text{for } r=1, \\
	O\left(q^{m \cdot \frac{(r-1)(p-1)}{p^2(p^{r-1} -1)}}\right) & \text{for } r >1,
\end{cases}
\end{align*}
for some coefficients $\tau_j \in \C$ which involve the residues at the respective poles $\alpha_j$. For the main term, we further compute
\begin{align*}
\sum\limits_{j=1}^{L_{\text{loc}}(C_p^r)} \tau_j \alpha_j^{-m} = \sum\limits_{j=1}^{L_{\text{loc}}(C_p^r)} \tau_j \left(R\xi^{-j}\right)^{-m} = \left(\sum\limits_{j=1}^{L_{\text{loc}}(C_p^r)} \tau_j \xi^{jm}\right) \cdot X^{\frac{r(p-1)}{p(p^r-1)}}
\end{align*}
which yields $c(\F_q((t)), C_p^r, m \bmod L_{\text{loc}}(C_p^r)) = \sum\limits_{j=1}^{L_{\text{loc}}(C_p^r)} \tau_j \xi^{jm}$. We note that this sum depends on the residue class of $m \bmod L_{\text{loc}}(C_p^r)$ since $\xi$ is a primitive $L_{\text{loc}}(C_p^r)$-th root of unity. 

Finally, we have $\sum\limits_{j=1}^{L_{\text{loc}}(C_p^r)} \tau_j \xi^{jm} \in \R_{\geq 0}$ as the coefficients $c_m$ are non-negative integers and the constant is positive for at least one residue class because $\phi(\Fqtl, C_p^r;z)$ actually has a simple pole at $z=R$.
\end{proof} 

\begin{remark} \label{MTLR}
Note that we could also apply \cite[Appendix]{La2}, which gives a Tauberian theorem for power series generalising \cite[Theorem 17.4]{Ro}. Indeed, Lagemann also provides an expression for the constant $c(\F_q((t)), C_p^r, m \bmod L_{\text{loc}}(C_p^r))$ which allows an explicit computation. However, we did not use his theorem in the previous proof due to the worse error bound compared to \cite[Theorem IV.9]{FS}. 

Note that a precise computation of the constant $c(\F_q((t)), C_p^r, m \bmod L_{\text{loc}}(C_p^r))$ involves a computation of the residue at each possible pole $\alpha_j = R \xi^{-j}$ and a complicated and tedious summation over these residues and the $L_{\text{loc}}(C_p^r)$-th roots of unity. %For the residue $p_j$ at the pole $\alpha_j$, one can compute $p_j = -e_r \cdot \frac{\alpha_j}{L_{\text{loc}}(C_p^r)} \cdot \Upsilon(\alpha_j)$ where $\Upsilon(\alpha_j)$ is the value of the function $\Upsilon(z) := \prod\limits_{i=1}^{r-1} \delta_{i, (t)}^{-1}(t) \cdot \Psi_{r, (t)}(t)$ at $\alpha_j$. Note that $\Upsilon(z)$ originates from Corollary \ref{PhiRational} and we applied the substitution $z = q^{-s}$. By Lagemann's formula in \cite[Theorem A.5]{La2}, this yields $\tau_j = \frac{e_r}{L_{\text{loc}}(C_p^r)} \cdot \Upsilon(\alpha_j)$ and thus
For the constant, one can compute $\tau_j = \frac{e_r}{L_{\text{loc}}(C_p^r)} \cdot \Upsilon_{\text{loc}}(\alpha_j)$ and thus
\begin{align*}
c(\F_q((t)), C_p^r, m \bmod L_{\text{loc}}(C_p^r)) = \frac{e_r}{L_{\text{loc}}(C_p^r)} \cdot \sum\limits_{j=1}^{L_{\text{loc}}(C_p^r)} \xi^{jm} \cdot \Upsilon_{\text{loc}}(\alpha_j)
\end{align*}
where $\Upsilon_{\text{loc}}(\alpha_j)$ is the value of the function $\Upsilon_{\text{loc}}(z) := \prod\limits_{i=1}^{r-1} \delta_{i, (t)}^{-1}(z) \cdot \Psi_{r, (t)}(z)$ at $\alpha_j$. The function $\Upsilon_{\text{loc}}(z)$ originates from Corollary \ref{PhiRational} and we applied the substitution $z = q^{-s}$.
\end{remark}

\begin{remark}
The result presented in Theorem \ref{LocalAsymptotics} restates Lagemann's results for elementary-abelian extensions of local function fields \cite{La1} using a different approach. 

Additionally, we have shown that each composition of $r$ which satisfies the condition $a_i \leq n$ for all blocks where $q = p^n$ contributes to the main term of the asymptotics. Accordingly, all those families of totally ramified extensions corresponding to a composition of $r$ satisfying the previous conditions have a positive density. However, extensions which are not totally ramified thus contain an unramified $C_p$-extension occur in the error term.
\end{remark}

\subsection{Meromorphic continuation of the Euler products and the Dirichlet series of $\Fqtg$} \label{GlobalContinuation}  \hfill \\
In the next stage, we want to determine the asymptotics of elementary-abelian extensions of the rational function field $\Fqtg$. Among other things, we compute the asymptotic behaviour of $C_p$-extensions, which has originally been figured out by Lagemann in \cite{La1} and \cite{La2}. 

To this end, we first consider the Euler products $\Phi_{f}(F, C_p^r; s)$ presented in Definition \ref{DefinitionDirichletSeries} and compute a meromorphic continuation in each case. Let us recall again that, with a view to Section \ref{ArbitraryBaseFields}, we will analyse these Euler products for an arbitrary global function field $F$, even if we are currently interested in the special case $F = \Fqtg$.

Furthermore, we recall the zeta function $\zeta_F(s)$ of $F$ and its meromorphic continuation
\begin{align*}
\zeta_F(s) = \frac{L_F\left(q^{-s}\right)}{\left(1-q^{-s}\right)\left(1-q^{1-s}\right)}
\end{align*}
to the whole complex plane where $L_F(X) \in \Z[X]$ is a polynomial of degree $2g_F$, the so called \emph{$L$-polynomial of $F$}, which satisfies the equations $L_F(0) = 1$ and $L_F(1) = \vert {\Pic^0_F} \vert$ (e.g.\ \cite[Theorem 5.9]{Ro}). Furthermore, all roots of $L_F\left(q^{-s}\right)$ satisfy $\re(s) = \frac{1}{2}$ by the famous Hasse-Weil Theorem which is sometimes referred to as the Riemann Hypothesis for function fields (e.g.\ \cite[Theorem 5.10]{Ro}). All poles of $\zeta_F(s)$ correspond to zeros of the denominator and lie on the axes $\re(s) = 0$ and $\re(s) = 1$. 

By the meromorphic continuations of the Euler factors $\Phi_{f, \idp}(F, C_p^r;s)$ given in Proposition \ref{LocalMeromorphicContinuation},
the zeta functions 
\begin{align*}
	\prod\limits_{\idp \in \PP_F} \delta_{j,\idp}^{-1}(s) = \zeta_{F}\left(p^{r+1-j}(p^j -1)s - j(p-1)\right)
\end{align*}
for $1 \leq j \leq f$ will arise naturally in the meromorphic continuation of $\Phi_{f}(F, C_p^r; s)$. Hence, we first compare the convergence abscissas $\frac{1+j(p-1)}{p^{r+1-j}\left(p^j-1\right)}$ of these zeta functions for $1 \leq j \leq r$ in the following Lemma \ref{GlobalInequality} since they are candidates for the convergence abscissa of the meromorphic continuations of the Euler product $\Phi_{f}(F, C_p^r;s)$ for $1 \leq f \leq r$. Additionally, we give a more precise estimate for the axes $\re(s) = \frac{1 - \frac{1}{p} + j(p-1)}{p^{r+1-j}\left(p^j-1\right)}$ until which we continue meromorphically in Proposition \ref{ContinuationEulerProduct}.

\begin{lemma} \label{GlobalInequality}
	Let $2 \leq j \leq r$. Then, the inequalities
	\begin{align*}
		\frac{1+(j-1)(p-1)}{p^{r+1-(j-1)}\left(p^{j-1}-1\right)} \leq \frac{1 - \frac{1}{p} + j(p-1)}{p^{r+1-j}\left(p^j-1\right)} < \frac{1 + j(p-1)}{p^{r+1-j}\left(p^j-1\right)}
	\end{align*}
	hold except for the case $j=p=2$. Equality on the left-hand side applies if and only if $p=2$ and $j=3$. 
	
	In the case $j=p=2$, we have
	\begin{align*}
		\frac{1+(j-1)(p-1)}{p^{r+1-(j-1)}\left(p^{j-1}-1\right)} = \frac{1}{2^{r-1}} = \frac{1+j(p-1)}{p^{r+1-j}\left(p^j-1\right)}.
	\end{align*}
\end{lemma}

\begin{proof}
The inequality on the right-hand side and the case $j=p=2$ are clear. For the inequality on the left-hand side, we find
	\begin{align*}
		& \frac{1 + (j-1)(p-1)}{p^{r+1 - (j-1)} (p^{j-1} -1)} \leq \frac{1 - \frac{1}{p} + j(p-1)}{p^{r+1-j}(p^j -1)} \\
		\Longleftrightarrow \quad & p^j - 1 + (p^j - 1) \cdot (j-1) \cdot (p - 1) + (p^{j-1} - 1) \leq p^j - p + j \cdot (p-1) \cdot (p^j - p) \\
		\Longleftrightarrow \quad & 1 + (p^j - 1) \cdot (j-1) + \sum\limits_{i=1}^{j-1} p^{j-1-i} \leq j \cdot (p^j - p) \\
		\Longleftrightarrow \quad & 2 + (p-1) \cdot j + \sum\limits_{i=1}^{j-1} p^{j-1-i} \leq p^j = 1 + (p-1)p^{j-1} + (p-1) \sum\limits_{i=1}^{j-1} p^{j-1-i} \\
		\Longleftrightarrow \quad & 1 + j(p-1) \leq (p-1)p^{j-1} + \underbrace{(p-2)\sum\limits_{i=1}^{j-1} p^{j-1-i}}_{\geq 0}.
	\end{align*}
	By induction on $j$, we finally prove the inequality $\frac{1}{p-1} + j \leq p^{j-1}$. We start with $j=3$ in the case $p=2$ and with $j=2$ otherwise. Obviously, the equality only applies in the case $p=2$ and $j=3$.
\end{proof}

Next, we define a product $\Lambda_{F,f,r,p}(s)$ of zeta functions of $F$ depending on $f, r$ and $p$ which we will compare with the Euler product $\Phi_{f}(F, C_p^r; s)$. %In fact, we will see that the rightmost pole(s) of $\Phi_{f}(F, C_p^r; s)$ of maximum order coincide(s) with the rightmost pole(s) of maximum order of $\Lambda_{F,f,r,p}(s)$.

\begin{definition} \label{MeromorphicFunctionGlobal}
For $1 \leq f \leq r$, we define the function
\begin{align*}
\Lambda_{F,f,r,p}(s) := \begin{cases}
\prod\limits_{\ell=1}^{p-1} \zeta_{F}(p^{r-1}(\ell+1)(p-1)s - \ell) \ &\text{for} \ f=1,\\
\zeta_{F}(3 \cdot 2^{r-1} s - 2) \cdot \zeta_{F}(2^{r}s - 1)^{2^f - 1} \ &\text{for} \ f=p=2,\\
\zeta_{F}(p^{r+1-f}(p^f-1)s - f(p-1)) \ &\text{else}.
\end{cases}
\end{align*}
We remark that the zeta function belonging to $\ell = p-1$ in the first case and the first zeta function in the second case coincide with the zeta function in the third case. The case distinction for $f = p = 2$ has already been motivated by Lemma \ref{GlobalInequality}.
\end{definition}

\begin{proposition} \label{GlobalPoles}
	The function $\Lambda_{F,f,r,p}(s)$ has the following candidates for poles on its convergence abscissa, which is the axis with real part $\frac{1 + f(p-1)}{p^{r+1-f}(p^f -1)}$. 
	
	\begin{enumerate}
		\item[a)] Let $f = 1$. Then, all candidates for such poles are given by
		\begin{align*}
			\frac{1}{(p-1) p^{r-1}} + \frac{2 \pi i}{\log(q)} \cdot \frac{j}{(p-1) p^{r-1} M} + \frac{2 \pi i}{\log(q)} \cdot \Z
		\end{align*}
		with $0 \leq j < (p-1) p^{r-1} M$ where $M$ is the least common multiple of the numbers $2, \ldots, p$. We note that candidates do not necessarily have to be actual poles, which we allow to simplify the description. 
		
		For $p=2$, all candidates are simple poles. For $p \neq 2$, the candidates with $j= M j'$ for any $0 \leq j' < (p-1)p^{r-1}$ are actual poles and have exactly the order $p-1$. The remaining candidates have smaller order, which may also be 0 or negative. For $p = 3$ and $r=1$, $\Lambda_{F,1,1,3}(s)$ is actually holomorphic for $j \in \{1,5,7,11\}$, for example.
		\item[b)] Let $f = p = 2 \leq r$. Then, all possible candidates for such poles are given by 
		\begin{align*}
			\frac{1}{2^{r-1}} + \frac{2 \pi i}{\log(q)} \cdot \frac{j}{3 \cdot 2^r} + \frac{2 \pi i}{\log(q)} \cdot \Z
		\end{align*}
		with $0 \leq j < 3 \cdot 2^r$. We note that candidates do not necessarily have to be actual poles, which we allow again to simplify the description. 
		
		The candidates with $j = 6 j'$ for any $0 \leq j' < 2^{r-1}$ are actual poles and have exactly the order $4$. The order of each other candidate is at most 3 and can also be 0 or negative.
		\item[c)] In the remaining cases not covered by a) and b), we have simple poles at
		\begin{align*}
			\frac{1 + f(p-1)}{p^{r+1-f}(p^f-1)} + \frac{2 \pi i}{\log(q)} \cdot \frac{j}{p^{r+1-f}(p^f-1)} + \frac{2 \pi i}{\log(q)} \cdot \Z
		\end{align*}
		with $0 \leq j < p^{r+1-f}(p^f-1)$. 
		
		We note that this case is easier to handle because $\Lambda_{F,f,r,p}(s)$ only consists of a single zeta function, so we can directly list all poles and their orders instead of talking about possible candidates.
	\end{enumerate}
\end{proposition}

\begin{proof}
	By \cite[Theorem 5.9]{Ro}, the poles of $\zeta_{F}(s)$ are given by the solutions of the equations $q^s = 1$ and $q^s = q$. By Definition \ref{MeromorphicFunctionGlobal}, we have to compute the respective poles of the zeta functions being included in $\Lambda_{F,f,r,p}(s)$.
	\begin{enumerate}
		\item[a)] See \cite[Proposition 6.6(a)]{La2}, apply an index shift by $1$ and consider the weighting by the factor $(p - 1)p^{r-1}$ for the count by discriminant and for the group $C_p^r$.
		\item[b)] Considering the periodicity of the complex exponential function, the zeta function $\zeta_{F}(3 \cdot 2^{r-1}s - 2)$ has simple poles exactly at the points
		\begin{align*} 
			s &= \frac{1}{3 \cdot 2^{r-2}} + \frac{2 \pi i}{\log(q)} \cdot \frac{j}{3 \cdot 2^{r-1}} + \frac{2 \pi i}{\log(q)} \cdot \Z \\
			s &= \frac{1}{2^{r-1}} + \frac{2 \pi i}{\log(q)} \cdot \frac{j}{3 \cdot 2^{r-1}} + \frac{2 \pi i}{\log(q)} \cdot \Z
		\end{align*}
		with $0 \leq j < 3 \cdot 2^{r-1}$ corresponding to the solutions of the equations $q^{3 \cdot 2^{r-1}s} = q^2$ and $q^{3 \cdot 2^{r-1}s} = q^3$. Analogously, the zeta function $\zeta_{F}(2^{r}s - 1)$ has simple poles exactly at the points
		\begin{align*} 
			s &= \frac{1}{2^r} + \frac{2 \pi i}{\log(q)} \cdot \frac{j}{2^r} + \frac{2 \pi i}{\log(q)} \cdot \Z \\
			s &= \frac{1}{2^{r-1}} + \frac{2 \pi i}{\log(q)} \cdot \frac{j}{2^r} + \frac{2 \pi i}{\log(q)} \cdot \Z
		\end{align*}
		with $0 \leq j < 2^r$ corresponding to the solutions of the equations $q^{2^rs} = q$ and $q^{2^rs} = q^2$. The poles with real part $\frac{1}{2^r}$ and $\frac{1}{3 \cdot 2^{r-2}}$ are not on the critical axis with $\re(s) = \frac{1}{2^{r-1}}$. 
		
		Since the greatest common divisor, respectively the least common multiple, of $2^r$ and $3 \cdot 2^{r-1}$ is $2^{r-1}$, respectively $3 \cdot 2^r$, the possible candidates for poles of $\Lambda_{F,2,r,2}(s)$ with $\re(s) = \frac{1}{2^{r-1}}$ are as claimed and the poles of order exactly $4$ are at the points
		\begin{align*}
			s &= \frac{1}{2^{r-1}} + \frac{2 \pi i}{\log(q)} \cdot \frac{j}{2^{r-1}} + \frac{2 \pi i}{\log(q)} \cdot \Z
		\end{align*}
		with $0 \leq j < 2^{r-1}$. Note that these are the common poles of $\zeta_{F}(3 \cdot 2^{r-1}s - 2)$ and $\zeta_{F}(2^{r}s - 1)$ on the convergence abscissa. %For $j \in \{0,6\}$, the zeta functions have common poles and we obtain pole order $4$. Poles of order $3$ corresponding to $\zeta_{\Fqtg}(4s - 1)^{3}$ are obtained for $j \in \{3,9\}$ and simple poles corresponding to $\zeta_{\Fqtg}(6s - 2)$ are obtained for $j \in \{2,4,8,10\}$.
		\item[c)] Analogous to the approach in b), the function $\zeta_{F}(p^{r+1-f}(p^f-1)s - f(p-1))$ has simple poles exactly at the points
		\begin{align*}
			s &= \frac{f(p-1)}{p^{r+1-f}(p^f - 1)} + \frac{2 \pi i}{\log(q)} \cdot \frac{j}{p^{r+1-f}(p^f-1)} + \frac{2 \pi i}{\log(q)} \cdot \Z \\
			s &= \frac{1 + f(p-1)}{p^{r+1-f}(p^f - 1)} + \frac{2 \pi i}{\log(q)} \cdot \frac{j}{p^{r+1-f}(p^f-1)} + \frac{2 \pi i}{\log(q)} \cdot \Z
		\end{align*}
		with $0 \leq j < p^{r+1-f}(p^f - 1)$ corresponding to the solutions of the equations $q^{p^{r+1-f}(p^f - 1)s} = q^{f(p-1)}$ and $q^{p^{r+1-f}(p^f - 1)s} = q^{1 + f(p-1)}$. %Due to $\frac{r(p-1)}{p(p^r - 1)} < \frac{1 - \frac{1}{p} + r(p-1)}{p(p^r - 1)} $, the poles with real part $\frac{r(p-1)}{p(p^r - 1)}$ are also not in the domain of interest.
	\end{enumerate}
\end{proof}

\begin{proposition} \label{LambdaZeros}
$\Lambda_{F,f,r,p}(s)$ has no zeros in the domain $\re(s) > \frac{\frac{1}{2} + f(p-1)}{p^{r+1-f}(p^f -1)}$. 
\end{proposition}

\begin{proof}
By the Hasse-Weil Theorem (e.g.\ \cite[Theorem 5.10]{Ro}), all zeros of $\zeta_{F}(s)$ satisfy $\re(s) = \frac{1}{2}$. 

Let $f = 1$. Then, all zeros of $\Lambda_{F,1,r,p}(s)$ lie on the axes $\re(s) = \frac{\frac{1}{2} +k}{p^{r-1}(k+1)(p-1)}$ for $1 \leq k \leq p-1$. By differentiation with respect to $k$, we obtain that the fraction is maximised for $k=p-1$, which yields the claim in this case. 

Let $f = p = 2$. Then, all zeros of $\Lambda_{F,2,r,2}(s)$ lie either on the axis $\re(s) = \frac{3}{2^{r+1}}$ for $\zeta_F(2^rs-1)$ or on the axis $\re(s) = \frac{5}{3 \cdot 2^r}$ for $\zeta_F(3 \cdot 2^{r-1} s -2)$. Due to $\frac{3}{2^{r+1}} < \frac{5}{3 \cdot 2^r} = \frac{\frac{1}{2} + 2}{3 \cdot 2^{r-1}}$, we are done in this case.

In the remaining cases, we have $\Lambda_{F,f,r,p}(s) = \zeta_{F}(p^{r+1-f}(p^f-1)s - f(p-1))$ and all zeros satisfy $\re(s) = \frac{\frac{1}{2} + f(p-1)}{p^{r+1-f}(p^f -1)}$, which completes the proof.
\end{proof}

\begin{proposition} \label{ContinuationEulerProduct}
	For $1 \leq f \leq r$, the complex function 
	\begin{align*}
		\psi_{F, f, r, p}(s) := \Phi_{f}(F, C_p^r; s) \cdot \Lambda_{F, f,r, p}^{-1}(s)
	\end{align*}
	has a holomorphic continuation to the domain $\re(s) > \frac{1+f(p-1)}{p^{r+1-f}\left(p^f-1\right)} - \epsilon$ with $\epsilon = \frac{1}{p^{r+2-f} (p^f -1)} > 0$. 
	
	Furthermore, $\psi_{F,f, r, p}(s)$ does not vanish at $s = \frac{1+f(p-1)}{p^{r+1-f}\left(p^f-1\right)}$, so $\Phi_{f}(F, C_p^r;s)$ has a pole at this point. The order of this pole is given by the order of vanishing of $\Lambda_{F,f,r, p}^{-1}(s)$ at this point which is 
	\begin{align*}
		\begin{cases}
		p-1 \ &\text{for} \ f=1,\\
			4 \ &\text{for} \ f=p=2,\\
			1 \ &\text{else}.
		\end{cases}
	\end{align*}
\end{proposition}

\begin{proof}
%It remains to show that the function $\psi_{F, f, r, p}(s)$ is holomorphic in the domain $\re(s) > \frac{1 - \frac{1}{p} + f(p-1)}{p^{r+1-f} (p^f - 1)}$ and does not vanish at the point $s = \frac{1+f(p-1)}{p^{r+1-f}\left(p^f-1\right)}$. \\
According to Proposition \ref{LocalMeromorphicContinuation}, we have
\begin{align*} %\label{Holomorphy}
	\psi_{F, f, r, p}(s) = \Lambda_{F, f,r, p}^{-1}(s) \cdot \prod\limits_{j=1}^{f} \prod\limits_{\idp \in \PP_F}  \delta_{j, \idp}^{-1}(s) \cdot \prod\limits_{\idp \in \PP_F} \Psi_{f, \idp}(s) \notag \\
	%= \quad & \prod\limits_{j=2}^{f-1} \zeta_{\Fqtg}(p^{r+1-j}(p^j -1)s - j(p-1)) \cdot \prod\limits_{\idp} \left(\delta^{\left(2^f - 2\right)}_{1,\idp}(s) \cdot \Psi_{f, \idp}(s)\right)
\end{align*}
where 
\begin{align*}
\Lambda_{F, f,r, p}^{-1}(s) \cdot \prod\limits_{j=1}^{f} \prod\limits_{\idp \in \PP_F}  \delta_{j, \idp}^{-1}(s) = \begin{cases}
	\prod\limits_{\ell=1}^{p-2} \zeta_{F}^{-1}(p^{r-1}(\ell+1)(p-1)s - \ell) \ &\text{for} \ f=1,\\
	\zeta_{F}^{-2}(2^{r}s - 1) \ &\text{for} \ f=p=2,\\
	\prod\limits_{j=1}^{f-1} \zeta_F(-\eta_j(s)) \ &\text{else}.
\end{cases}
\end{align*}
Note that we have $\prod\limits_{\idp \in \PP_F} \delta_{j, \idp}^{-1}(s) = \zeta_F(-\eta_j(s))$ by Definition \ref{LocalMeromorphicFunction}. The product $\prod\limits_{j=1}^{f-1} \zeta_F(-\eta_j(s))$ of zeta functions is holomorphic in the domain under consideration for $2 \leq f$ and $p \neq 2$ if $f=2$ by Lemma \ref{GlobalInequality}. Additionally, it does not have any zeros in our domain of interest by the Hasse-Weil Theorem (e.g.\ \cite[Theorem 5.10]{Ro}). Hence, it is sufficient to show the holomorphy of
\begin{align} \label{eq103}
\Upsilon_{f, p}(s) := 
\begin{cases}
	\prod\limits_{\idp \in \PP_F} \left(\prod\limits_{\ell=1}^{p-2} \left(1 - \np^{\ell - (p-1)p^{r-1}(\ell+1)s} \right) \cdot \Psi_{1, \idp}(s) \right) \ &\text{for} \ f=1,\\
	\prod\limits_{\idp \in \PP_F} \left( \left(1 - \np^{1 - 2^rs} \right)^2 \cdot \Psi_{2, \idp}(s) \right) \ &\text{for} \ f=p=2,\\
	\prod\limits_{\idp \in \PP_F} \Psi_{f, \idp}(s) \ &\text{else}.
\end{cases}
\end{align}
in the domain $\re(s) > \frac{1+f(p-1)}{p^{r+1-f}\left(p^f-1\right)} - \epsilon$ and the non-vanishing of $\Upsilon_{f, p}(s)$ at $s = \frac{1+f(p-1)}{p^{r+1-f}\left(p^f-1\right)}$. According to Definition \ref{Constants}, Definition \ref{LocalMeromorphicFunction} and Proposition \ref{LocalMeromorphicContinuation}, the Euler factors of $\Upsilon_{f, p}(s)$ have the shape $1 + G_{\idp}(\np^{-s})$ where $G_{\idp}(\np^{-s})$ is a finite sum of the form $\sum \gamma_i \np^{\upsilon_i(s)}$ with coefficients $\gamma_i \in \Z$ not depending on $\idp$ and linear functions $\upsilon_i(s)$ of $s$. 

It is well-known that an infinite product $\prod(1+a_n)$ is absolutely convergent if and only if the corresponding series $\sum a_n$ is absolutely convergent. Since the series $\sum\limits_{\idp \in \PP_F} \np^{-s}$ is absolutely convergent for $\re(s) > 1$, we have to check $\upsilon_i(s) < -1$ for all $\np$-exponents $\upsilon_i(s)$ of $G_{\idp}(\np^{-s})$ and all $\frac{1- \frac{1}{p} + f(p-1)}{p^{r+1-f}\left(p^f-1\right)} < s \in \R$. We have to carry out a case distinction for the cases given in equation \eqref{eq103}. 

\underline{Case 1: Let $f=1$.} By Definition \ref{Constants}, Definition \ref{LocalMeromorphicFunction} and Proposition \ref{LocalMeromorphicContinuation}, we have
\begin{align} \label{eq102}
& \prod\limits_{\ell=1}^{p-2} \left(1 - \np^{\ell - (p-1)p^{r-1}(\ell+1)s} \right) \cdot \Psi_{1, \idp}(s) \notag \\
%& = \prod\limits_{\ell=1}^{p-2} \left(1 - \np^{\ell - (p-1)p^{r-1}(\ell+1)s} \right) \cdot \left( 1 - \np^{p-1 -(p-1)p^{r}s} +  \sum\limits_{\ell=1}^{p-1} (\np-1) \np^{\ell-1 -(p-1)p^{r-1}(\ell +1)s} \right)
= &\prod\limits_{\ell=1}^{p-1} \left(1 - \np^{\ell - (p-1)p^{r-1}(\ell+1)s} \right) + \prod\limits_{\ell=1}^{p-2} \left(1 - \np^{\ell - (p-1)p^{r-1}(\ell+1)s} \right)  \sum\limits_{\ell=1}^{p-1} (\np-1) \np^{\ell-1 -(p-1)p^{r-1}(\ell +1)s}.
\end{align}
For $s \in \R$, the linear functions $\ell - (p-1)p^{r-1}(\ell+1)s$ are monotonically decreasing and satisfy $\ell - (p-1)p^{r-1}(\ell+1) s < -\frac{1}{2}$ for $\frac{\frac{1}{2}+p-1}{p^{r}(p-1)} < s$ and $1 \leq \ell \leq p-1$. Hence, we obtain \eqref{eq102}=
\begin{align*}
1 - \sum\limits_{\ell = 1}^{p-1} \np^{\ell - (p-1)p^{r-1}(\ell+1)s} + \sum\limits_{\ell = 1}^{p-1} \np^{\ell - (p-1)p^{r-1}(\ell+1)s} + O\left(\np^{-1 - \epsilon'} \right) = 1 + O\left(\np^{-1 - \epsilon'} \right)
\end{align*}
in the domain $\re(s) > \frac{1 - \frac{1}{p}+p-1}{p^{r}(p-1)}$ for some sufficiently small $\epsilon' >0$ depending on $s$ but not on $\idp$. In particular, the infinite product $\Upsilon_{1, p}(s)$ is absolutely convergent in this domain. Hence, the product does not vanish at $s = \frac{1 +p-1}{p^{r}(p-1)}$ if all Euler factors do not vanish at this point which is easy to verify by considering the evaluation of equation \eqref{eq102} at this point. 

\underline{Case 2: Let $f = p = 2$.} Again, by Definition \ref{Constants}, Definition \ref{LocalMeromorphicFunction} and Proposition \ref{LocalMeromorphicContinuation}, we have
\begin{align} \label{eq104}
	& \left(1 - \np^{1 - 2^r s} \right)^2 \cdot \Psi_{2, \idp}(s) \notag \\
	%& = \prod\limits_{\ell=1}^{p-2} \left(1 - \np^{\ell - (p-1)p^{r-1}(\ell+1)s} \right) \cdot \left( 1 - \np^{p-1 -(p-1)p^{r}s} +  \sum\limits_{\ell=1}^{p-1} (\np-1) \np^{\ell-1 -(p-1)p^{r-1}(\ell +1)s} \right)
	= \ & \left(1 - \np^{2 - 3 \cdot 2^{r-1} s} \right) \left(1 - \np^{1 - 2^r s} \right)^3 + 3 \cdot \left(1 - \np^{2 - 3 \cdot 2^{r-1} s} \right) \left(1 - \np^{1 - 2^r s} \right)^2 \cdot (\np -1) \cdot \np^{-2^rs}  \\
	\quad \quad & + 3 \cdot \left(1 - \np^{1 - 2^r s} \right)^2 \cdot \np^{1-2^rs} \cdot (\np-1)^2 \cdot \np^{-3 \cdot 2^{r-1}s} \notag \\
	\quad \quad & + \left(1 - \np^{1 - 2^r s} \right)^3 \cdot (\np-1) \cdot (\np-p) \cdot \np^{-3 \cdot 2^{r-1}s} \notag.
\end{align}
For $s \in \R$, the linear functions $1 - 2^r s$ and $2 - 3 \cdot 2^{r-1} s$ are monotonically decreasing and satisfy $1 - 2^r s < -\frac{2}{3}$ and $2 - 3 \cdot 2^{r-1} s < -\frac{1}{2}$ for $\frac{\frac{1}{2}+2}{3 \cdot 2^{r-1}} < s$. Hence, we obtain \eqref{eq104}=
\begin{align*}
	1 - \np^{2 - 3 \cdot 2^{r-1} s} - 3 \cdot \np^{1 - 2^{r} s} + 3 \cdot \np^{1 - 2^{r} s} +  \np^{2 - 3 \cdot 2^{r-1} s} + O\left(\np^{-1 - \epsilon'} \right) = 1 + O\left(\np^{-1 - \epsilon'} \right)
\end{align*}
in the domain $\re(s) > \frac{\frac{1}{2}+2}{3 \cdot 2^{r-1}}$ for some sufficiently small $\epsilon' >0$ depending on $s$ but not on $\idp$. Again, an easy evaluation of equation \eqref{eq104} at the point $s = \frac{1+2}{3 \cdot 2^{r-1}}$ yields the non-vanishing of $\Upsilon_{2, 2}(s)$ at this point. 

\underline{Case 3: Let $2 \leq f$ and $p \neq 2$ if $f=2$.} By Proposition \ref{LocalMeromorphicContinuation}, we have $\Psi_{f, \idp}(s) =$
\begin{align} \label{Eq1}
	\prod\limits_{j=1}^{f} \delta_{j, \idp}(s)
	+  \sum\limits_{h = 1}^{f} \binom{f}{h}_p \quad \sum\limits_{\theta \in \mathfrak{Z}^{ad}_h} \gamma_{\theta_{\omega}} \prod\limits_{\substack{j=1 \\
			j \not \in \{B_1, \ldots, B_{\lambda_{\rho}}\} }}^{f} \delta_{j, \idp}(s) \cdot \prod\limits_{j=1}^{\lambda_{\rho} - 1} \np^{\eta_{B_j}(s)} \cdot \psi_{\ell}(\theta).
\end{align}
For $s \in \R$, the linear functions $\eta_j(s) = j(p-1) -  p^{r+1-j} (p^j-1) s$ are monotonically decreasing and satisfy $\eta_j(s) < -1$ for $\frac{1+j(p-1)}{p^{r+1-j} (p^j-1)} < s$ and $1 \leq j \leq f$. By Lemma \ref{GlobalInequality}, we also have $\eta_j(s) < -1$ for $\frac{1 -\frac{1}{p}+ f(p-1)}{p^{r+1-f} (p^f-1)} < s$ and $1 \leq j \leq f-1$. Furthermore, we have $\eta_f(s) < - 1 + \frac{1}{p} \leq -\frac{1}{2}$ for $\frac{1 -\frac{1}{p}+ f(p-1)}{p^{r+1-f} (p^f-1)} < s$. Hence, we have
\begin{align} \label{eq105}
\prod\limits_{j=1}^{f} \delta_{j, \idp}(s) = 1 - \np^{\eta_f(s)} + O\left(\np^{-1 - \epsilon'} \right)
\end{align}
in the domain $\re(s) > \frac{1 -\frac{1}{p}+ f(p-1)}{p^{r+1-f} (p^f-1)}$ for some sufficiently small $\epsilon' > 0$ depending on $s$ but not on $\idp$. 

Next, we consider the second summand in equation \eqref{Eq1} and fix a $\theta \in \mathfrak{Z}^{ad}_h$. The largest $\np$-exponents appear for terms of the form $\prod\limits_{j=1}^{\lambda_{\rho} - 1} \np^{\eta_{B_j}(s)} \cdot \np^h \cdot \prod\limits_{i=1}^{\Lambda_{\lambda_{\rho}}} \np^{\beta_i(s) }$ by Definition \ref{Constants} and the fact that $g_{\theta}(\np)$ has degree $h$ by Definition \ref{2-level-composition}. For the exponent, we obtain
\begin{align} \label{Eq2}
	\sum\limits_{j=1}^{\lambda_{\rho} -1} \sum\limits_{i=1}^{j} \left(b_i(p-1) - (p-1) \sum\limits_{m=1+B_{i-1}}^{B_i} p^{r+1-m} s\right) + \sum\limits_{i=1}^{\Lambda_{\lambda_{\rho}}} \left(a_i \ell_{A_i} -(p-1) \sum\limits_{j=1+A_{i-1}}^{A_i} p^{r-j} (\ell_{A_i} + 1)s\right)
\end{align}
by applying the definition of $\eta_{B_j}(s)$ in Definition \ref{LocalMeromorphicFunction}, the definition of $\beta_i(s)$ in Definition \ref{Constants} and using the identity $h = \sum\limits_{i=1}^{\Lambda_{\lambda_{\rho}}} a_i$. We can simplify the notation in the second part of the exponent by applying a new indexing. In fact, we consider $h$ different $\ell$-values $\ell_1, \ldots, \ell_h$ such that $\ell_j = \ell_{A_i}$ for $A_{i-1} < j \leq A_i$ and all $1 \leq i \leq \Lambda_{\lambda_{\rho}}$. Then, we have \eqref{Eq2} = 
\begin{align} \label{eq500}
	\sum\limits_{i=1}^{h} \ell_i + (p-1) \sum\limits_{i=1}^{\lambda_{\rho} - 1} (\lambda_{\rho} -  i) b_i - (p-1) \sum\limits_{i=1}^{h} p^{r-i} (\ell_i + 1)s - (p-1) \sum\limits_{i=1}^{\lambda_{\rho} - 1} (\lambda_{\rho} -  i) \sum\limits_{j=1 + B_{i-1}}^{B_i} p^{r+1-j}s .
\end{align}

If the real part of this $\np$-exponent reaches the value $-1$, we have 
\begin{align} \label{eq600}
	\re(s) = \frac{1 + \sum\limits_{i=1}^{h} \ell_i + (p-1) \sum\limits_{i=1}^{\lambda_{\rho} - 1} (\lambda_{\rho} -  i) b_i}{(p-1) \sum\limits_{i=1}^{h} p^{r-i} (\ell_i + 1) + (p-1) \sum\limits_{i=1}^{\lambda_{\rho} - 1} (\lambda_{\rho} -  i) \sum\limits_{j=1 + B_{i-1}}^{B_i} p^{r+1-j}}.
\end{align}
Note that equation \eqref{eq500} also describes a monotonically decreasing function for $s \in \R$.
Hence, it is sufficient to compare all fractions of this form with $\frac{1 - \frac{1}{p} + f(p-1)}{p^{r+1-f}(p^f - 1)}$. If a fraction of the form given in equation \eqref{eq600} is smaller than $\frac{1 - \frac{1}{p} + f(p-1)}{p^{r+1-f}(p^f - 1)}$, we obtain the holomorphy for the corresponding term in equation \eqref{Eq1} in the desired domain. We are going to carry out the computations in the following Lemmas \ref{Lemma1} and \ref{Lemma2}. 

Note that the case $h = 1$ follows directly by our considerations for the case $f=1$ and Lemma \ref{GlobalInequality}. For $1 < h < f$, we can directly apply the Lemmas \ref{Lemma1} and \ref{Lemma2}. By Lemma \ref{GlobalInequality}, we obtain the desired comparison with $\frac{1 - \frac{1}{p} + f(p-1)}{p^{r+1-f}(p^f - 1)}$.
Furthermore, the case $h = f$ is also directly covered by the Lemmas and we find that only the fraction in the case $\rho = (f)$ with $\ell_1 = \ldots = \ell_f = p-1$ exceeds $\frac{1 - \frac{1}{p} + f(p-1)}{p^{r+1-f}(p^f - 1)}$. In this case, the summand is $\np^{f(p-1) - p^{r+1-f}(p^f-1)s} = \np^{\eta_f(s)}$, which has the coefficient $\binom{f}{f}_p \cdot \gamma_{\theta_{\omega}} = 1$. Consequently, this summand is canceled by its additive inverse in equation \eqref{eq105} which yields
\begin{align*}
\Psi_{f, \idp}(s) = 1 + O\left(\np^{-1 - \epsilon'} \right)
\end{align*}
in the domain $\re(s) > \frac{1 - \frac{1}{p} + f(p-1)}{p^{r+1-f}(p^f - 1)}$ for some sufficiently small $\epsilon' >0$ depending on $s$ but not on $\idp$. 

It remains to show the non-vanishing of $\psi_{F,f, r, p}(s)$ at $s = \frac{1 + f(p-1)}{p^{r+1-f}(p^f - 1)}$. Due to the non-vanishing of $\prod\limits_{j=1}^{f-1} \zeta_F(-\eta_j(s))$ at this point which we discussed earlier, we only have to check the non-vanishing of $\Psi_{f, \idp}(s)$ at this point by absolute convergence. However, this is guaranteed by Proposition \ref{LocalMeromorphicContinuation}.
\end{proof}

\begin{lemma} \label{Lemma1}
	Let $2 \leq h \leq r$. We consider a 2-level composition with the outer composition $\rho = (h)$. Accordingly, we have set the restriction $p-1 \geq \ell_1 \geq \ldots \geq \ell_h \geq 1$. Then, the inequality
	\begin{align} \label{eq300}
		\frac{1 + \sum\limits_{i=1}^{h} \ell_i}{(p-1) \sum\limits_{i=1}^{h} p^{r-i} (\ell_i+1)} < \frac{1 - \frac{1}{p} + h(p-1)}{p^{r+1-h}(p^h - 1)}
	\end{align}
	applies except for the case $\ell_1 = \ldots = \ell_h = p-1$ in which we have
	\begin{align*}
	\frac{1 + \sum\limits_{i=1}^{h} \ell_i}{(p-1) \sum\limits_{i=1}^{h} p^{r-i} (\ell_i+1)} = \frac{1+h(p-1)}{p^{r+1-h}(p^h-1)}.
	\end{align*}
\end{lemma}

\begin{proof}
	We multiply equation \eqref{eq300} with the factor $p^{r-h}$ and we have
	\begin{align} \label{Ineq}
		& \frac{1 + \sum\limits_{i=1}^{h} \ell_i}{(p-1) \sum\limits_{i=1}^{h} p^{h-i} (\ell_i+1)} < \frac{1 - \frac{1}{p} + h(p-1)}{p(p^h - 1)} \notag \\
		\Longleftrightarrow \quad & 0 < \left(1 - \frac{1}{p} + h (p-1) \right) \sum\limits_{i=1}^{h} p^{h-i} (\ell_i + 1) - \left(1 + \sum\limits_{i=1}^{h} \ell_i \right) \sum\limits_{i=0}^{h-1} p^{h-i} \notag \\
		\Longleftrightarrow \quad & 0 < p^h \left(h-1 + \sum\limits_{j=1}^{h} (\ell_1 - \ell_j)\right) + \sum\limits_{i=1}^{h-1} p^{h-i} \left(h(\ell_{i+1} - \ell_i) - \sum\limits_{\substack{j=1 \\ j \neq i}}^{h} \ell_j \right) \\
		& \quad \quad \quad - \left(h - 1 \right) (\ell_h + 1) - \sum\limits_{i=1}^{h} p^{h-1-i} (\ell_i + 1) \notag.
	\end{align}

	By induction on $h$, we show the inequality 
	\begin{align*}
		0 \leq p^h \sum\limits_{j=1}^{h} (\ell_1 - \ell_j) + \sum\limits_{i=1}^{h-1} p^{h-i} h(\ell_{i+1} - \ell_i).
	\end{align*}
	In the initial case $h=2$, we find $0 \leq (\ell_1 - \ell_2) p (p - 2)$ which holds. 
	%Since we have
	%\begin{align*}
	%&p^r \sum\limits_{j=1}^{r} (\ell_1 - \ell_j) + \sum\limits_{i=1}^{r-1} p^{r-i} r(\ell_{i+1} - \ell_i) \\
	% = & p^r \sum\limits_{j=1}^{r-1} (\ell_1 - \ell_j) + \sum\limits_{i=1}^{r-2} p^{r-i} (r-1) (\ell_{i+1} - \ell_i) + p^r (\ell_1 - \ell_r) + \sum\limits_{i=1}^{r-2} p^{r-i} (\ell_{i+1} - \ell_i) + r p (\ell_r - \ell_{r-1}),
	%\end{align*}
	By the induction hypothesis, it is sufficient to show $0 \leq p^h (\ell_1 - \ell_h) + \sum\limits_{i=1}^{h-2} p^{h-i} (\ell_{i+1} - \ell_i) + h p (\ell_h - \ell_{h-1})$. This inequality, however, holds since we have $p^h > p^{h-i}$ for $1 \leq i$ and $p^h \geq hp$ for $2 \leq h$ and we assumed $\ell_i \geq \ell_{i+1}$ for $1 \leq i \leq h-1$. 
	
	By equation \eqref{Ineq}, it remains to check
	\begin{align*}
		0 < p^h (h-1) - \sum\limits_{i=1}^{h-1} p^{h-i} \sum\limits_{\substack{j=1 \\ j \neq i}}^{h} \ell_j - (h-1)(\ell_h + 1) - \sum\limits_{i=1}^{h} p^{h-1-i} (\ell_i + 1).
	\end{align*}
	The right-hand side is minimised in the case $\ell_1 = \ldots = \ell_{h-1} = p-1, \ell_h = p-2$. We obtain
	\begin{align*}
		&p^h (h-1) - \sum\limits_{i=1}^{h-1} p^{h-i} (h-1)(p-1) + \sum\limits_{i=1}^{h-1} p^{h-i} - (h-1)(p-1) - \sum\limits_{i=1}^{h} p^{h-i} + \frac{1}{p} \\
		= \  &h - 2 + \frac{1}{p} > 0.
	\end{align*}
	in this case which concludes the proof.
\end{proof}

\begin{lemma} \label{Lemma2}
Let $2 \leq h \leq r$. We consider a 2-level composition of $h$ with the outer composition $\rho = (b_1, \ldots, b_{\lambda_{\rho}})$ and $\lambda_{\rho} \geq 2$ and its related restrictions to the $\ell$-values. Then, the following inequality applies:
	\begin{align} \label{Ineq2Blocks}
		\frac{1 + \sum\limits_{i=1}^{h} \ell_i + (p-1) \sum\limits_{i=1}^{\lambda_{\rho} -1 } (\lambda_{\rho} -i) b_i}{(p-1) \sum\limits_{i=1}^{h} p^{r-i} (\ell_i+1) + (p-1) \sum\limits_{i=1}^{\lambda_{\rho} -1 } (\lambda_{\rho} -i) \sum\limits_{j=1 + B_{i-1}}^{B_{i}} p^{r+1-j}} < \frac{1 - \frac{1}{p} + h(p-1)}{p^{r+1-h}(p^h - 1)}.
	\end{align}
\end{lemma}

\begin{proof}
	Initially, we multiply the fractions by the factor $p^{r-h}$ which yields
	\begin{align} \label{Ineq5}
		\frac{1 + \sum\limits_{i=1}^{h} \ell_i + (p-1) \sum\limits_{i=1}^{\lambda_{\rho} -1 } (\lambda_{\rho} -i) b_i}{(p-1) \sum\limits_{i=1}^{h} p^{h-i} (\ell_i+1) + (p-1) \sum\limits_{i=1}^{\lambda_{\rho} -1 } (\lambda_{\rho} -i) \sum\limits_{j=1 + B_{i-1}}^{B_{i}} p^{h+1-j}} < \frac{1 - \frac{1}{p} + h(p-1)}{p(p^h - 1)}.
	\end{align}

	Our approach is to prove equation \eqref{Ineq5} by induction on $h$. In fact, we have to distinguish two different cases in the induction step. A composition of $h+1$ can arise from a composition of $h$ in two ways. On the one hand we can extend the last block of a composition of $h$ by $1$ whereby the number of blocks does not change and on the other hand we can add a new block of length $1$. We prove the former case first and use it for the remaining one. 
	
	In the initial case $h=2$, we have to check the inequality for the composition $(1,1)$ which leads to
	\begin{align*}
		\frac{1 + \ell_1 + \ell_2 + p-1}{(p-1) \left(p (\ell_1 + 1) + \ell_2 + 1 \right) + (p-1) p^{2}} < \frac{1 - \frac{1}{p} + 2(p-1)}{p(p^2 - 1)}.
	\end{align*}
	By differentiation with respect to $\ell_1$ and $\ell_2$, we find that the fraction on the left-hand side is maximised for $\ell_1 = 1$ and $\ell_2 = p-1$. Accordingly, it is sufficient to check
	\begin{align*}
		& \frac{2p}{(p-1) (p^2 + 3p)} < \frac{1 - \frac{1}{p} + 2(p-1)}{p(p^2 - 1)} \\
		\Longleftrightarrow \quad & 2p^2 + 2p + \left(1 + \frac{1}{p}\right) (p+3) < 2p^2 + 6p \\
		\Longleftrightarrow \quad & 3 + 4p <  3p^2
	\end{align*}
	which holds because of $2p \leq p^2$ and $3 < p^2$. 
	
	Next, we carry out the induction step in the case in which the number of blocks of the composition does not increase. To abbreviate the notation, we write $\frac{A}{B} < \frac{C}{D}$ for \eqref{Ineq5}. Therefore, we have to prove the inequality
	\begin{align} \label{Ineq3}
		\frac{A + \ell_{h+1}}{pB + (p-1)(\ell_{h+1} +1)} < \frac{C + p -1}{pD + (p-1)p}.
	\end{align}
	Differentiating with respect to $\ell_{h+1}$ reveals that the fraction on the left-hand side is maximised for $\ell_{h+1} = p-1$. By multiplying with the denominators and using the induction hypothesis, it remains to check
	\begin{align} \label{Ineq2}
		& (p-1)p A + (p-1)(pD + (p-1)p) < (p-1)p B + (p-1)p (C + p-1) \notag \\
		\Longleftrightarrow \quad & A + D < B + C \notag \\
		\Longleftrightarrow \quad & 1 + \sum\limits_{i=1}^{h} \ell_i + p^{h+1} < 1 + h(p-1) - \frac{1}{p}  + p -1 + p^h + (p-1) \sum\limits_{i=1}^{h} p^{h-i} \ell_i \\
		& \quad \quad \quad \quad \quad \quad \quad \quad \quad \quad + (p-1) \sum\limits_{i=1}^{\lambda_{\rho} -1 } (\lambda_{\rho} -i) \sum\limits_{j=1 + B_{i-1}}^{B_{i}} (p^{h+1-j}-1). \notag
	\end{align}
	The right-hand side is minimised in the case $\lambda_{\rho} = 2$ and $b_1 = 1$. Therefore, we have $(p-1) \sum\limits_{i=1}^{\lambda_{\rho} -1 } (\lambda_{\rho} -1) \sum\limits_{j=1 + B_{i-1}}^{B_{i}} (p^{h+1-j}-1) \geq (p-1) (p^{h} - 1)$. To prove \eqref{Ineq2}, it is hence sufficient to check
	\begin{align*}
		& \sum\limits_{i=1}^{h} \ell_i + p^{h+1} < h(p-1) + p-1 - \frac{1}{p} + p^h + (p-1) \sum\limits_{i=1}^{h} p^{h-i} \ell_i + (p-1) (p^h - 1) \\
		\Longleftrightarrow \quad & \sum\limits_{i=1}^{h} \ell_i  < h(p-1) - \frac{1}{p}  + (p-1) \sum\limits_{i=1}^{h} p^{h-i} \ell_i
	\end{align*}
	which holds due to $\sum\limits_{i=1}^{h} \ell_i \leq h (p-1)$ and $\frac{1}{p} < (p-1) \sum\limits_{i=1}^{h} p^{h-i} \ell_i$. Accordingly, equation \eqref{Ineq3} holds which completes the first case. 
	
	Finally, we consider the case in which we add a new block of length $1$. We have to prove the inequality
	\begin{align*}
		\frac{A + \ell_{h+1} + (p-1) \sum\limits_{i=1}^{\lambda_{\rho}} b_i}{pB + (p-1)(\ell_{h+1} +1) + (p-1) \sum\limits_{i=1}^{\lambda_{\rho}} \sum\limits_{j=1 + B_{i-1}}^{B_i} p^{h+2-j} } < \frac{C + p -1}{pD + (p-1)p}.
	\end{align*}
	With respect to \eqref{Ineq3}, it remains to verify
	\begin{align} \label{Ineq4}
		(p-1) \sum\limits_{i=1}^{\lambda_{\rho}} b_i \cdot (pD + (p-1)p) < (p-1) \sum\limits_{i=1}^{\lambda_{\rho}} \sum\limits_{j=1 + B_{i-1}}^{B_i} p^{h+2-j} \cdot (C+ p-1).
	\end{align}
	As we have $\sum\limits_{i=1}^{\lambda_{\rho}} b_i = h$ and $\sum\limits_{i=1}^{\lambda_{\rho}} \sum\limits_{j=1 + B_{i-1}}^{B_i} p^{h+2-j} = \sum\limits_{j=1}^{h} p^{h+2-j}$, equation \eqref{Ineq4} is equivalent to
	\begin{align*}
		& h (p-1) \sum\limits_{i=1}^{h+1} p^{h+2-i} < \sum\limits_{j=1}^{h} p^{h+2-j} \cdot \left(1 - \frac{1}{p} + (h+1)(p-1)\right) \\
		\Longleftrightarrow \quad & h (p-1) p < \left( p - \frac{1}{p} \right) \sum\limits_{j=1}^{h} p^{h+2-j}.
	\end{align*}
	Due to $p^{h+2-j} \geq p^2$, we have $\left( p - \frac{1}{p} \right) \sum\limits_{j=1}^{h} p^{h+2-j} \geq h (p^2-1)p > h(p-1)p$.
\end{proof}

%\subsection{Meromorphic continuation of the Dirichlet series of $\Fqtg$} \hfill \\

\begin{theorem} \label{ContinuationDirichletSeries}
	The complex function 
	\begin{align*}
		\Psi_{F,r,p}(s) := \left( \sum\limits_{f=0}^r e_f \cdot \Phi_{f}(F, C_p^r; s) \right) \cdot \Lambda_{F,r,r,p}^{-1}(s)
	\end{align*}
	has a holomorphic continuation to the domain $\re(s) > \frac{1+r(p-1)}{p \left(p^r-1\right)} - \epsilon$ with $\epsilon = \frac{1}{p^{2} (p^r -1)} > 0$.
	
	Furthermore, $\Psi_{F,r,p}(s)$ does not vanish at $s = \frac{1+r(p-1)}{p\left(p^r-1\right)}$, so $\sum\limits_{f=0}^r e_f \cdot \Phi_{f}(F, C_p^r; s)$ has a pole at this point. The order of this pole is given by the order of vanishing of $\Lambda_{F,r,r, p}^{-1}(s)$ at this point which is $b_p(C_p^r)$ defined in Theorem \ref{MainResult}.
%	\begin{align} \label{PoleOrder}
%		\begin{cases}
%			p-1 \ &\text{for} \ r=1,\\
%			4 \ &\text{for} \ r=p=2,\\
%			1 \ &\text{else}.
%		\end{cases}
%	\end{align}

In particular, this Theorem applies to the Dirichlet series $\Phi(\Fqtg, C_p^r; s)$ by Proposition \ref{LCEuler}.
\end{theorem}

\begin{proof}
%By Proposition \ref{LCEuler}, we have $	\Phi(\Fqtg, C_p^r;s) = \sum\limits_{f=0}^{r} e_f \cdot \Phi_f(\Fqtg, C_p^r;s)$. \\
By Proposition \ref{ContinuationEulerProduct}, we obtain
\begin{align} \label{MC}
\Psi_{F,r,p}(s) = e_r \cdot \psi_{F,r,r,p}(s) + \Lambda^{-1}_{F,r,r,p}(s) \cdot \sum\limits_{f=1}^{r-1} e_f \cdot \Lambda_{F,f,r,p}(s) \cdot \psi_{F,f,r,p}(s) + \Lambda^{-1}_{F,r,r,p}(s) \cdot e_0.
\end{align}
It remains to prove the holomorphy of the function $\Psi_{F,r,p}(s)$ in the suggested domain and the non-vanishing of this function at the point  $s = \frac{1+r(p-1)}{p\left(p^r-1\right)}$. Proposition \ref{ContinuationEulerProduct} directly yields the holomorphy for $e_r \cdot \psi_{F,r,r,p}(s)$. 

By Proposition \ref{LambdaZeros}, $\Lambda_{F,r,r,p}(s)$ does not have any zeros in the domain $\re(s) > \frac{1 - \frac{1}{p} +r(p-1)}{p \left(p^r-1\right)}$. Hence, $\Lambda^{-1}_{F,r,r,p}(s)$ is holomorphic in this domain. In particular, we are done with the case $r=1$. 

In the case $r=p=2$, we have 
\begin{align*}
\Lambda^{-1}_{F,2,2,2}(s) \cdot e_1 \cdot \Lambda_{F,1,2,2}(s) \cdot \psi_{F,1,2,2}(s) = \zeta^{-1}_{F}(6s - 2) \cdot \zeta^{-2}_{F}(4s - 1) \cdot e_1 \cdot \psi_{F,1,2,2}(s).
\end{align*}
Again, $\zeta^{-1}_{F}(6s - 2) \cdot \zeta^{-2}_{F}(4s - 1)$  does not have poles in the considered domain and $\psi_{F,1,2,2}(s)$ is holomorphic for $\re(s) > \frac{1}{2} - \frac{1}{8} = \frac{3}{8}$ by Proposition \ref{ContinuationEulerProduct}. Due to $\frac{3}{8} < \frac{5}{12} = \frac{1 - \frac{1}{2} + 2}{2(2^2-1)}$, we get the claim in this case. 

Let $r \geq 2$, $1 \leq f \leq r-1$ and $p \neq 2$ if $r=2$. Since the Euler product $\Phi_f(F, C_p^r;s)$ has convergence abscissa $\frac{1+f(p-1)}{p^{r+1-f}(p^f -1)}$ by Proposition \ref{GlobalPoles} and Proposition \ref{ContinuationEulerProduct}, we are done by Lemma \ref{GlobalInequality}. 

Finally, we consider equation \eqref{MC} and we have 
\begin{align*}
	\Psi_{F,r,p}\left(\frac{1+r(p-1)}{p\left(p^r-1\right)} \right) = e_r \cdot \psi_{F,r,r,p}\left(\frac{1+r(p-1)}{p\left(p^r-1\right)} \right) \neq 0.
\end{align*}
The first equality follows by $\Lambda^{-1}_{F,r,r,p}\left(\frac{1+r(p-1)}{p \left(p^r-1\right)} \right) = 0$ and the holomorphy of $\Phi_{f}(F, C_p^r; s)$ in the domain $\re(s) > \frac{1+f(p-1)}{p^{r+1-f}(p^f -1)}$. Finally, the non-vanishing is ensured by Proposition \ref{ContinuationEulerProduct} and the claim about the pole order also follows from this Proposition.
\end{proof}

By Theorem \ref{ContinuationDirichletSeries}, all poles of $\Lambda_{\Fqtg,r,r,p}(s)$ on the critical line with $\re(s) = \frac{1 + r(p-1)}{p(p^r -1)}$ are candidates for poles of the Dirichlet series $\Phi(\Fqtg, C_p^r;s)$ on its convergence abscissa. Whether $\Phi(\Fqtg, C_p^r;s)$ actually has a pole or is holomorphically continuable depends on the vanishing order of $\Psi_{\Fqtg, r,p}(s)$ at the respective point. At least, we have already proven that the rightmost candidate on the real line is definitely a pole of $\Phi(\Fqtg, C_p^r;s)$ whose order is $b_p(C_p^r)$ defined in Theorem \ref{MainResult}.

\begin{remark} \label{LocalConfig}
In Proposition \ref{ContinuationEulerProduct}, Lemma \ref{Lemma1} and Lemma \ref{Lemma2}, we have proven that the rightmost poles of $\sum\limits_{f=0}^r e_f \cdot \Phi_{f}(F, C_p^r; s)$ come from the summand $e_r \cdot \Phi_{r}(F, C_p^r; s)$ and are encoded in the outer composition $(1)$ with $1 \leq \ell \leq p-1$ in the case $r = 1$, in the outer compositions $(1)$ with $\ell_1 = 1$ and $(2)$ with $\ell_1 = \ell_2 = 1$ in the case $r = p = 2$ and in the outer composition $(r)$ with $\ell_1 = \ldots = \ell_r = p-1$ in the remaining case. 

Now, we recall the equations \eqref{EinsetzenConfig} - \eqref{eq1} from our local computations and we want to determine explicitly which (local) configurations of conductors give these rightmost poles. 

For $r=1$, we obtain $\frac{1 + \ell + (p-1)k}{(p-1)(pk+\ell+1)}$ for the real part of the global poles corresponding to the local conductor $c = pk+\ell+1$ with $k \geq 0$ and $1 \leq \ell \leq p-1$ by equation \eqref{SideCalculation1}. Clearly, this fraction is maximised for $k=0$ which gives $\frac{1 + \ell + (p-1)k}{(p-1)(pk+\ell+1)} = \frac{1 + \ell}{(p-1)(\ell+1)} = \frac{1}{p-1}$ for each $1 \leq \ell \leq p-1$. Hence, the rightmost pole of order $p-1$ originates from the totally ramified local $C_p$-extensions with conductors $2, \ldots, p$ where each of these local conductors contributes one pole order. 

Now, let $r \geq 2$ and we consider the outer composition $(r)$ with $\ell_1 = \ldots = \ell_r = p-1$. These conditions lead to the local configurations satisfying $c_1 = \ldots = c_r = pk + p$ for any $k \geq 0$. By equation \eqref{SideCalculation1}, the real part of the global poles corresponding to this configuration is $\frac{1+r(p-1)+r(p-1)k}{\left(p^r-1\right)(pk+p)}$. Due to $\frac{1+r(p-1)+r(p-1)k}{\left(p^r-1\right)(pk+p)} = \frac{\frac{1}{k+1} +r(p-1)}{p\left(p^r-1\right)}$, this fraction is again maximised for $k=0$. Hence, excluding the case $r=p=2$, the local configuration where all local conductors are equal to $p$ and the local Galois group is isomorphic to $C_p^r$ yields the rightmost simple poles of $e_r \cdot \Phi_{r}(F, C_p^r; s)$. 

In the case $r=p=2$, we have to combine the results from the case $r=1$ for the outer composition $(1)$ and the results from the case $r \geq 2$ for the outer composition $(2)$. Hence, the rightmost poles of $e_2 \cdot \Phi_{2}(F, C_2^2; s)$ of order $4$ also originate from the local configurations where all local conductors are equal to $2$. However, pole order $3$ is contributed by the local $C_2$-algebras and pole order $1$ is contributed by the totally ramified local $C_2^2$-extensions with conductor $2$. 

To put it another way, the rightmost global poles of $\sum\limits_{f=0}^r e_f \cdot \Phi_{f}(F, C_p^r; s)$ are shifted further to the left if we only consider global extensions which do not have the previously discussed optimal local configurations of the respective case at any place $\idp \in \PP_{\Fqtg}$.

\end{remark}

\subsection{Asymptotics of elementary-abelian extensions of $\Fqtg$} \label{ProofSpecialCase} \hfill\\
%\begin{theorem} \label{GlobalAsymptotics}
%	Let $a_p(C_p^r) := \frac{1+r(p-1)}{p(p^r -1)}$, $b_p(C_p^r) := 
%	\begin{cases}
%		p-1 & \text{for } r=1, \\
%		4 & \text{for } r=p=2,\\
%		1 & \text{else,}
%	\end{cases}$ and \\
%$L(C_p^r) := \begin{cases}
%	(p-1)M & \text{for } r=1, \\
%	12 & \text{for } r=p=2, \\
%	p(p^r - 1) & \text{else,}
%\end{cases}$ \quad where $M$ is the least common multiple of the numbers $2, \ldots, p$ defined in Proposition \ref{GlobalPoles}. Furthermore, let $X = q^m$ and $m \to \infty$.
%Then, there exists a real polynomial $P_{m \bmod L(C_p^r)}$ of degree at most $b_p(C_p^r) - 1$ depending on the residue class of $m$ modulo $L(C_p^r)$ such that the number of $C_p^r$-extensions of a rational function field $\Fqtg$ when counted by discriminants satisfies
%\begin{align} \label{GA}
%	Z(\Fqtg, C_p^r; X) = P_{m \bmod L(C_p^r)}(\log(X)) \cdot X^{a_p(C_p^r)} + O\left(X^{\frac{1 - \frac{1}{p} + r(p-1)}{p(p^{r} -1)} + \epsilon}\right)
%\end{align}
%for every $\epsilon > 0$. The polynomial $P_{m \bmod L(C_p^r)}$ has degree $b_p(C_p^r) - 1$ for at least one residue class and its leading coefficient is always positive if $P_{m \bmod L(C_p^r)}$ has this degree. 
%% and depends on the period $L_{\text{max}}(C_p^r) := \begin{cases}
%%1 & \text{for } r=1 \ \text{and} \ p \neq 2, \\
%%2 & \text{for } r \in \{1, 2\} \ \text{and} \ p=2, \\
%%p(p^r - 1) & \text{else.}
%%\end{cases}$
%\end{theorem}
\vspace{-0.5cm}
\begin{proof}[Proof of Theorem \ref{MainResult} for rational function fields $\Fqtg$] \hfill \\
%We proceed as in the proof of Theorem \ref{LocalAsymptotics}. However, we have not proven that $\Phi(\Fqtg, C_p^r;s)$ is a rational function in $q^{-s}$. Hence, we want to apply \cite[Theorem IV.10]{FS}. \\
Substituting $z = q^{-s}$, we can transform $\Phi(\Fqtg, C_p^r;s)$ into a power series in the variable $z$. By Theorem \ref{ContinuationDirichletSeries} and Proposition \ref{GlobalPoles}, $\Phi(\Fqtg, C_p^r;z) = \sum\limits_{m \geq 0} c_m z^m$ has radius of convergence $R = q^{- \frac{1+r(p-1)}{p(p^r - 1)}}$ and is meromorphically continuable beyond its circle of convergence with a discrete set of poles on the circle of convergence. By Proposition \ref{GlobalPoles}, the candidates for poles on the circle of convergence are given by $\alpha_j = R \xi^{-j}$ for some primitive $L(C_p^r)$-th root of unity $\xi$ and $1 \leq j \leq L(C_p^r)$. %The poles of maximal order $b_p(C_p^r)$ on the circle of convergence are given by $R\zeta^{-j}$ for some $L_{\text{max}}(C_p^r)-$th root of unity $\zeta$ and $1 \leq j \leq L_{\text{max}}(C_p^r)$ where $L_{\text{max}}(C_p^r) := \begin{cases}
	%1 & \text{for } r=1 \ \text{and} \ p \neq 2, \\
	%2 & \text{for } r \in \{1, 2\} \ \text{and} \ p=2, \\
	%p(p^r - 1) & \text{else. }
%\end{cases}$ 

Now, we can apply \cite[Theorem IV.10]{FS} to obtain
\begin{align*}
c_m = \sum\limits_{j=1}^{L(C_p^r)} \Pi_j(m) \alpha_j^{-m} + O\left(X^{\frac{1 - \frac{1}{p} + r(p-1)}{p(p^{r} -1)} + \epsilon}\right)
\end{align*}
for every $\epsilon > 0$ and a complex polynomial $\Pi_j$ for each $1 \leq j \leq L(C_p^r)$ which is the zero polynomial if $\Phi(\Fqtg, C_p^r;z)$ is holomorphic at $R \xi^{-j}$ and whose degree equals the order of the pole at $R \xi^{-j}$ minus one otherwise. For the main term, we compute
\begin{align*}
\sum\limits_{j=1}^{L(C_p^r)} \Pi_j(m) \alpha_j^{-m} = \sum\limits_{j=1}^{L(C_p^r)} \Pi_j(m) \left(R \xi^{-j} \right)^{-m} = \left(\sum\limits_{j=1}^{L(C_p^r)} \Pi_j(m) \xi^{jm}\right) \cdot X^{a_p(C_p^r)}.
\end{align*}
Note that we have $X = q^{m}$ and thus $m = \frac{\log(X)}{\log(q)}$. Hence, we obtain equation \eqref{GA} for a complex polynomial $P_{m \bmod L(C_p^r)}$ of degree at most $b_p(C_p^r)-1$ by Proposition \ref{GlobalPoles}. 

Clearly, the sum depends on the residue class of $m \bmod L(C_p^r)$ since $\xi$ is a primitive $L(C_p^r)$-th root of unity. Additionally, the coefficients of the polynomial $P_{m \bmod L(C_p^r)}$ are actually real since the $c_m$ are non-negative integers. As $\Phi(\Fqtg, C_p^r;z)$ actually has pole of degree $b_p(C_p^r)$ at $z=R$ by Theorem \ref{ContinuationDirichletSeries}, the polynomial has degree $b_p(C_p^r)-1$ for at least one congruence class modulo $L(C_p^r)$. 
Finally, the leading term of $P_{m \bmod L(C_p^r)}$ is positive if $P_{m \bmod L(C_p^r)}$ has degree $b_p(C_p^r)-1$, again because the $c_m$ are non-negative integers.
\end{proof}

\begin{remark}
We recall Remark \ref{MTLR} and note that \cite[Theorem IV.10]{FS} again gives a better error bound than \cite[Theorem A.5]{La2}. 
The estimate for the error term in Theorem \ref{MainResult} in the case $F = \Fqtg$ can even be improved by carrying out case distinctions depending on $p$ and $r$ and precisely calculating the radius of the circle containing the subsequent poles. The bound given above, however, is valid uniformly for all $p$ and $r$. 

For a remark on the leading term $c(\Fqtg, C_p^r, m \bmod L_{\text{max}}(C_p^r))$ of $P_{m \bmod L(C_p^r)}$ (that is the coefficient corresponding to the $b_p(C_p^r)-1$-th power), we refer to Remark \ref{C-Constant} where a statement for arbitrary global function fields $F$ as base fields is made. Note that $L_{\text{max}}(C_p^r)$ is the period of the poles of order exactly $b_p(C_p^r)$ which is also defined in Remark \ref{C-Constant}.
\end{remark}

\section{Arbitrary base fields} \label{ArbitraryBaseFields}
From now on, we drop the restriction on the base field and admit any global function field $F$ of characteristic $p$, whose field of constant is $\F_q$, as our base field. In this setting, we have to deal with a potentially non-trivial class group $\Pic^{0}_F$. In fact, $\Pic^{0}_F$ is closely connected to the $p$-Selmer group 
\begin{align*}
	S_p := \{x \in F^{\times} \ : \ (x) = \mathfrak{a}^p \ \text{for some divisor} \ \mathfrak{a}\}/\left(F^{\times}\right)^p
\end{align*}
as we have $\Pic^{0}_F[p] \cong S_p$ for its $p$-torsion. 

By the Riemann-Roch theorem for curves, it can be shown that $\Pic^{0}_F$ is a finite abelian group which is generated by all places of degree $\leq g_F$ (e.g.\ see \cite[Proposition 5.1.3]{St} and its proof). Additionally, $\Pic_F$ is isomorphic to $\Z \times \Pic^{0}_F$ and the $\Z$-component can be generated by any divisor of degree $1$. %If such a divisor is not already generated by the set of places generating $\Pic^{0}_F$, we have to add at most two places with a degree larger than $g_F$. This is a direct consequence of the Hasse-Weil bound (\cite[Theorem 5.2.3]{St}) as we have $\F_{q^n}$-points for all sufficiently large $n$. Thereby, we note that a constant field extension does not change the genus (\cite[Theorem 3.6.3 (b)]{St}). In particular, a degree $1$ divisor can be generated by two places of successive degrees. \\

To compute the asymptotics of the elementary-abelian extensions of $F$, we follow David Wright's approach presented in \cite{Wr}. The concepts discussed in his Sections 2-4 up to page 32 also apply in our setting of wildly ramified elementary-abelian extensions in characteristic $p$. In the following, we work along Wright's approach and apply it to our situation. %In the following, we use Wright's notations as far as possible. However, we make some adjustments for our consideration of elementary abelian $p$-extensions. \\

\begin{definition}
Let $\mathbb{A}_F$ denote the adele ring of $F$ and $\mathbb{A}_F^{\times}$ the idele group of $F$. The set of homomorphisms from the idele class group $\mathbb{A}_F^{\times}/F^{\times}$ to the elementary-abelian group $C_p^f$ which we refer to as \emph{$C_p^f$-characters of the idele class group} is denoted by $\mathcal{X}_{C_p^f}(\mathbb{A}_F^{\times}/F^{\times}) \cong \prod\limits_{j = 1}^f \mathcal{X}_{C_p}(\mathbb{A}_F^{\times}/F^{\times})$. Furthermore, for any $\chi = (\chi_1, \ldots, \chi_f) \in \mathcal{X}_{C_p^f}(\mathbb{A}_F^{\times}/F^{\times})$ with $\chi_j \in \mathcal{X}_{C_p}(\mathbb{A}_F^{\times}/F^{\times})$ we define $U_{\chi} := \ker(\chi) = \bigcap\limits_{j=1}^f \ker(\chi_j)$. 

We have $\mathbb{A}_F^{\times} / U_{\chi} \cong C_p^{k(\chi)}$ for some $0 \leq k(\chi) \leq f$ and, by the idelic version of class field theory, this open subgroup $U_{\chi}$ corresponds to a Galois extension $F_{\chi}/F$ such that $\Gal(F_{\chi}/F) \cong C_p^{k(\chi)}$. The \emph{discriminant} $\partial(\chi)$ corresponding to the $C_p^f$-character $\chi$ is 
\begin{align*}
	\partial(\chi) := \partial\left(F_{\chi}/F\right).%^{\frac{\vert C_p^r \vert}{\vert C_p^{k(\chi)} \vert}} = \partial\left(F_{\chi}/F\right)^{p^{r-k(\chi)}}
\end{align*}
%which can be computed by the conductor-discriminant formula. 
For further details, see \cite[Section 2]{Wr}.
\end{definition}

\begin{definition}
For any $f \in \N_0$, we define
\begin{align*}
\mathfrak{F}_{C_p^{f}}(s) := \sum\limits_{\chi \in \mathcal{X}_{C_p^{f}}(\mathbb{A}_F^{\times}/F^{\times})} \vert \vert {\partial\left(F_{\chi}/F\right)} \vert \vert^{-p^{f-k(\chi)}s},
\end{align*}
which we call \emph{generating $C_p^f$-character series}. Note that the weight $p^{f-k(\chi)}$ is due to the conductor-discriminant formula.
\end{definition}

\begin{proposition} \label{WrightDecom}
The discriminant series $\Phi(F, C_p^r; s)$ has a decomposition as a linear combination of generating $C_p^f$-character series
\begin{align} \label{DiscDecom}
	\Phi(F, C_p^r; s) = \frac{1}{\vert {\Aut(C_p^r)} \vert} \quad \sum\limits_{f = 0}^r \binom{r}{f}_p \mu(C_p^{r-f}) \mathfrak{F}_{C_p^{f}}\left(p^{r-f}s\right) = \sum\limits_{f=0}^r \frac{e_f}{p^f} \cdot \mathfrak{F}_{C_p^{f}}\left(p^{r-f}s\right). 
\end{align}
\end{proposition}

\begin{proof}
The first equation is (2.14) in \cite[p. 29]{Wr}. The derivation of this equation can be pursued in detail in \cite[Section 2]{Wr}. The second equation is the definition of $e_f$ in Definition \ref{varphifunction}. 
Note that we have to evaluate $\mathfrak{F}_{C_p^{f}}$ at $p^{r-f}s$ by the conductor-discriminant formula since we consider the group $C_p^r$. 

For a transfer between Wright's and our notations, note that $\omega$ in \cite{Wr} is an idelic principal character of the form 
\begin{align*}
	\omega(x) = \vert x \vert_{A}^s = \prod\limits_{\idp \in \PP_F} \vert x \vert_{\idp}^s = \prod\limits_{\idp \in \PP_F} \np^{-\nu_{\idp}(x)s}.
\end{align*}
Additionally, we have $\phi(C_p^r) = \vert {\Aut(C_p^r)} \vert$ and we recall the Möbius $\mu$-function for finite abelian groups and the fact that the number of subgroups of type $C_p^f$ in $C_p^r$ for some $0 \leq f \leq r$ is $\binom{r}{f}_p$. 
\end{proof}

As our next step, we consider Section 3 in \cite{Wr}.

\begin{definition}
Let $\OO_S \subseteq F$ be the ring of $S$-integers in $F$ and $\OO_S^{\times}$ its group of units (the $S$-units). Furthermore, let $\mathcal{E}_{C_p^f}(S) := \prod\limits_{j=1}^f \OO_S^{\times}/ \left(\OO_S^{\times}\right)^p$ and let $\mathcal{A}_{C_p^f}(S)$ be a complete system of representatives for the group $\prod\limits_{j=1}^f \left(\left(\mathbb{A}_F^{\times}\right)^p F^{\times} \bigcap \mathbb{A}_F^{\times}(S) \right)/\left(\mathbb{A}_F^{\times}(S)\right)^p$ where $\mathbb{A}_F^{\times}(S) = \prod\limits_{\idp \in S} F_{\idp}^{\times} \times \prod\limits_{\idp \not \in S} \OO_{\idp}^{\times}$ is the group of $S$-ideles of $F$.  

Note that we have isomorphisms $F_{\idp}^{\times} \cong \Z \times \OO_{\idp}^{\times} \cong \Z \times \F_{q_{\idp}}^{\times} \times (1 + \pi_{\idp} \OO_{\idp})$ where $\F_{q_{\idp}}$ denotes the field of constants of $F_{\idp}^{\times}$ and $\pi_{\idp}$ is a uniformiser for the place $\idp$. 

Let $\mathcal{X}_{C_p^{f}}(\mathbb{A}_F^{\times}(S))$ denote the set of \emph{$C_p^f$-characters of the group of $S$-ideles}.
\end{definition}

Note that any $\chi \in \mathcal{X}_{C_p^{f}}(\mathbb{A}_F^{\times}(S))$ is uniquely described by a tuple of local maps $\chi_p: X \rightarrow C_p^f$ where $X = F_{\idp}^{\times}$ for $\idp \in S$ and $X = \OO_{\idp}^{\times}$ for $\idp \not \in S$ and $\chi_{\idp}$ is the trivial map for all but finitely many places $\idp \in \PP_F$. 

For any $\varepsilon = (\varepsilon_1, \ldots, \varepsilon_f) \in \mathcal{E}_{C_p^f}(S)$, the diagonal embedding of $\varepsilon$ to $\prod\limits_{j=1}^f \mathbb{A}_F^{\times}/\left(\mathbb{A}_F^{\times}\right)^p$, whose $\idp$-part is denoted by $\varepsilon_{\idp} = (\varepsilon_{1, \idp}, \ldots, \varepsilon_{f, \idp})$, is already contained in $\prod\limits_{j=1}^f \mathbb{A}_F^{\times}(S)/\mathbb{A}_F^{\times}(S)^p$. Any $\chi= (\chi_1, \ldots, \chi_f) \in \mathcal{X}_{C_p^{f}}(\mathbb{A}_F^{\times}(S))$ with $\chi_j \in \mathcal{X}_{C_p}(\mathbb{A}_F^{\times}(S))$ for $1 \leq j \leq f$ has an order dividing $p$. Hence, the evaluation of $\chi_j$ at $\varepsilon_j$, which we denote by $\chi_j(\varepsilon_j)$, is well-defined and we define
\begin{align*}
\chi(\varepsilon) := \prod\limits_{j=1}^f \chi_j(\varepsilon_j) = \prod\limits_{j=1}^f \prod\limits_{\idp \in \PP_F} \chi_{j, \idp}(\varepsilon_{j, \idp}) = \prod\limits_{\idp \in \PP_F} \prod\limits_{j=1}^f \chi_{j, \idp}(\varepsilon_{j, \idp}) = \prod\limits_{\idp \in \PP_F} \chi_{\idp}(\varepsilon_{\idp}) \in C_p.
\end{align*}
Note that we identify $C_p$ with the group of $p$-th roots of unity in $\C^{\times}$ which we denote by $\mu_p$.

%Note that we have isomorphisms $F_{\idp}^{\times} \cong \Z \times \OO_{\idp} \cong \Z \times \F_{q_{\idp}}^{\times} \times (1 + \pi_{\idp} \OO_{\idp})$ where $\F_{q_{\idp}}$ denotes the field of constants of $F_{\idp}^{\times}$ and $\pi_{\idp}$ is a uniformiser for the place $\idp$. \\
By \cite[Lemma 3.1]{Wr} and Wright's subsequent comment on this Lemma, we choose a sufficiently large finite set $S$ of places of $F$ so that $S$ generates $\Pic_F$. In this case, the $S$-idele class group is trivial and the Lemma yields a bijection between $\mathcal{E}_{C_p^f}(S)$ and $\mathcal{A}_{C_p^f}(S)$. 

\begin{proposition} \label{PropDecom}
For the generating $C_p^f$-character series, we have the identity
\begin{align} \label{CondDecom}
	\mathfrak{F}_{C_p^{f}}(s) = \frac{1}{\vert {\mathcal{E}_{C_p^f}(S)} \vert} \cdot \sum\limits_{\varepsilon \in \mathcal{E}_{C_p^f}(S)} \mathfrak{F}_{C_p^{f}, S}(s; \varepsilon)
\end{align}
where
\begin{align} \label{EulerProd}
	&\mathfrak{F}_{C_p^{f}, S}(s; \varepsilon) := \sum\limits_{\chi \in \mathcal{X}_{C_p^{f}}(\mathbb{A}_F^{\times}(S))} \chi(\varepsilon) \vert \vert {\partial\left(F_{\chi}/F\right)} \vert \vert^{-p^{f-k(\chi)}s} \notag \\
	= &\prod\limits_{\idp \in S} \sum\limits_{\chi_{\idp} \in \mathcal{X}_{C_p^{f}}(F_{\idp}^{\times})} \chi_{\idp}(\varepsilon_{\idp}) \vert \vert {\partial\left(F_{\chi_{\idp}}/F\right)} \vert \vert^{-p^{f-k(\chi_{\idp})}s} \times \prod\limits_{\idp \not \in S} \sum\limits_{\chi_{\idp} \in \mathcal{X}_{C_p^{f}}(\OO_{\idp}^{\times})} \chi_{\idp}(\varepsilon_{\idp}) \vert \vert {\partial\left(F_{\chi_{\idp}}/F\right)} \vert \vert^{-p^{f-k(\chi_{\idp})}s}.
\end{align}
%Thereby, $\chi_{\idp}$ describes the $\idp$-part of $\chi$ and $\varepsilon_{\idp}$ corresponds to the $\idp$-part of the diagonal embedding of $\varepsilon$ to $\prod\limits_{i=1}^r \mathbb{A}_F^{\times}(S)$. In particular, we have $\chi(\varepsilon) = \prod\limits_{\idp \in \PP_F} \chi_{\idp}(\varepsilon_{\idp})$. 
\end{proposition}

\begin{proof}
Equation \eqref{CondDecom} is equation (3.2) in \cite[p. 31]{Wr} and equation \eqref{EulerProd} is equation (3.4) in \cite[p. 31]{Wr}. The derivation of these equations can be pursued in detail in \cite[Section 3]{Wr}. 
\end{proof}

In the following part, we have to study these Euler products for different classes $\varepsilon \in \mathcal{E}_{C_p^i}(S)$ to be able to deduce their asymptotics. But first we should take a closer look at the structure of the group $\mathcal{E}_{C_p^i}(S)$. We recall the so called S-unit theorem for this purpose.

\begin{proposition} \label{S-unitTheorem}
The structure of the $S$-unit group $\OO_S^{\times}$ can be characterised by an isomorphism of groups: %the following isomorphism
\begin{align*}
\OO_S^{\times} \cong \F_q^{\times} \times \Z^{\vert S \vert -1}.
\end{align*}
\end{proposition}

\begin{proof}
For the proof, see \cite[Proposition 14.2]{St}.
\end{proof}

\begin{corollary} \label{EpsilonCharacters}
The structure of the group $\mathcal{E}_{C_p^f}(S)$ can be characterised by an isomorphism of groups: %the isomorphism
\begin{align*} 
\mathcal{E}_{C_p^f}(S) \cong C_p^{f(\vert S \vert -1)}.
\end{align*}
\end{corollary}

\begin{proof}
By the previous Proposition \ref{S-unitTheorem}, we have 
\begin{align*}
\mathcal{E}_{C_p^f}(S) = \prod\limits_{j=1}^f \OO_S^{\times}/ \left(\OO_S^{\times}\right)^p \cong \prod\limits_{j=1}^f C_p^{(\vert S \vert -1)} \cong C_p^{f(\vert S \vert -1)}.
\end{align*}
\end{proof}

Initially, we consider the trivial class $\mathds{1}_f := \prod\limits_{i=1}^{f} \left( \OO_S^{\times} \right)^p$ in $\mathcal{E}_{C_p^f}(S)$ and compute the Euler factors in equation \eqref{EulerProd}. %Note that we have isomorphisms $F_{\idp}^{\times} \cong \Z \times \OO_{\idp} \cong \Z \times \F_{q_{\idp}}^{\times} \times (1 + \pi_{\idp} \OO_{\idp})$ where $\F_{q_{\idp}}$ denotes the field of constants of $F_{\idp}^{\times}$ and $\pi_{\idp}$ is a uniformiser for the place $\idp$. \\ % with $\chi(\varepsilon) = 1$ for all $\chi \in \mathcal{X}_{C_p^{f}}(\mathbb{A}_F^{\times}(S))$. Due to $F_{\idp}^{\times} \cong \Z \times \OO_{\idp}^{\times}$, we have

\begin{proposition} \label{TrivialCharacter}
Let $0 \leq f \leq r$ and $\mathds{1}_f = (1, \ldots, 1) \in \mathcal{E}_{C_p^f}(S)$ be the trivial class. Then, we have $\frac{1}{\vert {\mathcal{E}_{C_p^f}(S)} \vert} \cdot \mathfrak{F}_{C_p^{f}, S}\left(p^{r-f}s; \varepsilon\right) = p^f \cdot \Phi_{f}(F, C_p^r; s)$.
\end{proposition}

\begin{proof}
Since $\mathds{1}_f$ is the class of $p$-th powers in $\mathcal{E}_{C_p^f}(S)$, we have $\chi(\varepsilon) = 1$ for all characters $\chi \in \mathcal{X}_{C_p^{f}}(\mathbb{A}_F^{\times}(S))$. 

Let $\idp \in S$. By $F_{\idp}^{\times} \cong \Z \times \OO_{\idp}^{\times}$ and the fact that the $\Z$-component encodes the unramified part, we have
\begin{align} \label{eq700}
\sum\limits_{\chi_{\idp} \in \mathcal{X}_{C_p^{f}}(F_{\idp}^{\times})} \vert \vert {\partial\left(F_{\chi_{\idp}}/F\right)} \vert \vert^{-p^{r-k(\chi_{\idp})}s} = p^f \cdot \sum\limits_{\chi_{\idp} \in \mathcal{X}_{C_p^{f}}(\OO_{\idp}^{\times})} \vert \vert {\partial\left(F_{\chi_{\idp}}/F\right)} \vert \vert^{-p^{r-k(\chi_{\idp})}s}
\end{align}
for the Euler factor of $\mathfrak{F}_{C_p^{f}, S}\left(p^{r-f}s; \varepsilon\right)$ corresponding to $\idp$. Now, the equations \eqref{EulerProd} and \eqref{eq700} and Corollary \ref{EpsilonCharacters} yield
\begin{align*}
\frac{1}{\vert {\mathcal{E}_{C_p^f}(S)} \vert} \cdot \mathfrak{F}_{C_p^{f}, S}\left(p^{r-f}s; \varepsilon\right) = p^f \cdot \prod\limits_{\idp \in \PP_F} \sum\limits_{\chi_{\idp} \in \mathcal{X}_{C_p^{f}}(\OO_{\idp}^{\times})} \vert \vert {\partial\left(F_{\chi_{\idp}}/F\right)} \vert \vert^{-p^{r-k(\chi_{\idp})}s}.
\end{align*}

It remains to prove $\sum\limits_{\chi_{\idp} \in \mathcal{X}_{C_p^{f}}(\OO_{\idp}^{\times})} \vert \vert {\partial\left(F_{\chi_{\idp}}/F\right)} \vert \vert^{-p^{r-k(\chi_{\idp})}s} = \Phi_{f, \idp}(F, C_p^r; s)$ by Definition \ref{DefinitionDirichletSeries}. The equality follows since Proposition \ref{QuotientSpace}, Proposition \ref{RS} and Corollary \ref{VSIso} provide an explicit description of the class field theory for local function fields. In fact, we have both a bijective map between $\F_p$-subspaces of $F_{\idp}/\wp(F_{\idp})$ and Artin-Schreier extensions of $F_{\idp}$ by Proposition \ref{QuotientSpace} and between open subgroups of $F_{\idp}^{\times}$ containing $\left(F_{\idp}^{\times}\right)^p$ and Artin-Schreier extensions of $F_{\idp}$ by local class field theory. Additionally, all unramified Artin-Schreier extensions of $F_{\idp}$ correspond to an open subgroup $U$ with $p\Z \times \OO_{\idp}^{\times} \subseteq U \subseteq F_{\idp}^{\times}$. Hence, the characters $\chi_{\idp} \in \mathcal{X}_{C_p^{f}}(\OO_{\idp}^{\times})$ describe the ramification in the place $\idp$ in the same way as the subspace $F_{\idp}^{ram}/\wp(F_{\idp}^{ram})$ in Corollary \ref{VSIso}.
\end{proof}

\begin{example}
In the case $F = \Fqtg$, we have $\Pic_{\Fqtg} \cong \Z$ so that the Picard group is generated by the infinite place $\infty$ of degree $1$ and we have $\mathcal{E}_{C_p^r}(\infty) \cong C_1$ by Corollary \ref{EpsilonCharacters}. Hence, we recover
\begin{align*}
\Phi(\Fqtg, C_p^r; s) = \frac{1}{\vert {\Aut(C_p^r)} \vert} \quad \sum\limits_{f = 0}^r \vert C_p^f \vert \cdot \binom{r}{f}_p \mu\left(C_p^{r-f}\right) \cdot \Phi_{f}(\Fqtg, C_p^r; s)
\end{align*}
by Proposition \ref{WrightDecom}, Proposition \ref{PropDecom} and Proposition \ref{TrivialCharacter}. 

One might ask what happens if we include additional places to our set $S$ since this procedure increases $\mathcal{E}_{C_p^f}(S)$ so that we also have non-trivial classes. The computations for such an extended $S$, however, certainly have to reveal the same result. In fact, we actually prove in Proposition \ref{NoSelmer} that the Euler products for all non-trivial classes have to vanish in the case of $\Fqtg$.
\end{example}

In the following, we explore the cause of this phenomenon in more general terms by computing character sums in the Euler factors in equation \eqref{EulerProd}. First, we define some notation for those character sums.

\begin{definition}
Let $1 \leq f$, $1 \leq m$ and $X \in \{\OO_{\idp}^{\times}, F_{\idp}^{\times}\}$. Then, we define
\begin{align*}
	\Sigma(\mathcal{X}_{C_p^f}(X), \leq m; \varepsilon_{\idp}) := \sum\limits_{\substack{\chi_{\idp} \in \mathcal{X}_{C_p^f}(X): \\ \ker(\chi_{\idp}) \supseteq H}} \chi_{\idp}(\varepsilon_{\idp}),
\end{align*}
where $H = p\Z \times \F_{q_{\idp}}^{\times} \times (1 + \pi^m_{\idp} \OO_{\idp})$ for $X = F_{\idp}^{\times}$ and $H = \F_{q_{\idp}}^{\times} \times (1 + \pi^m_{\idp} \OO_{\idp})$ for $X = \OO_{\idp}^{\times}$, to be the sum over all $C_p^f$-characters of $X$ that have a conductor which divides $\mathfrak{p}^{m}$ evaluated at $\varepsilon_{\idp}$.

Furthermore, we define $\Sigma(\mathcal{X}_{C_p^f}(X), m; \varepsilon_{\idp})$ to be the sum over all $C_p^f$-characters of $X$ whose conductor is precisely $\mathfrak{p}^{m}$ for $m \geq 0$ evaluated at $\varepsilon_{\idp}$. As there are no $C_p^f$-characters with conductor $\mathfrak{p}$, which would imply tame ramification, we note $\Sigma(\mathcal{X}_{C_p^f}(X), 1; \varepsilon_{\idp}) = 0$.
For $m \geq 1$, we have 
\begin{align*}
	\Sigma(\mathcal{X}_{C_p^f}(X), m; \varepsilon_{\idp}) = \Sigma(\mathcal{X}_{C_p^f}(X), \leq m; \varepsilon_{\idp}) - \Sigma(\mathcal{X}_{C_p^f}(X), \leq m-1; \varepsilon_{\idp}).
\end{align*}

Analogously, we define character sums for chains. For any chain $\mathfrak{C} := \left(c_1, \ldots, c_h\right) \in \N^{h}$ with $0 \leq h \leq f$, we define $\Sigma(\mathcal{X}_{C_p^f}(X), \mathfrak{C}; \varepsilon_{\idp})$ to be the sum over all characters $\chi_{\idp}$ such that the chain of conductor exponents (Proposition \ref{ConductorChain}) of the corresponding extension $F_{\chi_{\idp}}/F_{\idp}$ is $\mathfrak{C}$. Again, we evaluate the characters at $\varepsilon_{\idp}$.  For a character $\chi_{\idp}$, we also denote the chain corresponding to $F_{\chi_{\idp}}/F_{\idp}$ by $\mathfrak{C}(\chi_{\idp})$. We recall that the discriminant of this chain can be computed by Lemma \ref{DiscExponent}. In particular, for $1 \leq h$, all characters $\chi_{\idp} \in \mathcal{X}_{C_p^f}(X)$ such that $\mathfrak{C}(\chi_{\idp}) = \mathfrak{C}$ have conductor $\idp^{c_1}$ and we have the identity
\begin{align*}
\Sigma(\mathcal{X}_{C_p^f}(X), m; \varepsilon_{\idp}) = \sum\limits_{\substack{\mathfrak{C}' \in \N^{h}: \\ 1 \leq h \leq f, \\m = c_1'.}} \Sigma(\mathcal{X}_{C_p^f}(X), \mathfrak{C}'; \varepsilon_{\idp}).
\end{align*}
For $h=0$, we have $\Sigma(\mathcal{X}_{C_p^f}(X), 0; \varepsilon_{\idp}) = \Sigma(\mathcal{X}_{C_p^f}(X), \{\}; \varepsilon_{\idp})$ for the empty chain.
\end{definition}

\begin{proposition} \label{NoSelmer}
Let $1 \leq f$ and $\varepsilon = (\varepsilon_1, \ldots, \varepsilon_f) \in \mathcal{E}_{C_p^f}(S)$ be a class such that there is at least one element $\varepsilon_j$ for some $1 \leq j \leq f$ that does not represent a class in the $p$-Selmer group. Then, we have $\mathfrak{F}_{C_p^{f}, S}(s; \varepsilon) = 0$.
\end{proposition}

\begin{proof}
If there is such an element $\varepsilon_j$ that does not represent a class in the $p$-Selmer group, then the principal divisor $(\varepsilon_j)$ is not a $p$-th power of any divisor. Equivalently, there is at least one place $\mathfrak{q} \in S$ such that $p \nmid \nu_{\mathfrak{q}}(\varepsilon_j)$.
We consider equation \eqref{EulerProd} and compute the Euler factor belonging to the place $\mathfrak{q} \in S$. Our aim is to show that this Euler factor vanishes. 

To this end, we look at the character sums $\Sigma(\mathcal{X}_{C_p^f}(F_{\mathfrak{q}}^{\times}), \leq m; \varepsilon_{\mathfrak{q}})$ for $m \geq 1$. The set of characters whose conductor divide $\mathfrak{q}^m$ forms a group. Hence, we can apply the orthogonality of characters which means that we have $\Sigma(\mathcal{X}_{C_p^f}(F_{\mathfrak{q}}^{\times}), \leq m; \varepsilon_{\mathfrak{q}}) = 0$ if there is at least one character $\chi_{\mathfrak{q}}$ such that $\ker(\chi_{\mathfrak{q}}) \supseteq p\Z \times \F_{q_{\mathfrak{q}}}^{\times} \times (1 + \pi^m_{\mathfrak{q}} \OO_{\mathfrak{q}})$ and $\chi_{\mathfrak{q}}(\varepsilon_{\mathfrak{q}}) \neq 1$. In particular, if $m_{\mathfrak{q}}$ denotes the minimal integer with these properties, then we have $\Sigma(\mathcal{X}_{C_p^f}(F_{\mathfrak{q}}^{\times}), \leq m_{\mathfrak{q}}; \varepsilon_{\mathfrak{q}}) = 0$ and $\Sigma(\mathcal{X}_{C_p^f}(F_{\mathfrak{q}}^{\times}), m; \varepsilon_{\mathfrak{q}}) = 0$ for all $m > m_{\mathfrak{q}}$. 

In the case $f = j = 1$, it is sufficient to verify $m_{\mathfrak{q}} = 1$. For this purpose, we directly compute $\Sigma(\mathcal{X}_{C_p}(F_{\mathfrak{q}}^{\times}), \leq 1; \varepsilon_{\mathfrak{q}})$. Let $\chi_{\mathfrak{q}}: F_{\mathfrak{q}}^{\times} \rightarrow C_p$ be a character satisfying $\ker(\chi_{\mathfrak{q}}) \supseteq p\Z \times \F_{q_{\mathfrak{q}}}^{\times} \times (1 + \pi_{\mathfrak{q}} \OO_{\mathfrak{q}}) \cong p\Z \times \OO_{\mathfrak{q}}^{\times}$ and $x \in F_{\mathfrak{q}}^{\times}$ be arbitrary. Due to the isomorphism $\Z \times \OO_{\mathfrak{q}}^{\times} \cong F_{\mathfrak{q}}^{\times}, (z, \alpha) \mapsto \pi_{\mathfrak{q}}^z \cdot \alpha$, we have $\chi_{\mathfrak{q}}(x) = \zeta_p^{\nu_{\mathfrak{q}}(x)}$ for some $p$-th root of unity $\zeta_p$. If we choose $\zeta_p$ to be primitive and $x = \varepsilon_{\mathfrak{q}}$, then we conclude $\chi_{\mathfrak{q}}(\varepsilon_{\mathfrak{q}}) \neq 1$ since $p \nmid \nu_{\mathfrak{q}}(\varepsilon_{\mathfrak{q}})$ by our assumption. 

In the case $f > 1$, we can certainly also conclude $\Sigma(\mathcal{X}_{C_p^f}(F_{\mathfrak{q}}^{\times}), \leq 1; \varepsilon_{\mathfrak{q}}) = 0$ and $\Sigma(\mathcal{X}_{C_p^f}(F_{\mathfrak{q}}^{\times}), m; \varepsilon_{\mathfrak{q}}) = 0$ for all $m >1$ by our previous considerations since the $C_p^f$-characters are tuples of $C_p$-characters. We, however, have to keep in mind that we should prove $\Sigma(\mathcal{X}_{C_p^f}(F_{\mathfrak{q}}^{\times}), \mathfrak{C}; \varepsilon_{\mathfrak{q}}) = 0$ for all chains $\mathfrak{C} \in \N^{h}$ for some $0 \leq h \leq f$ to obtain the vanishing of the Euler factor. First, we note that there is only one chain for $h=0$, namely the empty chain whose character sum is $\Sigma(\mathcal{X}_{C_p^f}(F_{\mathfrak{q}}^{\times}), \leq 1; \varepsilon_{\mathfrak{q}}) = 0$.
Therefore, we can assume $1 \leq h$ without loss of generality. Now, we choose some character $\sigma_{\mathfrak{q}} = (\sigma_{\mathfrak{q}, 1}, \ldots, \sigma_{\mathfrak{q}, f}) \in \mathcal{X}_{C_p^f}(F_{\mathfrak{q}}^{\times})$ such that $\ker(\sigma_{\mathfrak{q}}) \supseteq p\Z \times \F_{q_{\mathfrak{q}}}^{\times} \times (1 + \pi_{\mathfrak{q}} \OO_{\mathfrak{q}})$ and $\sigma_{\mathfrak{q}}(\varepsilon_{\mathfrak{q}}) \neq 1$. For instance, we could choose $\sigma_{\mathfrak{q}} = (\mathds{1}, \ldots, \sigma_{\mathfrak{q}, j}, \ldots, \mathds{1})$ where $\sigma_{\mathfrak{q}, j}$ is a $C_p$-character satisfying $\sigma_{\mathfrak{q}, j}(\varepsilon_{j, \mathfrak{q}}) \neq 1$ and all other components contain the trivial character. The existence of such a character is guaranteed by our considerations in the case $f = 1$. 

Now, we fix a chain $\mathfrak{C} \in \N^h$ with $1 \leq h \leq f$ and we compute the character sum $\Sigma(\mathcal{X}_{C_p^f}(F_{\mathfrak{q}}^{\times}), \mathfrak{C}; \varepsilon_{\mathfrak{q}})$. For any character $\chi_{\mathfrak{q}} \in \mathcal{X}_{C_p^f}(F_{\mathfrak{q}}^{\times})$ with $\mathfrak{C}(\chi_{\mathfrak{q}}) = \mathfrak{C}$, we note that the character $\chi_{\mathfrak{q}} \cdot \sigma_{\mathfrak{q}} \in \mathcal{X}_{C_p^f}(F_{\mathfrak{q}}^{\times})$ satisfies $\mathfrak{C}(\chi_{\mathfrak{q}} \cdot \sigma_{\mathfrak{q}}) = \mathfrak{C}$ since the character $\sigma_{\mathfrak{q}}$ has conductor $\mathfrak{1}$ and belongs to the empty chain. In fact, the kernels of $\chi_{\mathfrak{q}}$ and $\chi_{\mathfrak{q}} \cdot \sigma_{\mathfrak{q}}$ coincide on $p\Z\ \times F_q^{\times} \times (1 + \pi_{\mathfrak{q}} \OO_{\mathfrak{q}})$ as $\sigma_{\mathfrak{q}}$ is trivial on this subgroup. Therefore, we obtain
\begin{align} \label{ChiTrick}
	\Sigma(\mathcal{X}_{C_p^f}(F_{\mathfrak{q}}^{\times}), \mathfrak{C}; \varepsilon_{\mathfrak{q}})
	&= \sum\limits_{\substack{\chi_{\mathfrak{q}} \in \mathcal{X}_{C_p^f}(F_{\mathfrak{q}}^{\times}): \\ \mathfrak{C}(\chi_{\mathfrak{q}})  
			= \mathfrak{C}}} \chi_{\mathfrak{q}}(\varepsilon_{\mathfrak{q}})
	=  \sum\limits_{\substack{\chi'_{\mathfrak{q}} \in \mathcal{X}_{C_p^f}(F_{\mathfrak{q}}^{\times}): \\ \mathfrak{C}(\chi'_{\mathfrak{q}}) = \mathfrak{C}}} \left(\chi'_{\mathfrak{q}} \cdot \sigma_{\mathfrak{q}}\right)(\varepsilon_{\mathfrak{q}}) \notag \\
	&= \sigma_{\mathfrak{q}}(\varepsilon_{\mathfrak{q}}) \sum\limits_{\substack{\chi'_{\mathfrak{q}} \in \mathcal{X}_{C_p^i}(F_{\mathfrak{q}}^{\times}): \\ \mathfrak{C}(\chi'_{\mathfrak{q}}) = \mathfrak{C}}} \chi'_{\mathfrak{q}}(\varepsilon_{\mathfrak{q}})
\end{align}
where we make use of the substitution $\chi_{\mathfrak{q}} = \chi'_{\mathfrak{q}} \cdot \sigma_{\mathfrak{q}}$. Now, equation \eqref{ChiTrick} yields 
\begin{align*}
0 = (1 - \sigma_{\mathfrak{q}}(\varepsilon_{\mathfrak{q}})) \cdot \sum\limits_{\substack{\chi_{\mathfrak{q}} \in \mathcal{X}_{C_p^f}(F_{\mathfrak{q}}^{\times}): \\ \mathfrak{C}(\chi_{\mathfrak{q}})  
		= \mathfrak{C}}} \chi_{\mathfrak{q}}(\varepsilon_{\mathfrak{q}}). 
\end{align*}
Due to $\sigma_{\mathfrak{q}}(\varepsilon_{\mathfrak{q}}) \neq 1$, we conclude $\sum\limits_{\substack{\chi_{\mathfrak{q}} \in \mathcal{X}_{C_p^f}(F_{\mathfrak{q}}^{\times}): \\ \mathfrak{C}(\chi_{\mathfrak{q}})  
		= \mathfrak{C}}} \chi_{\mathfrak{q}}(\varepsilon_{\mathfrak{q}}) =0$. This completes the proof for $f > 1$.
\end{proof}

\begin{remark}
By Proposition \ref{WrightDecom}, Proposition \ref{PropDecom}, Proposition \ref{TrivialCharacter} and Proposition \ref{NoSelmer}, we obtain
\begin{align*}
	\Phi(F, C_p^r;s) = \sum\limits_{f=0}^{r} e_f \cdot \Phi_f(F, C_p^r;s)
\end{align*}
for all global function fields $F$ whose $p$-Selmer group, respectively $p$-torsion of the class group, is trivial. Note that this is the shape of the decomposition of the Dirichlet series of $\Fqtg$ that we proved in Proposition \ref{LCEuler}.
\end{remark}

%\begin{remark}
%Applying Proposition \ref{NoSelmer} in the case of a rational function field $\Fqtg$ yields that the Euler products for all non-trivial classes $\varepsilon \in \mathcal{E}_{C_p^f}(S)$ vanish since the $p$-Selmer group is trivial for the base field $\Fqtg$. Consequently, we recover Proposition \ref{LCEuler} for all finite sets $S$ of places that contain the infinite place $\infty$.
%\end{remark}

It remains to investigate the Euler products in the case where we consider a non-trivial class $\varepsilon = (\varepsilon_1, \ldots, \varepsilon_f) \in \mathcal{E}_{C_p^f}(S)$ such that all components of $\varepsilon$ represent a class in the $p$-Selmer group.

\begin{lemma} \label{EulerfactorO}
Let $1 \leq f \leq r$ and $\varepsilon$ a non-trivial class $\varepsilon = (\varepsilon_1, \ldots, \varepsilon_f) \in \mathcal{E}_{C_p^f}(S)$ such that all components of $\varepsilon$ represent a class in the $p$-Selmer group. Let further $\idp \in S$. Then, we have
\begin{align*}
	\sum\limits_{\chi_{\idp} \in \mathcal{X}_{C_p^{f}}(F_{\idp}^{\times})} \chi_{\idp}(\varepsilon_{\idp}) \vert \vert {\partial\left(F_{\chi_{\idp}}/F\right)} \vert \vert^{-p^{r-k(\chi_{\idp})}s} = p^f \cdot \sum\limits_{\chi_{\idp} \in \mathcal{X}_{C_p^{f}}(\OO_{\idp}^{\times})} \chi_{\idp}(\varepsilon_{\idp}) \vert \vert {\partial\left(F_{\chi_{\idp}}/F\right)} \vert \vert^{-p^{r-k(\chi_{\idp})}s}.
\end{align*}
for the Euler factor of $\mathfrak{F}_{C_p^{f}, S}(p^{r-f} s; \varepsilon)$ corresponding to $\idp$ in Proposition \ref{PropDecom}.
\end{lemma}

\begin{proof}
	We have $p \mid \nu_{\idp}(\varepsilon_{j, \idp})$ for $1 \leq j \leq f$ by our assumption on $\varepsilon$, so $\chi_{\idp}(\varepsilon_{\idp})$ is trivial on the $\Z$-component of $F_{\idp}^{\times}$ for all characters $\chi_{\idp} \in \mathcal{X}_{C_p^f}(F_{\idp}^{\times})$. This yields the claim.
\end{proof}

By Lemma \ref{EulerfactorO}, we can restrict our considerations to characters of $\OO_{\idp}^{\times}$ for all places $\idp \in \PP_F$ regardless of affiliation to the finite set $S$ of places.

\begin{proposition} \label{ChiCond}
Let $1 \leq f$ and $\varepsilon$ a non-trivial class $\varepsilon = (\varepsilon_1, \ldots, \varepsilon_f) \in \mathcal{E}_{C_p^f}(S)$ such that all components of $\varepsilon$ represent a class in the $p$-Selmer group. Then, there exists a minimal non-negative integer $m_{\idp}$ for each place $\idp \in \PP_F$ depending on $\varepsilon$ such that $\Sigma(\mathcal{X}_{C_p^f}(\OO_{\idp}^{\times}), \leq m_{\idp}; \varepsilon_{\idp}) = 0$ and $\Sigma(\mathcal{X}_{C_p^f}(\OO_{\idp}^{\times}), m; \varepsilon_{\idp}) = 0$ for all $m > m_{\idp}$. 

Additionally, we have $m_{\idp} = 2$ for all but finitely many places $\idp \in \PP_F$ and the finitely many places with the property $m_{\idp} > 2$ satisfy $\deg(\idp) \leq 2g_F - 2$.
\end{proposition}

\begin{proof}
Again, we first consider the special case $f = 1$ and let $\idp \in \PP_F$ be an arbitrary place. In this case, we have to show that there is an integer $m_{\idp}$ and a character $\sigma_{\idp}: \OO_{\idp}^{\times} \rightarrow C_p$ with conductor $\idp^{m_{\idp}}$ such that $\sigma_{\idp}(\varepsilon_{\idp}) \neq 1$ again due to the orthogonality of characters. 

By our assumption on $\varepsilon$, we have $p \mid \nu_{\idp}(\varepsilon)$ and hence $\varepsilon  = u_{\idp} \cdot \pi_{\idp}^{j_{\idp}p}$ for some unit $u_{\idp} \in \OO_{\idp}^{\times}$ and some $j_{\idp} \in \Z$. 

For any integer $m \geq 2$, we consider the set of $C_p$-characters of $\OO_{\idp}^{\times}$ whose conductor divides $\idp^m$. Such a character $\chi_{\idp}$ has to satisfy $\ker(\chi_{\idp}) \supseteq \F_{q_{\idp}}^{\times} \times (1 + \pi_{\idp}^m \OO_{\idp})$. Note that we have $\chi_{\idp}(\varepsilon_{\idp}) = \chi_{\idp}(u_{\idp}) = 1$ for all $C_p$-characters of $\OO_{\idp}^{\times}$ whose conductor divides $\idp^m$ if and only if $u_{\idp} \in  \left(\F_{q_{\idp}}^{\times} \times (1 + \pi_{\idp}^m \OO_{\idp})\right) \cdot \left(\OO_{\idp}^{\times}\right)^p$. Hence, we have to prove that there is a minimal $m_{\idp}$ such that $u_{\idp} \not \in \left(\F_{q_{\idp}}^{\times} \times (1 + \pi_{\idp}^{m_{\idp}} \OO_{\idp})\right) \cdot \left(\OO_{\idp}^{\times}\right)^p$. 

Now, let $u_{\idp} = \sum\limits_{i \geq 0} a_i \pi_{\idp}^i$ be the Laurent series expansion of $u_{\idp}$. Since $u_{\idp}$ is not a $p$-th power by our assumption on $\varepsilon$, there is a minimal index $n_{\idp}$ with $p \nmid n_{\idp}$ and $a_{n_{\idp}} \neq 0$. Note that we can erase all terms for $i < n_{\idp}$ by multiplying $u_{\idp}$ with a proper element of $\left(\OO_{\idp}^{\times}\right)^p$. Due to $p \nmid n_{\idp}$, we, however, can not erase the term for this index. Consequently, we set $m_{\idp} := n_{\idp}+1$ and we have $u_{\idp} \not \in \left(\F_{q_{\idp}}^{\times} \times (1 + \pi_{\idp}^{m_{\idp}} \OO_{\idp})\right) \cdot \left(\OO_{\idp}^{\times}\right)^p$. This yields the existence of $m_{\idp}$ and the statements about the character sums. 

For the more precise statement about the $m_{\idp}$, we consider the differential $\frac{d u_{\idp}}{d \pi_{\idp}} = \sum\limits_{\substack{i \geq 1 \\ p \nmid i}} i a_i \pi_{\idp}^{i-1}$ and take a close look at \cite[Lemma 4.4]{La2} and its proof. In the proof, it is shown that the logarithmic differential $\omega := \frac{d \varepsilon}{\varepsilon}$ is regular and its divisor is effective. Furthermore, the equation
\begin{align} \label{Lag44}
	2g_F - 2 = \deg(\omega) = \sum\limits_{\idp \in \PP_F} \nu_{\idp}(\omega) \deg(\idp) = \sum\limits_{\idp \in \PP_F} \nu_{\idp}\left(\frac{du_{\idp}}{d \pi_{\idp}} \right) \deg(\idp)
\end{align}
is shown. By our previous computation, we have $\nu_{\idp}\left(\frac{du_{\idp}}{d \pi_{\idp}} \right) = n_{\idp}-1$. Thus, we deduce $n_{\idp} = 1$ for all places $\idp \in \PP_F$ satisfying $\deg(\idp) > 2g_F -2$ by equation \eqref{Lag44}. Moreover, we deduce $n_{\idp} < 2g_F$ and, hence, $m_{\idp} \leq 2g_F$ for all places $\idp \in \PP_F$ satisfying $\deg(\idp) \leq 2g_F -2$. This concludes the proof for $f = 1$. 

For $f > 1$, we have to take into account all components $\varepsilon_i$ of $\varepsilon$ for $1 \leq i \leq f$. Clearly, the considerations we made in the case $f = 1$ directly apply to each $\varepsilon_i$. In particular, we obtain non-negative integers $m_{i,\idp}$ for each $1 \leq i \leq f$ such that $u_{i, \idp} \not \in \left(\F_{q_{\idp}}^{\times} \times (1 + \pi_{\idp}^{m_{i, \idp}} \OO_{\idp})\right) \cdot \left(\OO_{\idp}^{\times}\right)^p$ where $\varepsilon_i = u_{i, \idp} \cdot \pi^{j_{i, \idp} p}$ for some unit $u_{i, \idp} \in \OO_{\idp}^{\times}$ and some $j_{i, \idp} \in \Z$. 

Now, we consider the minimum of the $m_{i, \idp}$ for $1 \leq i \leq f$ and claim that this integer is $m_{\idp}$. For this integer, we have at least one index $1 \leq j \leq f$ and one $C_p$-character $\sigma_{1, \idp}$ with conductor $\idp^{m_{\idp}}$ such that $\sigma_{1, \idp}(\varepsilon_{j, \idp}) \neq 1$ by our results in the case $f=1$. In particular, we can construct a $C_p^f$-character $\sigma_{f, \idp} = (\mathds{1}, \ldots, \sigma_{1, \idp}, \ldots, \mathds{1})$ with $\sigma_{1, \idp}$ in the $j$-th component which has conductor $\idp^{m_{\idp}}$ and satisfies $\sigma_{f, \idp}(\varepsilon_{\idp}) \neq 1$. By the orthogonality of characters, the statements about the character sums follow. The remaining statements directly follow by the characterisation of $m_{\idp}$ as the minimum of all $m_{i, \idp}$ for $1 \leq i \leq f$ and the results for the case $f=1$.
\end{proof}

Note that the computation of the character sums for the conductors is already sufficient to obtain the structure of the Euler factors in the special case $f=1$.

\begin{corollary} \label{LagResult}
	Let $\idp \in \PP_F$ be a fixed place satisfying $\deg(\idp) > 2g_F -2$. Furthermore, let $1 = f \leq r$ and let $\varepsilon \in \mathcal{E}_{C_p}(S)$ a non-trivial class such that $\varepsilon$ represents a class in the $p$-Selmer group. Then, we have	
	\begin{align*}
		\sum\limits_{\chi_{\idp} \in \mathcal{X}_{C_p}(\OO_{\idp}^{\times})} \chi_{\idp}(\varepsilon_{\idp}) \vert \vert {\partial\left(F_{\chi_{\idp}}/F\right)} \vert \vert^{-p^{r-k(\chi_{\idp})}s} = 1 - \np^{-(p-1)p^{r-1}2s}
	\end{align*}
	for the Euler factor of $\mathfrak{F}_{C_p, S}(p^{r-1}s; \varepsilon)$ corresponding to the place $\idp$ in Proposition \ref{PropDecom}.
\end{corollary}

\begin{proof}
	This is a direct consequence by Proposition \ref{ChiCond}. Furthermore, we note that this result matches with the computation carried out in the proof of \cite[Lemma 6.1]{La2} when we apply the conductor-discriminant formula.
\end{proof}

\begin{remark} \label{Rk1}
	If we consider the situation from Corollary \ref{LagResult} and a place $\idp$ with $\deg(\idp) \leq 2g_F -2$, we still have that the corresponding Euler factor is a finite polynomial in $\np^{-s}$. In fact, this polynomial has degree $(p-1) p^{r-1} m_{\idp}$.
\end{remark}

Next, we have to compute the character sums for different chains of conductor exponents to obtain the structure of the Euler factors for arbitrary $f$.

\begin{proposition} \label{ChiVanish}
Let $1 \leq f$, $\idp \in \PP_F$ a fixed place, $m_{\idp}$ the minimal non-negative integer satisfying the properties in Proposition \ref{ChiCond} and $\varepsilon$ a non-trivial class $\varepsilon = (\varepsilon_1, \ldots, \varepsilon_f) \in \mathcal{E}_{C_p^f}(S)$ such that all components of $\varepsilon$ represent a class in the $p$-Selmer group.
\begin{enumerate}
	\item[a)] Let $\mathfrak{C} = (c_1, \ldots, c_f) \in \N^{f}$ be an arbitrary chain satisfying $c_f > m_{\idp}$.
	Then, we have $\Sigma(\mathcal{X}_{C_p^f}(\OO_{\mathfrak{p}}^{\times}), \mathfrak{C}; \varepsilon_{\mathfrak{p}}) = 0$. 
	\item[b)] Let $\mathfrak{C} \in \N^{h}$ for some $1 \leq h \leq f$ be an arbitrary chain and $c_f \leq m_{\idp}$ if $h=f$. Then, we have 
	\begin{align} \label{ZeroSum}
		\sum\limits_{\substack{\mathfrak{C}' \in \N^j: \\ 1 \leq j \leq f, \\ \mathfrak{C}'_{> m_{\idp}} = \mathfrak{C}_{> m_{\idp}}.}} \Sigma(\mathcal{X}_{C_p^f}(\OO_{\mathfrak{p}}^{\times}), \mathfrak{C}'; \varepsilon_{\mathfrak{p}}) = 0
	\end{align}
	where $\mathfrak{C}_{> m_{\idp}}$ denotes the subchain of $\mathfrak{C}$ consisting of all components of $\mathfrak{C}$ that are larger than $m_{\idp}$.
\end{enumerate}	
\end{proposition}

\begin{proof}
Our approach is quite analogous to that procedure we applied in the proof of Proposition \ref{NoSelmer}. 

First, we construct a $C_p^f$-character $\sigma_{\idp}$ which has conductor $\idp^{m_{\idp}}$ and satisfies $\sigma_{\idp}(\varepsilon_{\idp}) \neq 1$. The construction is completely analogous to that we carried out in the proof of the previous Proposition \ref{ChiCond}. 

Now, we consider any character $\chi_{\idp} \in \mathcal{X}_{C_p^f}(\OO_{\mathfrak{p}}^{\times})$ which belongs to some fixed chain $\mathfrak{C} \in \N^f$. If this chain satisfies the condition $c_f > m_{\idp}$, the character $\chi'_{\idp} = \chi_{\idp} \cdot \sigma_{\idp}$ belongs to the same chain $\mathfrak{C}$. Hence, we can mimic the computation carried out in equation \eqref{ChiTrick} and its conclusion to obtain $\Sigma(\mathcal{X}_{C_p^f}(\OO_{\mathfrak{p}}^{\times}), \mathfrak{C}; \varepsilon_{\mathfrak{p}}) = 0$. This is part a). 

If we consider a chain $\mathfrak{C} \in \N^h$ for some $1 \leq h \leq f$ satisfying $c_f \leq m_{\idp}$ if $h=f$, then the character $\chi'_{\idp} = \chi_{\idp} \cdot \sigma_{\idp}$ may not belong to the same chain, respectively, the twist by $\sigma_{\idp}$ may cause cancellations or blowups of conductor exponents in this case. However, only conductor exponents $\leq m_{\idp}$ can be modified which yields equation \eqref{ZeroSum}. 
\end{proof}

Our results so far are good enough to show that $\mathfrak{F}_{C_p^{f}, S}(p^{r-f}s; \varepsilon)$ for $1 \leq f \leq r$ has a holomorphic continuation to a sufficiently large domain. 
Note that we do not precisely compute the Euler factors for some given place $\idp$. If we have a place satisfying $\deg(\idp) \leq 2g_F - 2$, we might have smallest conductor exponents for the chains larger than $2$ by Proposition \ref{ChiCond}. In addition, we would have to determine exactly all $\idp$-valuations for all differentials of the components of $\varepsilon$ which is certainly very complicated and involving and can not be performed in general terms since it is highly depending on the arithmetic properties of the base field $F$ under consideration. 

Even if we consider a place satisfying $\deg(\idp) > 2g_F - 2$, the precise computation of the Euler factor is still quite hard due to part b) of Proposition \ref{ChiVanish} Again, it is rather complicated to gain this information and it is highly depending on the properties of $\varepsilon$ and therefore also on the arithmetic properties of $F$.

\begin{proposition} \label{AsympEst}
Let $1 \leq f \leq r$ and $\varepsilon$ be a non-trivial class $\varepsilon = (\varepsilon_1, \ldots, \varepsilon_f) \in \mathcal{E}_{C_p^f}(S)$ such that all components of $\varepsilon$ represent a class in the $p$-Selmer group. 
\begin{enumerate}
	\item[a)] Let $f = 1$. Then, $\mathfrak{F}_{C_p, S}\left(p^{r-1}s; \varepsilon\right)$ has a holomorphic continuation to the domain $\re(s) > \frac{1}{4(p-1)p^{r-1}}$
	\item[b)] Let $f \geq 2$. Then, $\mathfrak{F}_{C_p^{f}, S}\left(p^{r-f}s; \varepsilon\right)$ has a holomorphic continuation to the domain $\re(s) > \Gamma$ with $\Gamma := \max\left\{\frac{1+ (f-1)(p-1)}{p^{r+1-(f-1)} \left(p^{f-1} -1\right)}, \frac{f(p-1)}{p^{r+1-f} \left(p^{f} -1\right)}\right\}$.
	\item[c)] Let $f = r = p= 2$. Then, we have $\Gamma = \frac{1}{2}$ and $\zeta^{-3}_F(4s-1) \cdot\mathfrak{F}_{C_2^{2}, S}(s; \varepsilon)$ has a holomorphic continuation to the domain $\re(s) > \frac{5}{12}$.
\end{enumerate}
\end{proposition}

\begin{proof}
Let $f=1$. By equation \eqref{EulerProd} in Proposition \ref{PropDecom}, Lemma \ref{EulerfactorO}, Proposition \ref{ChiCond} and Corollary \ref{LagResult}, we have
\begin{align*}
& \mathfrak{F}_{C_p, S}\left(p^{r-1}s; \varepsilon\right) = p^{f \vert S \vert} \cdot \prod\limits_{\idp \in \PP_F} \sum\limits_{\chi_{\idp} \in \mathcal{X}_{C_p^{f}}(\OO_{\idp}^{\times})} \chi_{\idp}(\varepsilon_{\idp}) \vert \vert {\partial\left(F_{\chi_{\idp}}/F\right)} \vert \vert^{-p^{r-k(\chi_{\idp})}s} \\ 
	= \quad &p^{f \vert S \vert} \cdot \prod\limits_{\substack{\idp \in \PP_F: \\ \deg(\idp) > 2g_F - 2}} (1 - \np^{-(p-1)p^{r-1}2s}) \cdot \prod\limits_{\substack{\idp \in \PP_F: \\ \deg(\idp) \leq 2g_F - 2}} \sum\limits_{\chi_{\idp} \in \mathcal{X}_{C_p^{f}}(\OO_{\idp}^{\times})} \chi_{\idp}(\varepsilon_{\idp}) \vert \vert {\partial\left(F_{\chi_{\idp}}/F\right)} \vert \vert^{-p^{r-k(\chi_{\idp})}s}.
\end{align*}
The finite product at the back is holomorphic on the whole complex plane by Remark \ref{Rk1}. For the infinite product, we have
\begin{align*} 
	\prod\limits_{\substack{\idp \in \PP_F: \\ \deg(\idp) > 2g_F - 2}} (1 - \np^{-(p-1)p^{r-1}2s}) = \zeta_F((p-1)p^{r-1}2s)^{-1} \prod\limits_{\substack{\idp \in \PP_F: \\ \deg(\idp) \leq 2g_F - 2}} (1 - \np^{-(p-1)p^{r-1}2s})^{-1}.
\end{align*}
All poles of the inverse zeta function are on the axis with $\re(s) = \frac{1}{4(p-1)p^{r-1}}$ by the Hasse-Weil Theorem and all poles of the finite product at the back lie on the axis with $\re(s) = 0$. 

Now, let $2 \leq f$. Our approach is a comparison of the Dirichlet series $\mathfrak{F}_{C_p^{f}, S}\left(p^{r-f}s; \varepsilon\right)$ with the Dirichlet series $\Phi_{f-1}(F, C_p^r;s)$ of the $(f-1)$-dimensional correction of the trivial character which we have already investigated in Section \ref{CFFqtg}. In fact, we want to apply Proposition \ref{ContinuationEulerProduct}. 

We recall that all chains appearing in the series $\Phi_{f-1}(F, C_p^r;s)$ have the form $(c_1, \ldots, c_h) \in \N^h$ for some $0 \leq h \leq f-1$. Initially, the chains appearing in the series $\mathfrak{F}_{C_p^{f}, S}(s; \varepsilon)$ have the form $(c_1, \ldots, c_h) \in \N^h$ for some $0 \leq h \leq f$. 
By Proposition \ref{ChiVanish}, however, the character sums for all chains $\mathfrak{C} \in \N^f$ such that $c_f > m_{\idp}$ vanish. 

%Hence, it is sufficient to consider chains with the shape $\mathfrak{C} = (c_1, \ldots, c_J, c_{J+1}, \ldots, c_h) \in \N^h$ where $0 \leq J < h \leq f$ is the maximal index such that $c_J > m_{\idp}$ if available. If and only if $c_1 \leq m_{\idp}$, then we have $J=0$. \\
By Proposition \ref{ChiCond}, there are only finitely many places $\idp \in \PP_F$ with $m_{\idp} > 2$. We can give an upper bound for the Euler factors of $\mathfrak{F}_{C_p^{f}, S}(s; \varepsilon)$ corresponding to these places by considering the respective Euler factors $\Phi_{f, \idp}(F, C_p^r; s)$ for the trivial class. Hence, this finite product of Euler factors of $\mathfrak{F}_{C_p^{f}, S}(s; \varepsilon)$ has a holomorphic continuation to the domain $\re(s) > \frac{f(p-1)}{p^{r+1-f}\left(p^f-1\right)}$ by Proposition \ref{LocalMeromorphicContinuation}, Corollary \ref{LocalPoles} and Lemma \ref{LocalInequalityAbscissa}. 

Now, we consider the infinite Euler product of places $\idp \in \PP_F$ satisfying $\deg(\idp) > 2g_F-2$ and hence $m_{\idp} = 2$ by Proposition \ref{ChiCond}. We can give an upper estimate for its Euler factors with respect to equation \eqref{ZeroSum} in the following way: 
Let $\mathfrak{C} = (c_1, \ldots, c_J) \in \N^J$ be a chain with some $0 \leq J < f$ and $c_J > 2$. Now, we consider the set of all chains $\mathfrak{C'} \in \N^h$ for some $J \leq h \leq f$ such that $\mathfrak{C}'_{>2} = \mathfrak{C}$. Then, we have exactly one chain for each index $J \leq h \leq f$, namely $\mathfrak{C}_h := (c_1, \ldots, c_J, 2, \ldots 2) \in \N^h$ whose last $h-J$ components are $2$, and equation \eqref{ZeroSum} applies. For $h < f$, we naively estimate
\begin{align*}
\vert {\Sigma(\mathcal{X}_{C_p^f}(\OO_{\mathfrak{p}}^{\times}), \mathfrak{C}_h; \varepsilon_{\mathfrak{p}})} \vert = \vert {\sum\limits_{\substack{\chi_{\idp} \in \mathcal{X}_{C_p^f}(\OO_{\mathfrak{p}}^{\times}): \\ \mathfrak{C}(\chi_{\idp}) = \mathfrak{C}_h}} \chi_{\idp}(\varepsilon_{\idp})} \vert \leq \vert \{\chi_{\idp} \in \mathcal{X}_{C_p^f}(\OO_{\mathfrak{p}}^{\times}) \ : \ \mathfrak{C}(\chi_{\idp}) = \mathfrak{C}_h\} \vert
\end{align*}
which is precisely the contribution of the chain $\mathfrak{C}_h$ in the case of the trivial class $\varepsilon = \mathds{1}_f$. For the chain satisfying $h=f$, we have the estimate
\begin{align} \label{est}
	&\vert {\Sigma(\mathcal{X}_{C_p^f}(\OO_{\mathfrak{p}}^{\times}), \mathfrak{C}_f; \varepsilon_{\mathfrak{p}})} \vert \notag \\
	\leq \quad &\sum\limits_{h=J}^{f-1} \vert \{\chi_{\idp} \in \mathcal{X}_{C_p^f}(\OO_{\mathfrak{p}}^{\times}) \ : \ \mathfrak{C}(\chi_{\idp}) = \mathfrak{C}_h\} \vert \notag \\
	\leq \quad & (f-J) \cdot \vert \{\chi_{\idp} \in \mathcal{X}_{C_p^f}(\OO_{\mathfrak{p}}^{\times}) \ : \ \mathfrak{C}(\chi_{\idp}) = \mathfrak{C}_{f-1} \} \vert
\end{align}
by equation \eqref{ZeroSum} and Definition \ref{DefinitionCountingFunction}. Note that the chain $\mathfrak{C}_f$ yields the largest discriminants of all chains $\mathfrak{C}_h$ for $J \leq h \leq f$ by Lemma \ref{DiscExponent}. 

The Euler factor which is constructed by these estimates looks very similar to $\Phi_{f-1, \idp}(F, C_p^r;s)$. In fact, they coincide except for the chains with the property $h  = f$ which do not appear in $\Phi_{f-1, \idp}(F, C_p^r;s)$. This difference, however, does not cause a shift of the convergence abscissa further to the right since the character sum of these chains can be bounded by an integral multiple of the character sum of a shorter chain which has a smaller discriminant by equation \eqref{est}. Consequently, $\mathfrak{F}_{C_p^{f}, S}(s; \varepsilon)$ is holomorphic in the domain $\re(s) > \frac{1+ (f-1)(p-1)}{p^{r+1-(f-1)} \left(p^{f-1} -1\right)}$ by Proposition \ref{GlobalPoles} and Proposition \ref{ContinuationEulerProduct}. 

For $f = r = p = 2$, we take a closer look at the series $\mathfrak{F}_{C_2^{2}, S}(s; \varepsilon)$. By Proposition \ref{ChiVanish}, we only have to deal with the empty chain $()$ and chains of the form $(c_1)$ and $(c_1,2)$ for any $c_1 \geq 2$. In particular, equation \eqref{ZeroSum} yields 
\begin{align*}
\Sigma(\mathcal{X}_{C_2^2}(\OO_{\mathfrak{p}}^{\times}), (c_1); \varepsilon_{\mathfrak{p}}) = - {\Sigma(\mathcal{X}_{C_2^2}(\OO_{\mathfrak{p}}^{\times}), (c_1,2); \varepsilon_{\mathfrak{p}})}
\end{align*}
for all $c_1 > 2$. By Proposition \ref{ChiCond}, we also have 
\begin{align*}
1 + \Sigma(\mathcal{X}_{C_2^2}(\OO_{\mathfrak{p}}^{\times}), (2); \varepsilon_{\mathfrak{p}}) = - {\Sigma(\mathcal{X}_{C_2^2}(\OO_{\mathfrak{p}}^{\times}), (2,2); \varepsilon_{\mathfrak{p}})},
\end{align*}
since only the chains $(), (2)$ and $(2,2)$ have a conductor that divides $\idp^2$ and the chain $()$ only consists of the trivial character.
Therefore, our computations for the series $\Phi_2(F, C_2^2;s)$ in Proposition \ref{ContinuationEulerProduct} yield that $\zeta^{-3}_F(4s-1) \cdot\mathfrak{F}_{C_2^{2}, S}(s; \varepsilon)$ is holomorphic in the domain $\re(s) > \frac{5}{12}$.
\end{proof}

\begin{theorem} \label{ContinuationF}
		The complex function 
	\begin{align*}
		\Theta_{F,r,p}(s) := \Phi(F, C_p^r; s) \cdot \Lambda_{F,r,r,p}^{-1}(s)
	\end{align*}
	has a holomorphic continuation to the domain $\re(s) > \frac{1+r(p-1)}{p \left(p^r-1\right)} - \epsilon$ with $\epsilon = \frac{1}{p^{2} (p^r -1)}$. 
	
	%Here, $\Lambda_{F,r,r,p}(s)$ is as in Definition \ref{MeromorphicFunctionGlobal} except for the difference that we use the zeta function $\zeta_F(s)$ of the base field $F$ instead of the rational function field $\Fqtg$.\\
	Furthermore, $\Theta_{F,r,p}(s)$ does not vanish at $s = \frac{1+r(p-1)}{p\left(p^r-1\right)}$, so $\Phi(F, C_p^r;s)$ has a pole at this point. The order of this pole is given by the order of vanishing of $\Lambda_{F,r,r, p}^{-1}(s)$ at this value which is $b_p(C_p^r)$ defined in Theorem \ref{MainResult}.

%	The Dirichlet series $\Phi(F, C_p^r; s)$ has the meromorphic continuation
%	\begin{align*}
%		\Lambda_{F,r,r,p}(s) \cdot \Psi_{F,r,p}(s).
%	\end{align*}
%	The continuation applies in the domain $\re(s) > \frac{1+r(p-1)}{p \left(p^r-1\right)} - \epsilon$
%	for every $\epsilon$ satisfying $0 < \epsilon \leq \frac{1}{p^{2} (p^r -1)}$ and the function $\Psi_{F,r,p}(s) := \Phi(F, C_p^r; s) \cdot \Lambda^{-1}_{F,r,r,p}(s)$ is holomorphic in this domain. Here, $\Lambda_{F,r,r,p}(s)$ is as in Definition
%	\ref{MeromorphicFunctionGlobal} except for the difference that we use the zeta function $\zeta_F(s)$ of the base field $F$ instead of the rational function field $\Fqtg$.
\end{theorem}

\begin{proof}
	The proof is analogous to that of Theorem \ref{ContinuationDirichletSeries}. %We, however, must take into account that the zeta functions $\zeta_F(s)$ of an arbitrary base field might have zeros on the line with $\re(s) = \frac{1}{2}$ by \cite[Theorem 5.10]{Ro}. Therefore, all possible zeros of $\Lambda_{F,r,r,p}(s)$ satisfy $\re(s) = \frac{\frac{1}{2} + r(p-1)}{p \left(p^r-1\right)} \leq \frac{1 - \frac{1}{p} + r(p-1)}{p \left(p^r-1\right)}$ and $\Lambda^{-1}_{F,r,r,p}(s)$ does not have poles in the considered domain. \\
	%It remains to check the holomorphy of $\Psi_{F,r,p}(s)$. We have 
	We have
	\begin{align*}
		&\Theta_{F,r,p}(s) = \Phi(F, C_p^r; s) \cdot \Lambda^{-1}_{F,r,r,p}(s) \\
		= \quad & \left(\sum\limits_{f=0}^{r} e_f \cdot \Phi_f(F, C_p^r;s) + \sum\limits_{f=0}^r \frac{e_f}{p^f \vert {\mathcal{E}_{C_p^f}(S)} \vert} \cdot \sum\limits_{\substack{\varepsilon \in \mathcal{E}_{C_p^f}(S) \\ \varepsilon \neq \mathds{1}_f}} \mathfrak{F}_{C_p^{f}, S}\left(p^{r-f}s; \varepsilon\right) \right) \cdot \Lambda^{-1}_{F,r,r,p}(s).
	\end{align*}
	by Proposition \ref{WrightDecom} and Proposition \ref{PropDecom}.
	For $\left(\sum\limits_{f=0}^{r} e_f \cdot \Phi_f(F, C_p^r;s) \right) \cdot \Lambda^{-1}_{F,r,r,p}(s)$, the holomorphy in the considered domain and the non-vanishing at the point $s = \frac{1+r(p-1)}{p\left(p^r-1\right)}$ is ensured by Theorem \ref{ContinuationDirichletSeries}. Proposition \ref{LambdaZeros}, Proposition \ref{AsympEst} and Lemma \ref{GlobalInequality} yield the holomorphy of the remaining part. 
	
	By Proposition \ref{GlobalPoles}, the remaining part vanishes at $s = \frac{1+r(p-1)}{p\left(p^r-1\right)}$ in all cases except $r=p=2$ since $\Lambda^{-1}_{F,r,r,p}(s)$ has a zero and $\mathfrak{F}_{C_p^{f}, S}\left(p^{r-f}s; \varepsilon\right)$ is holomorphic at this point. 
	
	In the case $r=p=2$, we have to be more precise. We, however, still have a zero at this point since $\Lambda^{-1}_{F,2,2,2}(s)$ has a zero of order $4$ and $\mathfrak{F}_{C_2^{2}, S}\left(s; \varepsilon\right)$ has a pole of order at most $3$.
\end{proof}

%\begin{theorem} \label{AsymptoticsF}
%	We recall the notation provided in Theorem \ref{GlobalAsymptotics}. Then, there exists a real polynomial $P_{m \bmod L(C_p^r), F}$ of degree at most $b_p(C_p^r) - 1$ depending on the residue class of $m$ modulo $L(C_p^r)$ and the base field $F$ such that the number of $C_p^r$-extensions of a global function field $F$ when counted by discriminants satisfies
%	\begin{align} \label{GA2}
%		Z(F, C_p^r; X) = P_{m \bmod L(C_p^r), F}(\log(X)) \cdot X^{a_p(C_p^r)} + O\left(X^{\frac{1 - \frac{1}{p} + r(p-1)}{p(p^{r} -1)} + \epsilon}\right)
%	\end{align}
%	for every $\epsilon > 0$. The polynomial $P_{m \bmod L(C_p^r), F}$ has degree $b_p(C_p^r) - 1$ for at least one residue class and its leading coefficient is always positive if $P_{m \bmod L(C_p^r), F}$ has this degree. 
%\end{theorem}

\begin{proof}[Proof of Theorem \ref{MainResult} for all global function fields $F$ as base field] \hfill \\
	The proof is completely analogous to that for the special case $F = \Fqtg$ and follows by Theorem \ref{ContinuationF}.
\end{proof}

\begin{remark} \label{C-Constant}
	We denote by $c(F, C_p^r, m \bmod L_{\text{max}}(C_p^r))$ the leading coefficient of the polynomial $P_{m \bmod L(C_p^r), F}$ (that is the coefficient corresponding to the $b_p(C_p^r)-1$-th power) and note that we allow this coefficient to be $0$ and that it only depends on the residue class of $m \bmod L_{\text{max}}(C_p^r)$ where $L_{\text{max}}(C_p^r)$ is the period of the poles of order exactly $b(C_p^r)$ and is given by $L_{\text{max}}(C_p^r) := \begin{cases}
		p-1 & \text{for } r=1 \ \text{and} \ p \neq 2, \\
		2 & \text{for } r \in \{1, 2\} \ \text{and} \ p=2, \\
		p(p^r - 1) & \text{else. }
	\end{cases}$ 

	As in the case of local function fields (Theorem \ref{LocalAsymptotics} and Remark \ref{MTLR}), a precise computation of the constant $c(F, C_p^r, m \bmod L_{\text{max}}(C_p^r))$ is very tedious as we have to compute the residue at each possible pole $R \zeta^{-j}$ of order exactly $b(C_p^r)$ where $\zeta$ is a primitive $L_{\text{max}}(C_p^r)$-th root of unity and we again have to carry out a complicated summation over these residues and the $L_{\text{max}}(C_p^r)$-th roots of unity. At least for the groups $C_p$ and $C_2^2$, where we have $L_{\text{max}}(C_p^r) \in \{2,p-1\}$, the effort is still very limited and a precise computation can be achieved quite quickly. 
	
	In the case of $C_p$-extensions, the constant is calculated explicitly by Lagemann in \cite[Addendum 1.3, Addendum 1.4 and Theorem 2.1]{La2}. Note that we have to remove the factor $\frac{1}{1-R^{L_{\text{max}}(C_p^r)}}$ in Lagemann's constants as he considered the asymptotics of the partial sums of the coefficients instead of the asymptotics of the coefficients themselves. In fact, Lagemann's counting function is defined slightly differently and he counts the number of extensions whose norm of the conductor, respectively discriminant, is at most $X$ instead of equal to $X$. In the case $p\neq 2$, it is also worth noting that we only have $c(F, C_p, m \bmod p-1) \neq 0$ for the progression $0 \bmod p-1$, which is to be expected since all local discriminants are congruent to $0 \bmod p-1$.  
	
	For the group $C_2^2$, one can compute $c(F, C_2^2, m \bmod 2) =$
	\begin{align*}
		\frac{\vert C_2^2 \vert}{\vert {\Aut(C_2^2)} \vert} \cdot \frac{\log(q)}{144} \cdot \zeta_{F}(1)^4 \cdot \gamma(m \bmod 2) \cdot \prod\limits_{\idp \in \PP_{F}} \left(\left(1 + 4 \cdot \np^{-1} + \np^{-2} \right) \cdot \left(1 - \np^{-1} \right)^4 \right)
	\end{align*}
	where $\zeta_{F}(1)$ denotes the value of the residue of $\zeta_{F}(s)$ at the simple pole at $s = 1$ and we have the constant $\gamma(m \bmod 2) := \begin{cases}
		1 & \text{for } m \equiv 0 \mod 2, \\
		0 & \text{for } m \equiv 1 \mod 2.
	\end{cases}$ 

	In fact, $c(F, C_2, 1 \bmod 2) = 0$ as well as $c(F, C_2^2, 1 \bmod 2) = 0$ should be expected as the norm of a local conductor of a $C_p$-extension may not be congruent to $1 \bmod p$. Hence, in this case, each norm of a local conductor is congruent to $0 \bmod 2$ and thus the norms of global conductors or global discriminants are always congruent to $0 \bmod 2$. 
	
	For the remaining groups with $L_{\text{max}}(C_p^r) = p(p^r-1)$, one obtains
	\begin{align*}
	c(F, C_p^r, m \bmod L_{\text{max}}(C_p^r)) = \frac{e_r}{L_{\text{max}}(C_p^r)} \cdot \log(q) \cdot \zeta_F(1) \cdot \sum\limits_{j=1}^{L_{\text{max}}(C_p^r)} \zeta^{jm} \cdot \psi_{F,r,r,p}(R\zeta^{-j})
	\end{align*}
where we recall the definition of $\psi_{F,r,r,p}(s)$ from Proposition \ref{ContinuationEulerProduct} and we apply the substitution $z = q^{-s}$.
\end{remark}

\begin{remark}
	Finally, we want to discuss the relation to Lagemann's results for the conductor density in \cite[Theorem 1.1 and Theorem 1.2]{La2} and we recall Remark \ref{LocalConfig}. 
	
	For $r=1$, the count by conductors and by discriminants are essentially identical, except for the weight $p-1$ due to the conductor-discriminant formula. Hence, we focus on the case $r \geq 2$.
	In fact, in the case $r \geq 2$, Lagemann proved that the main term of the asymptotics of the count by conductor originates from the local configurations where the local conductor is precisely $p$. A closer look at his computations reveals that the leading term actually comes from the configuration where all the conductors of $C_p$-subextensions are equal to $p$. Consequently, we can interpret the factor $p^r-1$ occurring in the denominator of $a_p(C_p^r)$ as a weight coming from the local configuration of $p^r-1$ equal conductors of the $C_p$-subextensions. Both the main term of the asymptotics of the count by conductors as well as the main term of the asymptotics of the count by discriminants can be traced back to this local configuration where we additionally have to consider the contribution of the local $C_2$-algebras with conductor $2$ for the asymptotics of $C_2^2$-extensions counted by discriminant.
\end{remark}

%\newpage

%\section*{Acknowledgments}
%My special thanks go to Jürgen Klüners for his excellent supervision and helpful discussions about this work. I would also like to thank Fabian Gundlach for proofreading an earlier version and for very productive discussions on Wright's paper and the last Section of this article. \\
%The work was funded by the Deutsche Forschungsgemeinschaft (DFG, German Research Foundation) — SFB-TRR 358/1 2023 — 491392403.

\bibliographystyle{plain}
%\bibliography{mybib}

%\bibliographystyle{elsarticle-num}
%\section*{\refname}
%\bibliography{mybib} 

\begin{thebibliography}{CCC}	% Zitiertes und auch Nützliches/ sehr lesenswertes!!
%\bibitem[Bu]{Bu}
%L.~Butler. \emph{Subgroup Lattices and Symmetric Functions}, Memoirs of AMS~Number~539, 1994.






\bibitem[Bha04]{Bha1}
M.~Bhargava. \emph{Higher composition laws III: The parametrization of quartic rings}, Annals of Mathematics \textbf{159} (2004) 1329-1360


\bibitem[Bha08]{Bha2}
M.~Bhargava. \emph{Higher composition laws IV: The parametrization of quintic rings}, Annals of Mathematics \textbf{167} (2008) 53-94

\bibitem[CSV24]{CSV}
A.~Carnevale, M.~Schein, C.~Voll. \emph{Generalized Igusa functions and ideal growth in nilpotent Lie rings}, 
% 	\hyperref{https://doi.org/10.1016/j.jnt.2014.09.026}, 
Algebra \& Number Theory \textbf{18:3} (2024) 537-582



\bibitem[DH71]{DH71}
H.~Davenport, H.~Heilbronn. \emph{On the Density of Discriminants of Cubic Fields. II}, Proceedings of the Royal Society of London. Series A \textbf{322} (1971), No. 1551, 405-420

\bibitem[Del48]{De}
S.~Delsarte. \emph{Fonctions de M\"obius sur les groupes abeliens finis}, Annals of Mathematics \textbf{49}(3) (1948) 600-609
% 	, Brian~Mortimer; \emph{Permutation Groups}, Graduate Texts in Mathematics, Springer, 1996.
%\bibitem[DiMo]{DM}
%J.~Dixon, B.~Mortimer. \emph{Permutation Groups}, Graduate Texts in Mathematics, Springer, 1996.
%\bibitem[ElVe]{EV}
%J.S.~Ellenberg, A.~Venkatesh. \emph{Counting extensions of function fields with specified Galois group
%	and bounded discriminant}, Geometric Methods in Algebra and Number Theory, 235; pp. 151-168.
%\bibitem[FeVo]{FV}
%I.~B.~Fesenko, S.~Vostokov. \emph{Local fields and their extensions}, AMS, 1993.
%\bibitem[Ha]{Ha}
%H.~Hasse. \emph{Zahlentheorie}, Akademie Verlag, Berlin, 1969.
%\bibitem[Iw]{Iw} K.~Iwasawa. \emph{Local Class Field Theory}, Oxford University Press, 1986.
%\bibitem[Ja]{Ja} U.~Jannsen, K.~Wingberg,
%\emph{Einbettungsprobleme und Galoisstruktur lokaler K\"orper},
%J. Reine Angew. Math. \textbf{319} (1980), 196--212. 
%\bibitem[Le]{Le}
%A.~Ledet. \emph{Brauer Type Embedding Problems}, AMS, 2005.


\bibitem[ETW23]{ETW}
J.~S.~Ellenberg, T.~Tran, C.~Westerland. \emph{Fox-Neuwirth-Fuks cells, quantum shuffle algebras,
	and Malle’s conjecture for function fields}, arXiv (2023)



\bibitem[EV05]{EV}
J.~S.~Ellenberg, A.~Venkatesh. \emph{Counting extensions of function fields with bounded discriminant and specified Galois group}, Geometric Methods in Algebra and Number Theory, Birkhäuser (2005) 151-168


\bibitem[FS09]{FS} P.~Flajolet, R.~Sedgewick. \emph{Analytic Combinatorics}, Cambridge University Press (2009)




\bibitem[KM20]{KM}
J.~Klüners, R.~Müller. \emph{The conductor density of local function fields with Abelian Galois group}, 
% 	\hyperref{https://doi.org/10.1016/j.jnt.2014.09.026}, 
Journal of Number Theory \textbf{212} (2022) 311-322




\bibitem[Lag10]{La1}
T.~Lagemann. \emph{Asymptotik wild verzweigter Funktionenk\"orper}, 
% 	\hyperref{http://dx.doi.org/10.14279/depositonce-2654}, 
Logos Verlag Berlin, 2010
\bibitem[Lag12]{La2}
T.~Lagemann. \emph{Distribution of Artin-Schreier extensions}, 
% 	\hyperref{https://doi.org/10.1016/j.jnt.2014.09.026}, 
Journal of Number Theory \textbf{132} (2012) 1867-1887


\bibitem[Lag15]{La3}
T.~Lagemann. \emph{Distribution of Artin-Schreier-Witt extensions}, 
% 	\hyperref{https://doi.org/10.1016/j.jnt.2014.09.026}, 
Journal of Number Theory \textbf{148} (2015) 288-310





\bibitem[Lan00]{Lang} S.~Lang. \emph{Algebra}, Graduate Texts in Mathematics \textbf{211}, Springer (2000)




\bibitem[Mal02]{Malle1}
G.~Malle. \emph{On the distribution of Galois groups}, Journal of Number Theory \textbf{92} (2002) 315-329



\bibitem[Mal04]{Malle2}
G.~Malle. \emph{On the distribution of Galois groups, II}, Experimental Mathematics \textbf{13}(2) (2004) 129-135



% \bibitem[KlMa]{Ma} 
% J.~Kl\"uners, G.~Malle. \emph{Counting Nilpotent Galois Extensions}, J. reine angew. Math., 572, 2004, 1--26.
%\bibitem[M\"a]{Mae} S.~M\"aki. \emph{The conductor density of abelian number fields}, J. London Math. Soc. (2) 47, no. 1, pp. 18--30, 1993.


\bibitem[Neu99]{Neu} J.~Neukirch. \emph{Algebraic Number Theory}, Grundlehren der mathematischen Wissenschaften \textbf{322}, Springer (1999)

\bibitem[Ros02]{Ro} M.~Rosen. \emph{Number Theory in Function Fields}, Graduate Texts in Mathematics \textbf{210}, Springer (2002)



%\bibitem[Se1]{Se} 
%J.-P.~Serre. \emph{Une ``formule de masse''  pour les extensions totalement ramifi\'ees de degr\'e donn\'e d'un corps local}, C.R.~Acad.~Sc.~Paris~S\'erie A \textbf{286} (1978), 1031-1036.
%\bibitem[Se2]{Se2} 
%J.-P.~Serre. \emph{Linear Representations of Finite Groups}, Springer-Verlag, 1977.
%
%\bibitem[Sp]{Sp}
%A.~Speiser. \emph{Theorie der Gruppen von endlicher Ordnung}, Springer-Verlag, 1937.

\bibitem[Sti09]{St} H.~Stichtenoth. \emph{Algebraic Function Fields and Codes}, Graduate Texts in Mathematics \textbf{254}, Springer (2009)


\bibitem[Woo10]{Woo} M.~M.~Wood. \emph{On the probabilities of local behaviors in abelian field extensions}, Compositio Mathematica \textbf{146} (2010) 102-128


\bibitem[Wri89]{Wr} D.~J.~Wright. \emph{Distribution of discriminants of abelian extensions}, Proceedings of the London Mathematical Society (3) \textbf{58}(1) (1989) 17-50


\end{thebibliography}

\end{document}